\documentclass[10pt]{amsart}

\usepackage[utf8]{inputenc}
\usepackage[T1]{fontenc}
\usepackage{textcomp}
\usepackage{bm}
\usepackage{amsfonts}
\usepackage{mathtools}
\usepackage{mathrsfs}   % script capitals

\usepackage[square,compress,comma,numbers,sort]{natbib}
\usepackage{color}
\usepackage{float}
\usepackage{comment}

\usepackage[shortlabels]{enumitem}

\usepackage[top=2cm,bottom=2cm,left=2.5cm,right=2.5cm,marginparwidth=2cm]{geometry}

\allowdisplaybreaks[4]

\usepackage[unicode]{hyperref}
\hypersetup{
  colorlinks=true,
  citecolor=red,
  linkcolor=blue,
  pdftitle={Sojourns of Vector-Valued Stationary Gaussian Random Fields},
  pdfauthor={Krzysztof Dębicki; Enkelejd Hashorva; Zbigniew Michna},
  pdfsubject={Exact asymptotics for sojourn times of vector-valued stationary Gaussian random fields},
  pdfkeywords={exact asymptotics; high exceedance probability; sojourn times; occupation times; vector-valued Gaussian random fields; stationary Gaussian random fields; Berman constants; Pickands condition; varying domains},
}
\usepackage[nameinlink]{cleveref}

\usepackage{aliascnt}
\theoremstyle{plain}
\newtheorem{theorem}{Theorem}[section]

\newaliascnt{lemma}{theorem}
\newtheorem{lemma}[lemma]{Lemma}
\aliascntresetthe{lemma}

\newaliascnt{proposition}{theorem}
\newtheorem{proposition}[proposition]{Proposition}
\aliascntresetthe{proposition}

\newaliascnt{corollary}{theorem}
\newtheorem{corollary}[corollary]{Corollary}
\aliascntresetthe{corollary}

\theoremstyle{definition}
\newaliascnt{definition}{theorem}
\newtheorem{definition}[definition]{Definition}
\aliascntresetthe{definition}

\newaliascnt{example}{theorem}
\newtheorem{example}[example]{Example}
\aliascntresetthe{example}

\theoremstyle{remark}
\newaliascnt{remark}{theorem}
\newtheorem{remark}[remark]{Remark}
\aliascntresetthe{remark}

\crefname{theorem}{Theorem}{Theorems}
\Crefname{theorem}{Theorem}{Theorems}
\crefname{lemma}{Lemma}{Lemmas}
\Crefname{lemma}{Lemma}{Lemmas}
\crefname{proposition}{Proposition}{Propositions}
\Crefname{proposition}{Proposition}{Propositions}
\crefname{corollary}{Corollary}{Corollaries}
\Crefname{corollary}{Corollary}{Corollaries}
\crefname{definition}{Definition}{Definitions}
\Crefname{definition}{Definition}{Definitions}
\crefname{example}{Example}{Examples}
\Crefname{example}{Example}{Examples}
\crefname{remark}{Remark}{Remarks}
\Crefname{remark}{Remark}{Remarks}

\newcommand{\abs}[1]{\left\lvert #1 \right\rvert}

 \DeclarePairedDelimiterXPP\pk[1]{\mathbb{P}}\{ \}{}{ #1}
 \DeclarePairedDelimiterXPP\E[1]{\mathbb{E}}\{ \}{}{	#1}

\NewDocumentCommand{\ceil}{s O{} +m}{
	\IfBooleanTF{#1} % starred
	{$\left\lceil#3\right\rceil$} % \ceil*[..]{..}
	{#2\lceil#3#2\rceil} % \ceil[..]{..}
}
\NewDocumentCommand{\floor}{s O{} m}{
	\IfBooleanTF{#1} % starred
	{$\left\lfloor#3\right\rfloor$}
	{#2\lfloor#3#2\rfloor}
}

\DeclarePairedDelimiter{\oldnormaux}{\bracevert  }{\bracevert }

\NewDocumentCommand{\oldnorm}{som}{
	\IfBooleanTF{\mkern-15mu\vphantom#1 \mkern-15mu\vphantom}
	{\oldnormaux*{\mkern-15mu\vphantom #3}\mkern-15mu\vphantom}
	{\IfNoValueTF{\mkern-15mu\vphantom#2\mkern-15mu\vphantom}
		{\oldnormaux*{\mkern-15mu\vphantom{dq}#3}\mkern-15mu}
		{\oldnormaux[#2]{#3}}%
	}%
}

\definecolor{c20}{rgb}{0.,0.7,0.}
\definecolor{c30}{rgb}{0.,0.,1.}
\definecolor{c40}{rgb}{1,0.1,0.7}
\definecolor{c50}{rgb}{1,0,0}
\definecolor{c60}{rgb}{1,0.9,0.1}
\definecolor{c70}{rgb}{0.50,1.00,0.00}

\def\zm#1{{\textcolor{red}{#1}}}
\def\qq#1{{\textcolor{cyan}{#1}}}
\def\qq#1{#1}

\def\zm#1{#1}

\newcommand{\prooflem}[1]{{Proof of} Lemma \ref{#1}:}

\newcommand{\QED}{\hfill $\Box$}

\newcommand{\COM}[1]{}

\def\IF{\infty}

\newcommand{\R}{\mathbb{R}}

\newcommand{\N}{\mathbb{N}}

\newcommand{\BQN}{\begin{eqnarray}}
\newcommand{\EQN}{\end{eqnarray}}
\newcommand{\BQNY}{\begin{eqnarray*}}
	\newcommand{\EQNY}{\end{eqnarray*}}

\newcommand{\limit}[1]{\lim_{#1 \to   \infty}}

\newcommand{\kb}[1]{\boldsymbol{#1}}
\newcommand{\vk}[1]{\kb{#1}}

\def\bqn#1{\begin{eqnarray} #1 \end{eqnarray}}

\newcommand{\BS}{\begin{proposition}}
	\newcommand{\ES}{\end{proposition}}
\newcommand{\BT}{\begin{theorem}}
	\newcommand{\ET}{\end{theorem}}
\newcommand{\BK}{\begin{corollary}}
	\newcommand{\EK}{\end{corollary}}

\newcommand{\BEX}{\begin{example}}
	\newcommand{\EEX}{\end{example}}

\newcommand{\BD}{\begin{definition}}
	\newcommand{\ED}{\end{definition}}

\newcommand{\BRM}{\begin{remark}}
	\newcommand{\ERM}{\end{remark}}

\newcommand{\BEL}{\begin{lemma}}
	\newcommand{\EEL}{\end{lemma}}

\newcommand\ind[1]{\mathbb{I}{( #1)}}

\definecolor{c20}{rgb}{0.,0.7,0.}
\definecolor{c30}{rgb}{0.,0.,1.}
\definecolor{c40}{rgb}{1,0.1,0.7}
\definecolor{c50}{rgb}{1,0,0}
\definecolor{c60}{rgb}{1,0.9,0.1}

\def\qq#1{\textcolor{cyan}{#1}}
\def\qq#1{#1}

\DeclareMathOperator{\Var}{Var}
\DeclareMathOperator{\Cov}{Cov}

\newcommand{\normE}[1]{\lVert #1 \rVert}

\newcommand{\normF}[1]{\left\Vert#1\right\Vert }

\DeclareMathOperator*{\sgn}{sgn}

\DeclareMathOperator{\diag}{diag}

\def\id{\mathbf{1}}
\def\Z{\vk{Z}}

\def\IF{\infty}

\def\SAi#1#2{S_{i} ( #1, #2) }
\def\RR{\mathcal{R}}

\begin{document}
\title{Sojourns of Vector-Valued Stationary Gaussian Random Fields}

\author{Krzysztof D{\k{e}}bicki}
\address{Krzysztof D{\k{e}}bicki, Mathematical Institute, University of Wroc\l aw, pl. Grunwaldzki 2/4, 50-384 Wroc\l aw, Poland}
\email{Krzysztof.Debicki@math.uni.wroc.pl}

\author{Enkelejd Hashorva}
\address{Enkelejd Hashorva, Department of Actuarial Science, University of Lausanne\\
	Chamberonne 1015 Lausanne, Switzerland}
\email{Enkelejd.Hashorva@unil.ch}

\author{Zbigniew Michna}
\address{Zbigniew Michna, Department of Operations Research and Business Intelligence, Wrocław University of Science and Technology, Wybrzeże Stanisława Wyspiańskiego 27, 50-370 Wrocław, Poland}
\email{Zbigniew.Michna@pwr.edu.pl}
%\thanks{\textsuperscript{*}Corresponding author.}

\date{\today} 

\begin{abstract}
Let $\vk X(\vk t)$, $\vk t\in \R^k$, be a centered stationary
$\R^d$-valued Gaussian random field with covariance matrix function
$
R(\vk s,\vk t)
=
\E{\vk X(\vk s)\vk X(\vk t)^\top}.$
We study the exact asymptotics of
\[
\kappa_u(x)
=
\pk*{
\eta(u)\int_{[0,T]^k}
\ind{\vk X(\vk t)>u\vk b}\,d\vk t>x
},
\qquad
\vk b=(b_1,\ldots,b_d)^\top,
\]
as $u\to\infty$ for $x\ge0$ and $T>0$. Here $\eta(u)$ is the natural
scaling function associated with the local covariance expansion.

We extend Berman's approach and the uniform double-sum method to this
vector-valued random-field setting. Our formulation uses weaker assumptions
and permits varying  domains $T=T_u$ at an admissible rate. It also yields
exact tail asymptotics for the supremum functional of $\vk X$.
\end{abstract}

\bigskip

\maketitle


\bigskip


\section{Introduction}

Let
$
\vk X(\vk t)=\bigl(X_1(\vk t),\ldots,X_d(\vk t)\bigr)^\top,
  \vk t\in\R^k
$
be a jointly measurable centered stationary $\R^d$-valued Gaussian random
field (GRF) on a complete probability space. Its covariance matrix function is
\[
R(\vk s,\vk t)
=
\E{\vk X(\vk s)\vk X(\vk t)^\top},
 \qquad \vk s,\vk t\in\R^k.
\]
Given a deterministic threshold
\[
u\vk b=(ub_1,\ldots,ub_d)^\top,
\qquad
\vk b=(b_1,\ldots,b_d)^\top\in\R^d\setminus(-\infty,0]^d,
\qquad u>0,
\]
we define the corresponding sojourn functional of $\vk X$ on $[0,T]^k$
above $u\vk b$ by
\begin{equation}
L_u(T)
=
\int_{[0,T]^k}
\ind{\vk X(\vk t)>u\vk b}\,d\vk t,
\qquad T>0.
\label{loc.time}
\end{equation}

Sojourn functionals arise in insurance risk modelling
\cite{OKP3,novikov2025sojourns}, queueing theory
\cite{BorstNunezZwart2006,debicki2023sojourns}, financial mathematics
\cite{AKA,OKP4,GemanYor1996} and statistics
\cite{pham2013rate,leonenko2023sojourn,leonenko2025high}.

Motivated by Berman's seminal results \cite{Berman82,Berman92} we choose a
positive scaling function $\eta(u)$ and study as $u\to\infty$ the probabilities
\[
\kappa_u(x)
=
\pk{\eta(u)L_u(T)>x},
\qquad x\ge0.
\]

As in the one-dimensional setting
\cite{Berman82,Berman92,DLM20,DHLM23,debicki2023sojourns} the asymptotics of
$\kappa_u(x)$ for $x\ge0$ are essentially controlled by those of
$\kappa_u(0)$. For continuous sample paths,
\[
\kappa_u(0) 
=
\pk{\exists_{\vk t\in[0,T]^k}\,\vk X(\vk t)>u\vk b}.
\]

For intuition consider $d=k=1$ and fix $T>0$. Let $X(t)$, $t\in\R$, be a
centered stationary Gaussian process with covariance function $\RR(t)$,
unit variance and continuous sample paths. Assume that $1-\RR(t)>0$ for
$t>0$ and is regularly varying at $0$ with index $\alpha\in(0,2]$. By
\cite[Thm 3.3.1]{Berman92}
\begin{equation}
\lim_{u\to\infty}
\frac{\int_0^x y\,dF_u(y)}
     {\eta(u)\E{L_u(T)}}
=
\mathfrak F(x)
\label{eq:A1}
\end{equation}
for all continuity points $x>0$ of  
$\mathfrak F$, where $\eta(u)$ is a positive scaling function satisfying
\[
1-\RR(1/\eta(u))\sim u^{-2},
\qquad u\to\infty,
\]
and $F_u$ denotes the distribution function of $\eta(u)L_u(T)$. Throughout
$g(u)\sim h(u)$ as $u\to\infty$ means $\limit{u}g(u)/h(u)=1$.

The convergence in \eqref{eq:A1} yields the tail asymptotics of
$\eta(u)L_u(T)$. Indeed \cite[Lem 1.2.1]{Berman92} gives for almost all
$x>0$
\begin{equation}
\pk{\eta(u)L_u(T)>x}
\sim
T\left(
\int_x^\infty \frac1y\,d\mathfrak F(y)
\right)
\eta(u)\pk{X(0)>u},
\qquad u\to\infty .
\label{eq:A12}
\end{equation}
The extension of \eqref{eq:A12} to $x=0$ is not straightforward. It was
established in \cite[Thm 10.5.1]{Berman92} that
\begin{equation}
\pk*{\sup_{t\in[0,T]}X(t)>u}
\sim
T H\eta(u)\pk{X(0)>u},
\qquad u\to\infty,
\label{eq:B}
\end{equation}
where
\[
H
=
\int_{0+}^\infty \frac1y\,d\mathfrak F(y)
\in(0,\infty)
\]
is the Pickands constant.

Berman also obtained a compound Poisson approximation for
$\eta(u)L_u(T_u)$ with varying $T_u$ under the additional condition
\begin{equation}
\lim_{u\to\infty} \RR(u)\ln u=0,
\label{bermanCA}
\end{equation}
see \cite[Thm 7.4.1]{Berman92}. The approximation is determined by the index
$\alpha$ and the Berman function
\[
x\mapsto \int_x^\infty \frac1y\,d\mathfrak F(y).
\]

This approximation also applies to the supremum functional and relies on
extreme-value arguments. Its key ingredient is the Pickands lemma for short
intervals; see \Cref{lemma:Pickands}. This lemma determines the leading
constant but does not enter the normal approximation \zm{of}
\[
\frac{L_u(T_u)-\E{L_u(T_u)}}{\sqrt{\Var(L_u(T_u))}}
%\quad as\quad u\to\infty
\]
derived in \cite[Thm 8.2.1]{Berman92} for $T_u\to\infty$ at a suitable rate.
Central limit theorems and related results are available in other settings;
see e.g., 
\cite{bulinski2012central,pham2013rate,
spodarev2013limit,Spodarev2014,LeonenkoOlenko2014,kratz2018central,nualart2020continuous,
makogin2022limit,leonenko2023sojourn,leonenko2024limit,caponera2025sojourn,leonenko2025high}
and the references therein. Malliavin--Stein results such as
\cite{pham2013rate} address Gaussian fluctuations of sojourn times. We study
rare-exceedance asymptotics including the supremum functional. As emphasized
by Berman \cite[p.~142]{Berman92} these two regimes use fundamentally
different techniques. The compound Poisson regime relies on blocking from
extreme-value theory while the normal regime uses Gaussian methods related to
Hermite chaos expansions of sojourn functionals.

Berman's approach to sojourn tails was later complemented by the uniform
double-sum method in several one-dimensional settings; see e.g.,
\cite{debicki2019approximation,DLM20,novikov2025sojourns}. Related questions
for random fields were considered in \cite{chen2022extrema,DHLM23}.

We extend Berman's result \eqref{eq:A1} to sojourn functionals of
$\R^d$-valued GRFs. As a special case we derive asymptotic tail results for
the supremum functional $\kappa_u(0)$ extending \cite{PI23}. The
vector-valued setting is technically challenging because Slepian's inequality
is unavailable in higher dimensions. In tractable cases this can be overcome
by constructing suitable matrix-valued covariance functions \cite{IeN23}.

\textbf{Our contribution.}
The principal findings of this paper are:
\begin{enumerate}[(i)]
\item we extend Berman's approximation of sojourn functionals to
vector-valued GRFs and varying domains;

\item we extend the scalar uniform double-sum method of
\cite{debicki2023sojourns} to vector-valued GRFs;

\item we formulate conditions on the GRFs extending the applicability
of the technique in \cite{DHW20};

\item we represent the Pickands constants of \cite{DHW20} in terms of
Berman constants \zm{useful} for simulation.
\end{enumerate}

\textbf{Notation.}
All vectors in $\R^d$ are column vectors and are written in bold letters; for
example,
\[
\vk b=(b_1,\ldots,b_d)^\top,
\qquad
\vk0=(0,\ldots,0)^\top,
\qquad
\id=(1,\ldots,1)^\top .
\]
Vector operations are understood component-wise (in the Hadamard
sense); in particular, $\vk x\odot\vk y$ denotes the component-wise product
of two vectors of the same dimension. For two vectors $\vk x$ and $\vk y$,
we write $\vk x\le\vk y$ if
$x_i\le y_i$ for all $1\le i\le d$, and similarly for the other inequalities.
For \(c>0\) and $\vk\alpha=(\alpha_1,\ldots,\alpha_k)\in(0,2]^k$ we write
 
\[
   [\vk t]_{\vk\alpha}
   \coloneqq\sum_{i=1}^k |t_i|^{\alpha_i},
   \qquad 
   c^{1/\vk\alpha}\vk t
   \coloneqq
   \bigl(c^{1/\alpha_1}t_1,\ldots,c^{1/\alpha_k}t_k\bigr), \qquad \vk t\in\R^k.
\]
We denote by $\normE{\vk x}$ the Euclidean norm of a vector $\vk x$, and by
\[
\normF{A}
=
\sqrt{\sum_{1\le i,j\le d}a_{ij}^2}
\]
the Frobenius norm of a $d\times d$ real matrix $A=(a_{ij})$. All
unadorned matrix norms below are understood to be Frobenius norms; by
finite-dimensional norm equivalence, the particular choice is immaterial
for our estimates. We write
$\mathcal I_d$ for the $d\times d$ identity matrix and
\[
\diag (\vk x )
=
\diag(x_1,\ldots,x_d)
\]
for the diagonal matrix with entries $x_1,\ldots,x_d$ on the main diagonal.

\textbf{Brief organisation of the paper.}
Section~\ref{s:extreme} contains the main results and is followed by
illustrative examples. Section~\ref{s.proofs} contains the proofs and
auxiliary results. The Appendix collects technical lemmas.

\section{Main Results}
\label{s:extreme}

Throughout this section $\vk X(\vk t)$, $\vk t\in\R^k$, denotes a centered
$\R^d$-valued GRF with covariance matrix function
$R(\vk s,\vk t)$, $\vk s,\vk t\in\R^k$. Set
\[
\RR(\vk h)=R(\vk h,\vk0), \qquad \vk h\in\R^k,
\qquad
\Sigma\coloneq \RR(\vk0).
\]
We assume that $\Sigma$ is nonsingular and $\vk b\in\R^d$ satisfies the
{\it Savage condition}
\begin{equation}
\label{Savage}
\vk w\coloneq \Sigma^{-1}\vk b>\vk0 .
\end{equation}

If \eqref{Savage} fails the quadratic minimisation problem reduces to the
active coordinates and can be treated similarly to \cite{DHW20}. We omit
this case for notational simplicity.

The main assumptions used to extend \eqref{eq:A1} are the following:  

\begin{enumerate}[({B}0)]
\item\label{I:B0} The centered $\R^d$-valued GRF $\vk X$ is stationary,
the covariance matrix function
\[
   \RR(\vk h)=\Cov(\vk X(\vk h),\vk X(\vk0)),
   \qquad \vk h\in\R^k,
\]
is continuous, and the variance matrix $\Sigma $ is nonsingular.
\end{enumerate}

\begin{enumerate}[({B}1)]
\item\label{I:B1}
For all $\vk t\in\R^k\setminus\{\vk0\}$ we have 
\[
\vk w^\top(\Sigma-\RR(\vk t))\vk w>0 .
\]
\end{enumerate}
\begin{enumerate}[({B}2)]
\item\label{I:B2}
There is a positive function
$\vk v(u)\coloneqq(v_1(u),\ldots,v_k(u))$, $u>0$, whose components $v_i$
are regularly varying at infinity with indices $2/\alpha_i$ for
$\alpha_i\in(0,2]$, $i=1,\ldots,k$ \zm{and} there are also \(r>0\) and a
matrix-valued function \(V\) on \(\{\vk t:\|\vk t\|\le r\}\) such that
\begin{equation}
\label{e:stat-new}
\lim_{u\to\infty}
\sup_{\|\vk t\|\le r}
\left\|
u^2\bigl(\Sigma-\RR(\vk t/\vk v(u))\bigr)-V(\vk t)
\right\|
=0 .
\end{equation}
It also satisfies
\begin{equation}
\label{e:tangent-nondegeneracy}
   \vk w^\top V(\vk t)\vk w>0,
   \qquad
   0<\|\vk t\|\le r .
\end{equation}
\end{enumerate}
 
\begin{enumerate}[({B}3)]
\item\label{I:B3}
For every $\theta>0$ we have 
\[
m(\theta)
\coloneqq
\inf_{\|\vk t\|_\infty\ge \theta}
\vk w^\top(\Sigma-\RR(\vk t))\vk w
>0 .
\]
\end{enumerate}

By the matrix-valued Bochner--Cram\'er theorem
\cite{Cramer1940,Kleiber2017coherence}, there exists a finite Hermitian
nonnegative matrix-valued measure
$F=(F_{ij})_{1\le i,j\le d}$ on $\R^k$ such that (integration is understood
entrywise)
\begin{equation}
\label{eq:matrix-spectral-representation}
   \RR(\vk h)
   =
   \int_{\R^k}
   e^{\mathrm{i}\langle \vk h,\lambda\rangle}\,F(d\lambda),
   \qquad
   \vk h\in\R^k,
   \qquad
   \Sigma=F(\R^k).
\end{equation}
The diagonal measures $F_{ii}$ are the spectral measures of the component
fields, whereas the off-diagonal finite complex measures $F_{ij}$ are their
cross-spectral measures. For every real vector $\vk a\in\R^d$
\[
   F_{\vk a}(A)
   :=
   \vk a^\top F(A)\vk a,
   \qquad A\in\mathcal B(\R^k)
\]
is the finite nonnegative spectral measure of the scalar field
$\vk a^\top\vk X$. Since $F_{\vk a}$ is symmetric we have 
\begin{equation}
\label{eq:directional-spectral-loss}
   \vk a^\top\bigl(\Sigma-\RR(\vk h)\bigr)\vk a
   =
   \int_{\R^k}
   \bigl(1-\cos\langle\vk h,\lambda\rangle\bigr)
   F_{\vk a}(d\lambda).
\end{equation}
Taking $\vk a=\vk w$ gives direct spectral formulations of
Assumptions~\ref{I:B1} and~\ref{I:B3}. Moreover,
Assumption~\ref{I:B2} can be written equivalently as the locally uniform
Fourier--Stieltjes limit
\[
  \limit{u} u^2
   \int_{\R^k}
   \left(
      1-e^{\mathrm{i}\langle\vk t/\vk v(u),\lambda\rangle}
   \right)F(d\lambda)
   =
   V(\vk t).
\]
Thus the off-diagonal cross-spectral measures determine the cross-component
part, including any nonsymmetric part, of the tangent matrix function $V$.
Spectral tail or moment assumptions may be used to verify this limit in
specific models. We retain the covariance formulation because the proofs use
local covariance increments directly. \\
The convergence in Assumption~\ref{I:B2} is stronger than pointwise convergence.
 It is used as a transparent sufficient condition because the conditional tangent approximation, the vector Pickands lemma, the covariance increment estimates, and the two-block bound all require control that is uniform over compact sets of rescaled lags. In the matrix regular-variation examples below, local uniformity follows from the uniform convergence theorem away from the origin and from Potter bounds near the origin; in the concrete covariance constructions it is verified directly. 

Under Assumption~\ref{I:B0} the field $\vk X$ is $L^2$-continuous and its
components are stochastically continuous. Hence
\cite[Thm 5, p.~169 and Thm 1, p.~171]{MR636254} yields a jointly measurable
and separable version which we use below. If $\vk X$ has continuous sample
paths the continuity of $\mathcal{R}$ is automatic.

Assumption~\ref{I:B2} specifies \(V\) only locally.
Lemma~\ref{lemV} gives a unique anisotropically homogeneous extension to
\(\R^k\). We use \(V\) for this canonical extension. Then
\bqn{\label{RV2a}
   R_V(\vk t,\vk s)
   =
   V(\vk t)+V(-\vk s)-V(\vk t-\vk s),
   \qquad \vk s,\vk t\in\R^k,
}
is the covariance matrix function of a centered $\R^d$-valued GRF
$\vk Y$ with stationary increments and continuous sample paths \zm{(see Lemma~\ref{lemV} below)}.

Assumption~\ref{I:B3} is a weak global separation condition needed when
$T=T_u$ \zm{is changing} as $u\to\infty$.

Let $
\vk E=(E_1,\ldots,E_d)^\top$ 
be independent of \(\vk Y\) and have independent standard exponential
components. On a complete probability space define the possibly infinite
random variable (rv)
\[
\mathcal J(\vk Y,\vk w)
=
\int_{\R^k}
\ind{
\diag(\vk w)(\vk Y(\vk t)-V(\vk t)\vk w)
+
\vk E
>
\vk0
}
\,d\vk t .
\]
By Lemma~\ref{lemV} the field $\vk Y$ has a version with continuous sample
paths. Hence $\mathcal J(\vk Y,\vk w)$ is a well-defined
$[0,\infty]$-valued rv.
 
Let
\begin{equation}
\label{def.F}
\mathfrak F_{\vk Y,\vk w}(x)= \pk{\mathcal J(\vk Y,\vk w) \le x}, \qquad x\in [0, \IF]
\end{equation}
and set
\begin{equation}
\label{B.def.F}\mathcal B_{\vk Y,\vk w}(0)
=
\int_{0+}^\infty
\frac1y\,d\mathfrak F_{\vk Y,\vk w}(y), \qquad 
\mathcal B_{\vk Y,\vk w}(x)
=
\int_x^\infty \frac1y\,d\mathfrak F_{\vk Y,\vk w}(y), \qquad x>0.
\end{equation} 

The next result establishes continuity properties of both
$\mathfrak F_{\vk Y,\vk w}$ and $\mathcal B_{\vk Y,\vk w}$.

\begin{theorem}
\label{contBx} 
Under the Savage condition~\eqref{Savage} and
Assumptions~\ref{I:B0} and~\ref{I:B2},
$\mathfrak F_{\vk Y,\vk w}$ is continuous on $[0,\infty)$. If
$\mathcal B_{\vk Y,\vk w}(0)<\infty$ then
$ \mathcal B_{\vk Y,\vk w}$ is continuous on $[0,\infty)$.
\end{theorem}
 
In the remainder of the paper we set
\[
\eta(u)=\prod_{i=1}^k v_i(u),
\]
where $v_i$ are defined in Assumption~\ref{I:B2}. Recall that
for $T_u>0$ the sojourn functional $L_u(T_u)$ is defined in
\eqref{loc.time}. Let $F_u$ denote the distribution function of
$\eta(u)L_u(T_u)$.

\BT
\label{T:stat}
Assume the Savage condition~\eqref{Savage} and
Assumptions~\ref{I:B0}--\ref{I:B2}. Let $T_u>0$ satisfy
\bqn{\label{ta} \lim_{u\to\infty}
T_u v_i(u)=\infty,
\qquad i=1,\ldots,k 
} 
and
\bqn{\label{tb}
\lim_{u\to\infty} \frac{\ln T_u}{u^2}=0.
}
If $\limsup_{u\to\infty}T_u=\infty$ also assume
Assumption~\ref{I:B3}. Then for all $x\ge0$
\begin{equation}
\label{aH}
\lim_{u\to\infty}
\frac{
\int_0^x y\,dF_u(y)
}{
\eta(u)\E{L_u(T_u)}
}
=
\mathfrak F_{\vk Y,\vk w}(x).
\end{equation}
\ET

For the tail \qq{distribution asymptotics of $L_u(T_u) $,} we additionally assume \qq{sample paths continuity} of
$\vk X$ as required by the Borell--TIS argument.

\BT
\label{T:doublesum}
Assume the conditions of \Cref{T:stat} and that $\vk X$ has continuous sample
paths. Then
\begin{equation}
\label{doublesum}
\lim_{u\to\infty}
\frac{
\pk{
\eta(u) L_u(T_u) 
>x
}
}{
T_u^k\eta(u)\pk{\vk X(\vk0)>u\vk b}
}
=
\mathcal B_{\vk Y,\vk w}(x)
\end{equation}
holds for all $x\ge0$. The function
$\mathcal B_{\vk Y,\vk w}$ is continuous on $[0,\infty)$ and
\bqn{\label{kocheha}
\mathcal B_{\vk Y,\vk w}(x)
=
\lim_{S\to\infty}
\frac1{S^k}
\int_{\R^d}
\pk*{
\int_{[0,S]^k}
\ind{
\diag(\vk w)(\vk Y(\vk t)-V(\vk t)\vk w) >\vk z
}
\,d\vk t >x }
e^{\vk z^\top\vk1}\,d\vk z
\in(0,\infty).
}
 \ET

\begin{remark}\label{rem:cp-extension}
\begin{enumerate}[(i)]
\item
We assume that $\vk w=\Sigma^{-1}\vk b$ has positive components. This is the
Savage condition for the quadratic minimisation problem
$\Pi_{\Sigma}(\vk b)$ introduced in the Appendix. For general $\vk b$ with at
least one positive component let $\widetilde{\vk b}$ be the unique solution of
$\Pi_{\Sigma}(\vk b)$ and let $I$ be its determining index set. Put
\[
\vk w\coloneqq\Sigma^{-1}\widetilde{\vk b},
\qquad
\vk w_I=\Sigma_{II}^{-1}\vk b_I>\vk0_I,
\qquad
\vk w_J=\vk0_J,
\]
where $J=\{1,\ldots,d\}\setminus I$; see the Appendix. For this extension, retain the full $d\times d$ local covariance
assumption~\ref{I:B2}, with the vector $\vk w$ defined above. The limiting
tangent field and constants are then restricted to the coordinates indexed
by $I$. The normalising probability remains the full-vector tail
$\pk{\vk X(\vk0)>u\vk b}$. In particular, coordinates $j\in J$ for which
$\widetilde b_j=b_j$ may contribute a nontrivial Gaussian factor to this
one-point tail; they cannot simply be removed from the normalisation.
The separation assumptions~\ref{I:B1} and~\ref{I:B3} depend only on the
$I$-coordinates. For example, Assumption~\ref{I:B1} becomes
\[
\vk w_I^\top
(\Sigma-\RR(\vk t))_{II}
\vk w_I>0,
\qquad \vk t\in\R^k\setminus\{\vk0\}.
\]

\item
In the special case $x=0$ \Cref{T:doublesum} yields for every fixed
$T>0$
\begin{equation}
\label{Pavel}
\lim_{u\to\infty}
\frac{
\pk{\exists_{\vk t\in[0,T]^k}\,\vk X(\vk t)>u\vk b}
}{
T^k\eta(u)\pk{\vk X(\vk0)>u\vk b}
}
=
\mathcal B_{\vk Y,\vk w}(0).
\end{equation}
Thus $\mathcal B_{\vk Y,\vk w}(0)$ is the Pickands constant
appearing in \cite{DHW20}, where it was defined through the limit in
\eqref{kocheha}. The representation \eqref{B.def.F} is better suited 
\qq{for simulations} \zm{than representation \eqref{kocheha}}.

\item
The preceding results also yield a compound Poisson extension in the
vector-valued random-field setting. Assume all conditions of
\Cref{T:doublesum}, except that the growth condition~\eqref{tb} is replaced by
\[
\limit{u}T_u^k\eta(u)\pk{\vk X(\vk0)>u\vk b}
 = \vartheta\in(0,\infty) 
\]
and also assume the Berman-type decorrelation condition
\[
\lim_{\|\vk t\|_\infty\to\infty}
\bigl\|
\Sigma^{-1/2}\RR(\vk t)\Sigma^{-1/2}
\bigr\| 
\log \|\vk t\|_\infty
=0
\]
holds. A separated big-block/small-corridor argument in the spirit of
\cite[Thm 7.4.1]{Berman92} then implies
\[
\eta(u)L_u(T_u)\xrightarrow{d}\mathcal Z_\vartheta,
\qquad u\to\infty,
\]
where $\mathcal Z_\vartheta$ is a compound Poisson rv with intensity
$
\lambda
=
\vartheta\,\mathcal B_{\vk Y,\vk w}(0)>0
$
and jump-size survival function
$ 
\frac{\mathcal B_{\vk Y,\vk w}(x)}
{\mathcal B_{\vk Y,\vk w}(0)},\ x>0.
$
A proof is given in Section~\ref{s:compound-Poisson}. In particular,
\eqref{eq:cp-laplace-comparison} supplies a comparison estimate for the
joint Laplace transform of the retained block occupations. This controls
the cumulative dependence error over the growing number of blocks.
\end{enumerate}
\end{remark}

\section{Discussions and examples}

\qq{We first discuss the assumptions and connect them to existing literature.
Subsequently, we provide several illustrative examples.}

\medskip\noindent\textit{The scalar random-field case \qq{$d=1,k\ge1$}.}
Let \(d=1\) and write
\[
   \Sigma=\sigma^2,\qquad
   \gamma(\vk t)=\Sigma-\mathcal R(\vk t),
   \qquad
   w=\Sigma^{-1}b=b/\sigma^2>0 .
\]
Then Assumption~\ref{I:B1} reduces to
\[
   \gamma(\vk t)>0,\qquad \vk t\in\mathbb R^k\setminus\{\vk 0\}.
\]
In the unit-variance case this is the familiar assumption
\(\mathcal R(\vk t)<1\) for
\(\vk t\ne\vk 0\).
Assumption~\ref{I:B2} becomes  
\[
   u^2\gamma(\vk t/\vk v(u))\to V(\vk t), \qquad u \to \infty 
\]
locally uniformly in \(\vk t\), together with the non-degeneracy
$
   V(\vk t)>0,$ 
$
\vk t\ne\vk 0.
$
Thus after the normalisation
\(\sigma^2=b=1\), this is the analogue of the local condition
\(\limit{u}u^2(1-\mathcal R(\vk t/\vk v(u)))=h(\vk t)\) used in
\cite{chen2022extrema}, with \(h=V\).  
Assumption~\ref{I:B3} reduces to
\[
   \inf_{\|\vk t\|_\infty\ge\theta}\gamma(\vk t)>0,
   \qquad \theta>0,
\]
or in unit variance notation \qq{to}
$
   \sup_{\|\vk t\|_\infty\ge\theta}\mathcal R(\vk t)<1.
$
For fixed \(T\) this is a global version of the separation condition imposed
in \cite{chen2022extrema}. Our formulation is stronger because it also covers
varying domains \(T=T_u\). Thus Assumption~\ref{I:B3} implies
Assumption~\ref{I:B1}. Assumption~\ref{I:B1} is also imposed in
\cite[Condition~3]{kobelkov2019maximum}.

\medskip\noindent\textit{The vector-valued one-parameter case $d>1$, $k=1$.}
Consider a centered stationary $\R^d$-valued Gaussian process
$\vk X(t)$, $t\in\R$
and write
\[
   \Gamma(t)\coloneqq \Sigma-\RR(t),
   \qquad
   \vk w=\Sigma^{-1}\vk b>\vk0 .
\]

Assumption~\ref{I:B1} reduces to
\[
   \vk w^\top\Gamma(t)\vk w>0,
   \qquad t\ne0.
\]
 
Assumption~\ref{I:B2} has a particularly simple form. If it holds with
exponent \(\alpha\in(0,2]\), then
$
   u^2\Gamma(t/v(u))\to V(t)
$
locally uniformly, where $v$ is regularly varying at infinity with index \(2/\alpha\).  
The homogeneity of \(V\), together with \(V(-t)=V(t)^\top\), implies
\[
   V(t)
   =
   |t|^\alpha
   \left(
      A\ind{t\ge0}
      +
      A^\top\ind{t<0}
   \right),
   \qquad
   A\coloneqq V(1).
\]
Condition \eqref{e:tangent-nondegeneracy} then becomes
$
   \vk w^\top A\vk w>0.
$
A convenient sufficient condition for Assumption~\ref{I:B2} is the existence
of a positive function \(\ell\), slowly varying at zero, and a real
\(d\times d\) matrix \(A\) such that as \(t\to0\)
\begin{equation}
\label{eq:k1-matrix-rv}
\normF{
\Sigma-\RR(t)
-
|t|^\alpha\ell(|t|)
\left(
A\ind{t>0}
+
A^\top\ind{t<0}
\right)
}
=o\bigl(|t|^\alpha\ell(|t|)\bigr).
\end{equation}
\qq{Indeed,} let \(v(u)\) be specified asymptotically by
\[
   \lim_{u\to\infty}
   u^2v(u)^{-\alpha}\ell(1/v(u))=1 .
\]
By the asymptotic inversion theorem for regularly varying functions such a
function can be chosen with $v\in RV_{2/\alpha}$.
Then as \(u\to\infty\)
\[
u^2\Gamma(t/v(u))
\to
|t|^\alpha
\left(
A\ind{t\ge0}
+
A^\top\ind{t<0}
\right)
=
V(t)
\]
locally uniformly. Here uniformity away from $t=0$ follows from the uniform
convergence theorem, while Potter's bound and the positive index \(\alpha\)
control a neighbourhood of $t=0$.
If
$
   \vk w^\top A\vk w>0,$
then
\[
   \vk w^\top V(t)\vk w
   =
   |t|^\alpha\vk w^\top A\vk w>0,
   \qquad t\ne0,
\]
so Assumption~\ref{I:B2} holds. 
  
Assumption~\ref{I:B3} must be checked separately. In the one-parameter
vector-valued case our assumptions therefore reduce to one-dimensional
matrix regular variation and explicit global nondegeneracy conditions. Taking
\(x=0\) in the main results recover and extend \cite{DHW20}.

\medskip\noindent\textit{The vector-valued $k$-parameter additive case $d>1$.}
For simplicity suppose the Savage condition $\vk w>\vk0$ holds and write
\[
   S_{\alpha,A}(t)
   \coloneqq
   |t|^\alpha
   \left(A\ind{t\ge0}+A^\top\ind{t<0}\right),
   \qquad t\in\R .
\]
Consider the additive local structure
\[
   \Sigma-\RR(\vk t)
   =
   \sum_{i=1}^k S_{\alpha_i,A_i}(t_i)
   +
   o\left(\sum_{i=1}^k |t_i|^{\alpha_i}\right),
   \qquad \vk t\to\vk0,
\]
together with $\vk w^\top A_i\vk w>0$ for real $d\times d$ matrices $A_i$.
Then Assumption~\ref{I:B2} holds with
\[
   V(\vk t)=\sum_{i=1}^k S_{\alpha_i,A_i}(t_i).
\]
See \cite{PI23} for the related locally additive framework.

We next present examples illustrating the assumptions.

\begin{example} \label{ex:tangent-V}
Given \(\vk\alpha=(\alpha_1,\ldots,\alpha_k)\in(0,2]^k\), a useful class of
candidate tangent matrix functions is
\begin{equation}
\label{exV}
V(\vk t)
=
\sum_{i=1}^k \SAi{t_i}{V_i},
\qquad
\vk t=(t_1,\ldots,t_k)\in\R^k,
\end{equation}
where $V_i\in\R^{d\times d}$, $i=1,\ldots,k$, and
\[
\SAi{x}{V_i}
=
|x|^{\alpha_i}
\left(
V_i\ind{x\ge0}+V_i^\top\ind{x<0}
\right)
=
|x|^{\alpha_i}
\left(
V_i^+ + V_i^-\sgn(x)
\right),
\]
with
\[
V_i^+=\frac12(V_i+V_i^\top),
\qquad
V_i^-=\frac12(V_i-V_i^\top).
\]
We use the convention \(\sgn(x)=1\) for \(x\ge0\) and \(\sgn(x)=-1\) for
\(x<0\).

For every index \(i\) such that \(0<\alpha_i<2\) define
\begin{equation}
\label{e:posV}
V_i^\star
=
\sin\left(\frac{\pi\alpha_i}{2}\right)V_i^+
-
\sqrt{-1}\,
\cos\left(\frac{\pi\alpha_i}{2}\right)V_i^-.
\end{equation}
In view of  \cite[Prop.~9]{ACLP13} together with a direct treatment of the endpoint
\(\alpha_i=2\), the kernel
\[
   R_V(\vk t,\vk s)
   =
   V(\vk t)+V(-\vk s)-V(\vk t-\vk s)
\]
is a covariance matrix function if and only if the following conditions hold:
\begin{itemize}
\item if \(0<\alpha_i<2\), then \(V_i^\star\) in \eqref{e:posV} is
Hermitian nonnegative definite;
\item if \(\alpha_i=2\), then \(V_i=V_i^\top\) is nonnegative definite.
\end{itemize}
Indeed for \(\alpha_i=2\) the corresponding contribution to the covariance
kernel is \(2t_is_iV_i\). When $\alpha_i=1$, the antisymmetric term
$xV_i^-$ is additive and therefore cancels from $R_V$, see also
\cite{IeN23}. 
The identity
\[
   \vk w^\top V(\vk t)\vk w
   =
   \sum_{i=1}^k |t_i|^{\alpha_i}
   \vk w^\top V_i^+\vk w
\]
shows that the additional condition
$
   \vk w^\top V_i^+\vk w>0,
$
$i=1,\ldots,k,$
is necessary and sufficient for \eqref{e:tangent-nondegeneracy}. Here $V_i^+=V_i$ when $\alpha_i=2$.
To
construct stationary GRFs \(\vk X\) producing such structures, additional
positive-definiteness mechanisms are needed. One such mechanism is given in
the next example.
\end{example}

\begin{example} 
\label{ex:shifted-OU}
Let $Z(t)$, $t\in\R$ be a centered stationary Ornstein--Uhlenbeck (OU) process
with covariance
$
   \E{Z(t+h)Z(t)}=e^{-|h|},
   $
$h\in\R,$
and fix $\tau>0$. Define the stationary bivariate Gaussian process
\[
   \vk X(t)=\bigl(Z(t),Z(t+\tau)\bigr)^\top,
   \qquad t\in\R.
\]
Its covariance matrix function and variance matrix are
\[
   \RR(h)
   =
   \begin{pmatrix}
      e^{-|h|} & e^{-|h-\tau|}\\
      e^{-|h+\tau|} & e^{-|h|}
   \end{pmatrix},
   \qquad
   \Sigma
   =
   \begin{pmatrix}
      1 & e^{-\tau}\\
      e^{-\tau} & 1
   \end{pmatrix}, \qquad h\in \R.
\]
Since $\tau>0$ the matrix $\Sigma$ is positive definite. The process has
continuous sample paths and therefore satisfies Assumption~\ref{I:B0}.

Take $v(u)=u^2$, corresponding to $\alpha=1$. A direct expansion at the
origin gives, locally uniformly in $t$,
\[
   u^2\bigl(\Sigma-\RR(t/u^2)\bigr)
   \rightarrow
   V(t)
   =
   \begin{pmatrix}
      |t| & -e^{-\tau}t\\
      e^{-\tau}t & |t|
   \end{pmatrix}, \qquad u\to \infty.
\]
In particular $V(-t)=V(t)^\top$. For every $\vk w\ne\vk0$
\[
   \vk w^\top V(t)\vk w
   =
   |t|\,(w_1^2+w_2^2)>0,
   \qquad t\ne0,
\]
because the antisymmetric part of $V(t)$ disappears in the quadratic form.
Thus Assumption~\ref{I:B2} holds for every threshold $\vk b$ satisfying the
Savage condition $\vk w=\Sigma^{-1}\vk b>\vk0$.

For completeness Assumptions~\ref{I:B1} and~\ref{I:B3} can be verified
without a matrix positivity argument. Put
\[
   Y_{\vk w}(t)=w_1Z(t)+w_2Z(t+\tau),
   \qquad
   \Lambda(h)=\vk w^\top(\Sigma-\RR(h))\vk w.
\]
The OU covariance has the strictly positive spectral density
\[
   f(\lambda)=\frac{1}{\pi(1+\lambda^2)},
   \qquad \lambda\in\R.
\]
Hence for $h\ne0$ and $\vk w\ne\vk0$
\[
\begin{aligned}
   2\Lambda(h)
   &=\Var\bigl(Y_{\vk w}(h)-Y_{\vk w}(0)\bigr)\\
   &=\int_{\R}
      |e^{\sqrt{-1}\lambda h}-1|^2
      |w_1+w_2e^{\sqrt{-1}\lambda\tau}|^2
      f(\lambda)\,d\lambda
   >0.
\end{aligned}
\]
This proves Assumption~\ref{I:B1}. Moreover, $\Lambda$ is continuous and
$\Lambda(h)\to\vk w^\top\Sigma\vk w>0$ as $|h|\to\infty$. Its strict
positivity away from the origin therefore implies
\[
   \inf_{|h|\ge\theta}\Lambda(h)>0,
   \qquad \theta>0,
\]
which is Assumption~\ref{I:B3}. Thus this example satisfies all field-related
assumptions of \Cref{T:doublesum}. \qq{Moreover, $V$} belongs to the class \qq{of matrix valued functions considered} in
\Cref{ex:tangent-V}.
\end{example}

\begin{example} \label{ex:Hadamard-Schoenberg}
Let
$ 
   \Gamma:\R^k\to[0,\infty)^{d\times d}$ 
be continuous and conditionally negative definite in the matrix-valued
pseudo-cross-variogram sense. Assume that
$   \Gamma(\vk0)=\vk0_{d\times d},$
$   \Gamma(-\vk h)=\Gamma(\vk h)^\top,
$
and that, for some \(\vk\alpha=(\alpha_1,\ldots,\alpha_k)\in(0,2]^k\)
\[
   \Gamma(c^{1/\vk\alpha}\vk t)=c\Gamma(\vk t),
   \qquad c>0,\quad \vk t\in\R^k .
\]
By the matrix-valued Schoenberg theorem
\cite{IeN23,dorr2023characterization}
\[
   C(\vk h)
   =
   \exp^*(-\Gamma(\vk h))
   =
   \left(e^{-\Gamma_{ij}(\vk h)}\right)_{i,j=1}^d
\]
is a positive definite matrix-valued function. Let \(\Sigma\) be a positive
definite \(d\times d\) real matrix and define
\[
   \RR(\vk h)
   =
   \Sigma\zm{\odot} C(\vk h),
   \qquad \vk h\in\R^k,
\]
where \zm{\(\odot\)} denotes the Hadamard product. Since positive definite
matrix-valued functions are closed under Hadamard product, \(\RR\) is a
covariance matrix function of a centered stationary \(\R^d\)-valued GRF.

There is an important structural limitation in the present normalization.
Since $\Gamma(\vk0)=\vk0_{d\times d}$ we have
$C(\vk0)=\mathbf1\mathbf1^\top$. Hence all components of any field with
covariance $C$ agree in $L^2$ at each fixed location. It follows that
\[
   C_{ij}(\vk h)=c(\vk h),
   \qquad
   \Gamma_{ij}(\vk h)=\gamma(\vk h),
   \qquad 1\le i,j\le d
\]
for scalar functions \qq{$c(\vk h)=e^{-\gamma(\vk h)}$}. Consequently, we have 
\[
   \RR(\vk h)=c(\vk h)\Sigma,
   \qquad
   V(\vk t)=\gamma(\vk t)\Sigma.
\]
Thus the construction is valid but necessarily separable; unlike
\Cref{ex:shifted-OU}, it does not produce a genuinely nonseparable
matrix-valued function $V$.

Taking
$
   v_i(u)=u^{2/\alpha_i},
$
$i=1,\ldots,k 
$
by homogeneity, we obtain 
$
   \Gamma(\vk t/\vk v(u))=u^{-2}\Gamma(\vk t),
$
and therefore, locally uniformly in \(\vk t\)
\[
\begin{aligned}
u^2\bigl(\Sigma-\RR(\vk t/\vk v(u))\bigr)
&=
u^2\,
\Sigma\zm{\odot}
\left[
\mathbf 1\mathbf 1^\top
-
\exp^*(-u^{-2}\Gamma(\vk t))
\right]        \to
\Sigma\zm{\odot}\Gamma(\vk t), \qquad u \to \infty.
\end{aligned}
\]
Thus \qq{(\ref{e:stat-new}) in} Assumption~\ref{I:B2} holds with
$   V(\vk t)=\Sigma\zm{\odot}\Gamma(\vk t).$
Assume further that the Savage condition
$   \vk w=\Sigma^{-1}\vk b>\vk0
$
holds and that
$   \gamma(\vk t)>0,$
 $\vk t\ne\vk0$.
Then by continuity and homogeneity
\[
   \vk w^\top V(\vk t)\vk w
   =
   \gamma(\vk t)\vk w^\top\Sigma\vk w
   \asymp
   [\vk t]_{\vk\alpha},
   \qquad \vk t\in\R^k,
\]
so \eqref{e:tangent-nondegeneracy} holds. Furthermore, for
\(\vk h\ne\vk0\),
\[
\vk w^\top(\Sigma-\RR(\vk h))\vk w
=
\bigl(1-e^{-\gamma(\vk h)}\bigr)\vk w^\top\Sigma\vk w
>0,
\]
which proves Assumption~\ref{I:B1}. Finally, continuity, homogeneity, and
strict positivity of \(\gamma\) on the anisotropic unit sphere imply that
\(\gamma(\vk h)\) is bounded away from zero outside every neighbourhood of
the origin. Hence Assumption~\ref{I:B3} holds as well.

The corresponding Gaussian field has a continuous modification by
Kolmogorov's continuity theorem. 
\end{example}
 
\begin{example}\label{ex:operator-fOU}
Let \(H\) be a symmetric \(d\times d\) matrix with eigenvalues
$
   h_1,\ldots,h_d\in(0,1].
$ 
By the spectral theorem, there exists an orthogonal matrix
\(Q\in\mathbb{R}^{d\times d}\), whose columns form an orthonormal
eigenbasis of \(H\) such that
\[
   Q^\top Q=QQ^\top=\mathcal I_d,
   \qquad
   H=Q\diag(h_1,\ldots,h_d)Q^\top .
\]
If some eigenvalues are repeated, then \(Q\) need not be unique.
Set
$  
   h_*=\min_{1\le j\le d}h_j
$
and define
\[
   P_*
   =
   Q\diag\bigl(\ind{h_j=h_*},\,1\le j\le d\bigr)Q^\top .
\]
Although \(Q\) need not be unique, \(P_*\) is uniquely determined by
\(H\); it is the orthogonal projection onto the eigenspace associated
with the eigenvalue \(h_*\). 
Let \(A\) be a nonsingular \(k\times k\) matrix and define
\[
   \RR(\vk h)
   =
   Q\diag\left(
   e^{-\|A\vk h\|^{2h_1}},
   \ldots,
   e^{-\|A\vk h\|^{2h_d}}
   \right)Q^\top,
   \qquad \vk h\in\R^k .
\]
For every \(j\), the scalar function
\[
   \vk h\longmapsto e^{-\|A\vk h\|^{2h_j}}
\]
is positive definite on \(\R^k\), since \(2h_j\in(0,2]\).
Consequently, \(\RR\) is a covariance matrix function of a centered
stationary \(\R^d\)-valued GRF, with
\[
   \Sigma=\RR(\vk0)=\mathcal I_d.
\]

Let
\[
   \alpha=2h_*,
   \qquad
   v_i(u)=u^{2/\alpha}=u^{1/h_*},
   \qquad i=1,\ldots,k .
\]
Since all components of \(\vk v(u)\) are equal to \(u^{1/h_*}\), we have
\[
   \left\|A\frac{\vk t}{\vk v(u)}\right\|
   =
   u^{-1/h_*}\|A\vk t\|.
\]
Therefore, we have 
\[
\begin{split}
   &u^2\bigl(\mathcal I_d-\RR(\vk t/\vk v(u))\bigr)  =
   Q\diag\left(
   u^2\left(1-
   e^{-u^{-2h_1/h_*}\|A\vk t\|^{2h_1}}\right),
   \ldots,
   u^2\left(1-
   e^{-u^{-2h_d/h_*}\|A\vk t\|^{2h_d}}\right)
   \right)Q^\top .
\end{split}
\]
Using \(1-e^{-x}\sim x\) as \(x\downarrow0\), we obtain, locally
uniformly in \(\vk t\)
\[
   u^2\bigl(\mathcal I_d-\RR(\vk t/\vk v(u))\bigr)
   \rightarrow
   V(\vk t)
   =
   \|A\vk t\|^{2h_*}P_*,
   \qquad u\to\infty.
\]
Indeed, the \(j\)-th diagonal term converges to
\(\|A\vk t\|^{2h_*}\) when \(h_j=h_*\), and to \(0\) when
\(h_j>h_*\). Thus, the limit in \eqref{e:stat-new} from
Assumption~\ref{I:B2} holds. If further
\[
   \vk w>\vk0
   \qquad\text{and}\qquad
   \vk w^\top P_*\vk w>0,
\]
then
\[
   \vk w^\top V(\vk t)\vk w
   =
   \bigl(\vk w^\top P_*\vk w\bigr)
   \|A\vk t\|^{2h_*}.
\]
Since \(A\) is nonsingular, there exist constants \(c_1,c_2>0\) such
that
\[
   c_1\|\vk t\|
   \le \|A\vk t\|
   \le c_2\|\vk t\|,
   \qquad \vk t\in\R^k.
\]
Moreover, by equivalence of norms on \(\R^k\)
\[
   \|A\vk t\|^{2h_*}
   \asymp
   \sum_{i=1}^k |t_i|^{2h_*}
   =
   [\vk t]_{\vk\alpha},
   \qquad \vk t\to\vk0,
\]
establishing \eqref{e:tangent-nondegeneracy}. Assumption~\ref{I:B1} also holds. Indeed, for every
\(\vk h\ne\vk0\), nonsingularity of \(A\) gives \(\|A\vk h\|>0\), and
therefore
\[
   \mathcal I_d-\RR(\vk h)
   =
   Q\diag\left(
   1-e^{-\|A\vk h\|^{2h_1}},
   \ldots,
   1-e^{-\|A\vk h\|^{2h_d}}
   \right)Q^\top
\]
is positive definite. In fact, this verifies the corresponding
stronger full-matrix nondegeneracy condition.

Assumption~\ref{I:B3} follows similarly. For every \(\theta>0\) 
nonsingularity of \(A\) implies
\[
   \inf_{\|\vk h\|_\infty\ge\theta}\|A\vk h\|>0.
\]
Consequently, for each \(j\)
\[
   \sup_{\|\vk h\|_\infty\ge\theta}
   e^{-\|A\vk h\|^{2h_j}}<1,
\]
which yields the required uniform separation from the covariance
matrix at the origin.

Finally, the corresponding GRF admits a continuous modification.
Indeed, for \(\|\vk h\|\le1\)
\[
   1-e^{-\|A\vk h\|^{2h_j}}
   \le \|A\vk h\|^{2h_j}
   \le C\|\vk h\|^{2h_*}.
\]
Thus for every sufficiently large \(p\ge2\)  Gaussian moment
equivalence yields
\[
   \E{\bigl\|\vk X(\vk t+\vk h)-\vk X(\vk t)\bigr\|^p}
   \le C_p\|\vk h\|^{p h_*}.
\]
Choosing \(p\) such that \(p h_*>k\), Kolmogorov's continuity theorem
yields a continuous modification.

If \(h_1=\cdots=h_d\), then \(P_*=\mathcal I_d\), and the condition
\(\vk w^\top P_*\vk w>0\) is automatic for every \(\vk w\ne\vk0\).
\end{example} 

\begin{remark} 
One can pass from one-parameter examples to additive \(k\)-parameter examples
by taking independent coordinate processes. Namely, let
$
   \vk X^{(i)}(t_i)$, $t_i\in\R$, $i=1,\ldots,k,
$
be independent centered stationary \(\R^d\)-valued Gaussian processes, and set
\[
   \vk X(\vk t)=\sum_{i=1}^k \vk X^{(i)}(t_i),
   \qquad \vk t=(t_1,\ldots,t_k).
\]
Then
$
   \RR(\vk h)=\sum_{i=1}^k \RR_i(h_i),
$
and   $\Sigma=\sum_{i=1}^k\Sigma_i$.
Assume that $\Sigma$ is nonsingular and choose a threshold $\vk b$ for which
$\vk w=\Sigma^{-1}\vk b>\vk0$.
If
\[
   \Sigma_i-\RR_i(t)
   =
   |t|^{\alpha_i}
   \left(
   B_i\ind{t\ge0}+B_i^\top\ind{t<0}
   \right)
   +o(|t|^{\alpha_i}),
   \qquad t\to0,
\]
then with $v_i(u)=u^{2/\alpha_i}$ the tangent matrix of $\vk X$ is
\[
   V(\vk t)
   =
   \sum_{i=1}^k
   |t_i|^{\alpha_i}
   \left(
   B_i\ind{t_i\ge0}+B_i^\top\ind{t_i<0}
   \right).
\]
If \(\vk w^\top B_i\vk w>0\) for every \(i\), then
\(\vk w^\top V(\vk t)\vk w>0\) for every \(\vk t\ne\vk0\), verifying
\eqref{e:tangent-nondegeneracy}. The global nondegeneracy
assumptions~\ref{I:B1} and~\ref{I:B3} have to be checked for the resulting
covariance matrix $\RR$. If the coordinate processes have continuous sample
paths, then so does $\vk X$.
\end{remark}

\begin{example} \label{ex:biv-Matern}
Let \(d=2\), \(k\ge1\), and let \(M(r\mid\nu,a)\) denote the Matérn correlation
function
\[
   M(r\mid\nu,a)
   =
   \frac{2^{1-\nu}}{\Gamma(\nu)}(ar)^\nu K_\nu(ar),
   \qquad r\ge0 ,
\]
where \(a>0\), \(\nu>0\), and \(K_\nu\) is the modified Bessel function of the
second kind. Let \(A\) be a nonsingular \(k\times k\) matrix, fix
\[
   \nu\in(0,1),
   \qquad
   \nu_{12}>1,
   \qquad
   \rho\in(0,1)
\]
and define
\[
   \RR(\vk h)
   =
   \begin{pmatrix}
   M(\|A\vk h\|\mid\nu,1)
   &
   \rho M(\|A\vk h\|\mid\nu_{12},1)
   \\[1mm]
   \rho M(\|A\vk h\|\mid\nu_{12},1)
   &
   M(\|A\vk h\|\mid\nu,1)
   \end{pmatrix}.
\]
For this
model, the bivariate Matérn admissibility condition is
\begin{equation}
\label{eq:Matern-rho-bound}
   0<\rho
   \le
   \frac{
      \Gamma(\nu+k/2)\Gamma(\nu_{12})
   }{
      \Gamma(\nu)\Gamma(\nu_{12}+k/2)
   };
\end{equation}
see \cite[Thm.~3]{GneitingKleiberSchlather2010}. Assume
\eqref{eq:Matern-rho-bound}, so that $\RR$ is a positive definite
matrix-valued covariance function, then we have 
\[
   \Sigma=\RR(\vk0)
   =
   \begin{pmatrix}
   1 & \rho\\
   \rho & 1
   \end{pmatrix}.
\]
Assume the Savage condition
$
   \vk w=\Sigma^{-1}\vk b>\vk0,
$
and for \(\vk h\ne\vk0\) put \(r=\|A\vk h\|>0\). Since every Mat\'ern
correlation is strictly smaller than one away from the origin,
\[
\begin{aligned}
\vk w^\top(\Sigma-\RR(\vk h))\vk w
={}&
\bigl(1-M(r\mid\nu,1)\bigr)(w_1^2+w_2^2)\\
&+2\rho\bigl(1-M(r\mid\nu_{12},1)\bigr)w_1w_2
>0.
\end{aligned}
\]
Thus Assumption~\ref{I:B1} holds without any additional restriction on
\(\rho\) beyond the bivariate Mat\'ern admissibility condition. Let next 
\[
   \alpha=2\nu,
   \qquad
   v_i(u)=u^{2/\alpha}=u^{1/\nu},
   \quad i=1,\ldots,k .
\]
Using the standard Matérn expansion
\[
   M(r\mid\nu,1)
   =
   1-c_\nu r^{2\nu}+o(r^{2\nu}),
   \qquad r\downarrow0,
\]
where 
$
   c_\nu=\frac{\Gamma(1-\nu)}{2^{2\nu}\Gamma(1+\nu)},
$
and using $1-M(r\mid\nu_{12},1)=O(r^2)$, we obtain, locally uniformly in
$\vk t$,
\[
   u^2\bigl(\Sigma-\RR(\vk t/\vk v(u))\bigr)
   \to
   V(\vk t)
   =
   c_\nu\|A\vk t\|^{2\nu}\mathcal I_2 .
\]
Thus \qq{limit (\ref{e:stat-new}) in} Assumption~\ref{I:B2} holds. Moreover, we have 
\[
   \vk w^\top V(\vk t)\vk w
   =
   c_\nu\|A\vk t\|^{2\nu}\|\vk w\|^2>0,
   \qquad \vk t\ne\vk0 .
\]
Since \(A\) is nonsingular
\[
   \|A\vk t\|^{2\nu}
   \asymp
   \sum_{i=1}^k |t_i|^{2\nu}
   =
   [\vk t]_{\vk\alpha}
\]
as $\vk t\to\vk0$, implying \eqref{e:tangent-nondegeneracy}. Hence
Assumption~\ref{I:B2} is satisfied.
Assumption~\ref{I:B3} also holds. Indeed, the positive quadratic form
displayed above is continuous and tends to
$\vk w^\top\Sigma\vk w>0$
as \(\|\vk h\|\to\infty\). Hence, it is bounded away from zero on every set
$
   \{\vk h:\|\vk h\|_\infty\ge\theta\},
$
$\theta>0$.
The local Matérn expansion yields a continuous
modification of the corresponding Gaussian field. 
\end{example}

\section{Proofs}\label{s.proofs}
We first give several auxiliary lemmas.

\subsection{Auxiliary lemmas}\label{s.aux}

\qq{In what follows, functions $v_i(u)$, $u>0$ are defined as in~\ref{I:B2}.}
\begin{lemma}\label{lemV}
If Assumptions~\ref{I:B0} and~\ref{I:B2} hold, then the restriction of $V$
to the neighbourhood of the origin appearing in \eqref{e:stat-new} admits a unique anisotropically
homogeneous extension to \(\R^k\). In the sequel, \(V\) denotes this canonical
extension. With this convention, for every finite $R>0$
\begin{equation}\label{e:stat-all-r}
\lim_{u\to\infty}
\sup_{\|\vk t\|\le R}
\left\|
u^2\bigl(\Sigma-\RR(\vk t/\vk v(u))\bigr)-V(\vk t)
\right\|
=0.
\end{equation}
Moreover, \(V\) is continuous on \(\R^k\) and satisfies
\[
   V(\vk0)=\vk0_{d\times d},
   \qquad
   V(-\vk h)=V(\vk h)^\top,
   \qquad
   V(c^{1/\vk\alpha}\vk h)=cV(\vk h),
   \qquad c>0,\quad \vk h\in\R^k .
\]
Furthermore, $R_V$ defined in \eqref{RV2a} is the covariance matrix function
of a centered $\R^d$-valued GRF $\vk Y(\vk t)$, $\vk t\in\R^k$, with
stationary increments and cross-variogram
\[
   D(\vk h)
   =
   \frac12\bigl(V(\vk h)+V(-\vk h)\bigr),
   \qquad \vk h\in\R^k.
\]
The GRF $\vk Y$ admits a version with continuous sample paths.
\end{lemma}
\begin{proof}
Let \(r>0\) be as in \eqref{e:stat-new}, and put
\(
   B_r=\{\vk t\in\R^k:\|\vk t\|\le r\}.
\)
Since \(\RR\) is continuous the function
\(
   \vk t\mapsto
   u^2\bigl(\Sigma-\RR(\vk t/\vk v(u))\bigr)
\)
is continuous. Hence \eqref{e:stat-new} implies that \(V\) is continuous on
\(B_r\). Also, since \(\RR(\vk0)=\Sigma\), we have
\(
   V(\vk0)=\vk0_{d\times d}.
\)

We first construct the extension. Fix \(\vk t\in\R^k\), and choose \(a>1\)
such that
\[
   \vk s\coloneqq a^{-1/\vk\alpha}\vk t\in B_{r/2}.
\]
Let
\[
   w_u\coloneqq \frac{u}{\sqrt a},
   \qquad
   d_{u,i}\coloneqq
   a^{1/\alpha_i}\frac{v_i(w_u)}{v_i(u)},
   \qquad i=1,\ldots,k .
\]
Since \(v_i\in RV_{2/\alpha_i}\), we have \(d_{u,i}\to1\), \qq{as $u\to \infty$}. Moreover,
componentwise we have 
\[
   \frac{\vk t}{\vk v(u)}
   =
   \frac{\vk d_u\odot\vk s}{\vk v(w_u)} .
\]
For all large \(u\), \qq{we have} \(\vk d_u\odot\vk s\in B_r\). Therefore by
\eqref{e:stat-new} and the continuity of \(V\) on \(B_r\),
\[
\begin{aligned}
u^2\bigl(\Sigma-\RR(\vk t/\vk v(u))\bigr)
&=
a w_u^2
\bigl(
\Sigma-\RR((\vk d_u\odot\vk s)/\vk v(w_u))
\bigr)      \to
aV(\vk s), \qquad u\to\infty.
\end{aligned}
\]
The limit is independent of the admissible choice of \(a\), since the
left-hand side does not depend on \(a\). We define
\[
   V(\vk t)\coloneqq 
   aV(a^{-1/\vk\alpha}\vk t),
\]
where \(a>1\) is any number such that
\(a^{-1/\vk\alpha}\vk t\in B_{r/2}\). This defines \(V\) \qq{for all}
\(\R^k\). If \(\vk t\in B_r\), then \eqref{e:stat-new} gives the same limit
with the original local value \(V(\vk t)\), so the extension agrees with the
local function on \(B_r\). The extension is unique among anisotropically
homogeneous extensions, because any such extension is determined by its values
on \(B_r\).

The preceding argument also gives the locally uniform convergence. Indeed,
fix \(R>0\), and choose \(a_R>1\) so large that
\(
   a_R^{-1/\vk\alpha}\vk t\in B_{r/2},\)
   \(\|\vk t\|\le R.\)
Writing
\[
   w_u=u/\sqrt{a_R},
   \qquad
   \vk s=a_R^{-1/\vk\alpha}\vk t,
   \qquad
   d_{u,i}=a_R^{1/\alpha_i}\frac{v_i(w_u)}{v_i(u)},
\]
we have \(\vk d_u\to\vk1\) and
\[
   \frac{\vk t}{\vk v(u)}
   =
   \frac{\vk d_u\odot\vk s}{\vk v(w_u)} .
\]
The possible \(\vk s\)'s form a compact subset of \(B_{r/2}\). Hence, using
\eqref{e:stat-new} and the uniform continuity of \(V\) on compact subsets of
\(B_r\)
\[
\begin{aligned}
&\sup_{\|\vk t\|\le R}
\left\|
u^2\bigl(\Sigma-\RR(\vk t/\vk v(u))\bigr)-V(\vk t)
\right\|                                                   \\
&\quad \le
a_R
\sup_{\vk s}
\left\|
w_u^2
\bigl(
\Sigma-\RR((\vk d_u\odot\vk s)/\vk v(w_u))
\bigr)
-
V(\vk d_u\odot\vk s)
\right\|                                                   \\
&\qquad
+
a_R
\sup_{\vk s}
\left\|
V(\vk d_u\odot\vk s)-V(\vk s)
\right\|
\to0, \qquad u\to \IF.
\end{aligned}
\]
This proves \eqref{e:stat-all-r}. Since the convergence is locally uniform and
the approximating functions are continuous, \(V\) is continuous on \(\R^k\).

To verify homogeneity fix \(c>0\) and choose
\(b>\max\{1,c^{-1}\}\) so large that
\(b^{-1/\vk\alpha}\vk h\in B_{r/2}\). Then \(cb>1\) is admissible for
\(c^{1/\vk\alpha}\vk h\), and the definition of the extension gives
\[
   V(c^{1/\vk\alpha}\vk h)
   =cbV(b^{-1/\vk\alpha}\vk h)
   =cV(\vk h),
   \qquad c>0,\quad \vk h\in\R^k .
\]
Also, stationarity gives
\(    \RR(-\vk h)=\RR(\vk h)^\top \), and passing to the limit in \eqref{e:stat-all-r} yields
\(  V(-\vk h)=V(\vk h)^\top\), $\vk h\in\R^k$.

We now prove that \(R_V\) is a covariance matrix function. Let
\(\vk t_1,\ldots,\vk t_n\in\R^k\). The block covariance matrix of
\[
   \left(
   u(\vk X(\vk t_1/\vk v(u))-\vk X(\vk0)),
   \ldots,
   u(\vk X(\vk t_n/\vk v(u))-\vk X(\vk0))
   \right)
\]
has \((i,j)\)-block
\[
u^2
\left[
\RR((\vk t_i-\vk t_j)/\vk v(u))
-
\RR(\vk t_i/\vk v(u))
-
\RR(-\vk t_j/\vk v(u))
+
\Sigma
\right].
\]
By \eqref{e:stat-all-r} this block converges to
\[
   V(\vk t_i)+V(-\vk t_j)-V(\vk t_i-\vk t_j)
   =
   R_V(\vk t_i,\vk t_j).
\]
Thus every finite block matrix generated by \(R_V\) is nonnegative definite 
as a limit of covariance matrices. Hence \(R_V\) is a covariance matrix
function and defines a centered \(\R^d\)-valued GRF   \(\vk Y\).
Since $V(\vk0)=0$, we have $\vk Y(\vk0)=\vk0$ almost surely. Moreover, we have 
\[
   \E*{
   (\vk Y(\vk h)-\vk Y(\vk0))
   (\vk Y(\vk h)-\vk Y(\vk0))^\top
   }
   =
   V(\vk h)+V(-\vk h),
\]
so the cross-variogram is
\(
   D(\vk h)=\frac12\bigl(V(\vk h)+V(-\vk h)\bigr).
\)
A direct substitution in the formula for \(R_V\) yields the stationary increments of \(\vk Y\).
Next, by continuity \qq{and homogeneity of $V$}, and
compactness of
\(
   \{\vk h\in\R^k:[\vk h]_{\vk\alpha}=1\},
\)
there exists \(C>0\) such that
\(
   \|V(\vk h)\|
   \le
   C[\vk h]_{\vk\alpha},\)  \(\vk h\in\R^k\).
Therefore, for each coordinate process \(Y_j\)
\(
   \Var(Y_j(\vk t)-Y_j(\vk s))
   \le
   C[\vk t-\vk s]_{\vk\alpha}.
\)
\qq{Since} locally
\[
   [\vk t-\vk s]_{\vk\alpha}
   \le
   C\|\vk t-\vk s\|^{\alpha_*},
   \qquad
   \alpha_*:=\min_{1\le i\le k}\alpha_i>0, 
\]
we have for every \(p\ge2\) that
\[
   \E{|Y_j(\vk t)-Y_j(\vk s)|^p}
   \le
   C_p\|\vk t-\vk s\|^{p\alpha_*/2}
\]
locally. Choosing \(p\) so large that \(p\alpha_*/2>k\), Kolmogorov's
continuity theorem gives a continuous version of each coordinate process.
Hence \(\vk Y\) admits a version with continuous sample paths.
\end{proof}
\begin{lemma}\label{l.V.regularity}
If Assumptions~\ref{I:B0} and~\ref{I:B2} hold and $V$ denotes the canonical
homogeneous extension given by \Cref{lemV}, then  
\[
   \vk w^\top V(\vk t)\vk w>0,
   \qquad
   \vk t\in\R^k\setminus\{\vk0\}
\]
and there exist constants \(0<c_V<C_V<\infty\) such that
\begin{equation}
\label{additional}
   c_V[\vk t]_{\vk\alpha}
   \le
   \vk w^\top V(\vk t)\vk w
   \le
   C_V[\vk t]_{\vk\alpha},
   \qquad
   \vk t\in\R^k .
\end{equation}
\end{lemma}

\begin{proof}
Positivity on the punctured neighbourhood on which \(V\) is initially
defined is part of Assumption~\ref{I:B2}. By the homogeneity of the canonical
extension, it holds for every \(\vk t\ne\vk0\).

Now consider the compact set
   $S_{\vk\alpha}
   :=
   \{\vk t\in\R^k:[\vk t]_{\vk\alpha}=1\}$.
The function
$
   \vk t\mapsto \vk w^\top V(\vk t)\vk w
$
is continuous and strictly positive on \(S_{\vk\alpha}\). Therefore there
exist constants \(0<\zm{c}_V\le \zm{C}_V<\infty\) such that
\[
   c_V
   \le
   \vk w^\top V(\vk t)\vk w
   \le
   C_V,
   \qquad
   \vk t\in S_{\vk\alpha}.
\]
For \(\vk t\ne\vk0\) write
$   r_\alpha=[\vk t]_{\vk\alpha}$, $\vk e=r_\alpha^{-1/\vk\alpha}\vk t$.
Then \([\vk e]_{\vk\alpha}=1\), and by homogeneity,
$V(\vk t)=r_\alpha V(\vk e)$.
Hence
\[
   c_V[\vk t]_{\vk\alpha}
   \le
   \vk w^\top V(\vk t)\vk w
   \le
   C_V[\vk t]_{\vk\alpha}.
\]
The same inequality is trivial at \(\vk t=\vk0\). This proves
\eqref{additional}.
\end{proof}
Put
\[
   \Gamma(\vk h)\coloneqq\Sigma-\RR(\vk h),
   \qquad
   \Lambda(\vk h)\coloneqq \vk w^\top\Gamma(\vk h)\vk w.
\]
\begin{lemma} 
\label{lem:anisotropic-Potter}
Under Assumptions~\ref{I:B0} and~\ref{I:B2} the following assertions hold.
\begin{enumerate}[(i)]
\item There exist \(C_0>0\) and \(\theta_0>0\) such that
\begin{equation}
\label{eq:Gamma-Lambda-comparable}
   \|\Gamma(\vk h)\|
   \le C_0\Lambda(\vk h),
   \qquad
   0<\|\vk h\|<\theta_0 .
\end{equation}

\item There exist indices
\[
   \vk\beta=(\beta_1,\ldots,\beta_k),
   \qquad
   0<\beta_i<\alpha_i,\quad i=1,\ldots,k,
\]
such that, for every \(R>0\), there exist positive constants \(C_R,u_R\)
satisfying
\begin{equation}
\label{e:local-holder-matrix}
u^2
\left\|
\Sigma-\RR(\vk t/\vk v(u))
\right\|
\le
C_R[\vk t]_{\vk\beta},
\qquad
[\vk t]_{\vk\alpha}\le R,
\qquad
u\ge u_R .
\end{equation}

\item For every \(p\in(0,1)\), there exist finite positive constants
\(c_p,\theta_p,u_p\) such that
\begin{equation}
\label{e:scalar-potter-lower}
u^2
\vk w^\top
\bigl(\Sigma-\RR(\vk t/\vk v(u))\bigr)
\vk w
\ge
c_p[\vk t]_{\vk\alpha}^{p}
\end{equation}
whenever
\[
u\ge u_p,
\qquad
[\vk t]_{\vk\alpha}\ge1,
\qquad
\|\vk t/\vk v(u)\|<\theta_p .
\]

\item For every \(R>0\),
\begin{equation}
\label{e:relative-matrix-increment}
\lim_{a\to\infty}
\lim_{\theta\downarrow0}
\limsup_{u\to\infty}
\sup_{\mathcal D(u,\theta,a,R)}
\frac{
\left\|
\RR((\vk h+\vk r)/\vk v(u))
-
\RR((\vk h+\vk r')/\vk v(u))
\right\|
}{
\vk w^\top
\bigl(\Sigma-\RR(\vk h/\vk v(u))\bigr)
\vk w
}
=0,
\end{equation}
where
\[
\mathcal D(u,\theta,a,R)
=
\left\{
(\vk h,\vk r,\vk r')\in(\R^k)^3:
\begin{array}{l}
[\vk h]_{\vk\alpha}\ge a,
[\vk r]_{\vk\alpha}\le R,\quad [\vk r']_{\vk\alpha}\le R,\\[1mm]
\|\vk h/\vk v(u)\|<\theta,
\|(\vk h+\vk r)/\vk v(u)\|<\theta,
\|(\vk h+\vk r')/\vk v(u)\|<\theta
\end{array}
\right\}.
\]
\end{enumerate}
\end{lemma}

\begin{proof}
We prove \textup{(ii)}, \textup{(i)}, \textup{(iii)}, and \textup{(iv)}, in
that order. We first prove the local matrix bound in \textup{(ii)}. Choose
\(b\in(0,1/2)\), and then
choose \(\varepsilon>0\) so small that
\[
   1-b(2+\alpha^*\varepsilon)>0,
   \qquad
   \alpha^*\coloneqq\max_{1\le i\le k}\alpha_i .
\]
The assertion is trivial at \(\vk t=\vk0\). For
\(0<[\vk t]_{\vk\alpha}\le1\), put
\[
   q=[\vk t]_{\vk\alpha},
   \qquad
   u'=uq^{-b},
   \qquad
   t_i'=t_i\frac{v_i(u')}{v_i(u)},
   \qquad i=1,\ldots,k .
\]
Then \(\vk t'/\vk v(u')=\vk t/\vk v(u)\). Since \(u'\ge u\),
Potter's bounds for the regularly varying functions \(v_i\) apply uniformly
to the pair \((u,u')\) once \(u\) exceeds a fixed threshold independent of
\(\vk t\). Hence
\[
\begin{aligned}
   [\vk t']_{\vk\alpha}
   &=
   \sum_{i=1}^k |t_i|^{\alpha_i}
   \left(\frac{v_i(u')}{v_i(u)}\right)^{\alpha_i}                     
   \le
   C\sum_{i=1}^k |t_i|^{\alpha_i}
   q^{-b(2+\alpha_i\varepsilon)}                                     \le
   Cq^{1-b(2+\alpha^*\varepsilon)}
   \le C .
\end{aligned}
\]
Thus \(\vk t'\) ranges over a fixed compact set. Since \(u'\ge u\), the
locally uniform convergence
\eqref{e:stat-all-r} yields
\[
   u^2\|\Gamma(\vk t/\vk v(u))\|
   =q^{2b}
     u'^2\|\Gamma(\vk t'/\vk v(u'))\|
   \le Cq^{2b}.
\]
Set \(\beta_i=2b\alpha_i\). Since \(2b\in(0,1)\), subadditivity of
\(x\mapsto x^{2b}\) on \([0,\infty)\) yields
\[
   q^{2b}
   =
   \left(\sum_{i=1}^k|t_i|^{\alpha_i}\right)^{2b}
   \le
   \sum_{i=1}^k|t_i|^{2b\alpha_i}
   =
   [\vk t]_{\vk\beta}.
\]
Thus the required estimate follows when \([\vk t]_{\vk\alpha}\le1\). On
\(1<[\vk t]_{\vk\alpha}\le R\), the set of possible \(\vk t\)'s is contained in a
compact set, so \eqref{e:stat-all-r} gives a uniform bound for the
left-hand side, while the same subadditivity inequality yields
\[
   [\vk t]_{\vk\beta}
   =
   \sum_{i=1}^k\bigl(|t_i|^{\alpha_i}\bigr)^{2b}
   \ge
   \left(\sum_{i=1}^k|t_i|^{\alpha_i}\right)^{2b}
   =
   [\vk t]_{\vk\alpha}^{2b}
   >1.
\]
This proves \eqref{e:local-holder-matrix}.

We next prove \textup{(i)} by constructing an intrinsic scale. Since
\(v_i^{\alpha_i}\) is regularly varying with index 2, there exist \(U>0\) and constants
\(1<a_0<A_0<\infty\) such that
\begin{equation}
\label{eq:fixed-scale-ratio}
   a_0
   \le
   \frac{v_i(2s)^{\alpha_i}}{v_i(s)^{\alpha_i}}
   \le A_0,
   \qquad s\ge U,\quad i=1,\ldots,k .
\end{equation}
For \(\vk x\ne\vk0\) let
$H_{\vk x}(s)
   \coloneqq [\vk v(s)\odot\vk x]_{\vk\alpha}$.
Choose \(\theta_0>0\) so small that
\(H_{\vk x}(U)<1\) whenever \(0<\|\vk x\|<\theta_0\).
Moreover, \eqref{eq:fixed-scale-ratio} implies
\[
   a_0H_{\vk x}(s)
   \le H_{\vk x}(2s)
   \le A_0H_{\vk x}(s),
   \qquad s\ge U.
\]
Hence \(H_{\vk x}(2^nU)\) is strictly increasing and tends to infinity.
Let \(s=2^nU\) be its first point at which \(H_{\vk x}(s)\ge1\).
Since \(H_{\vk x}(U)<1\), we have \(n\ge1\), and the minimality of the
crossing together with the preceding display yields
\begin{equation}
\label{eq:intrinsic-shell}
   1
   \le H_{\vk x}(s)
   \le A_0H_{\vk x}(s/2)
   < A_0 .
\end{equation}
Since 
$K:=\{\vk z:1\le [\vk z]_{\vk\alpha}\le A_0\}$
is compact,  we can choose \(S_0\ge U\) so that the convergence in
\eqref{e:stat-all-r} satisfies
\[
   \sup_{\vk z\in K}
   \left\|
   s^2\Gamma(\vk z/\vk v(s))-V(\vk z)
   \right\|
   \le\varepsilon_0,
   \qquad s\ge S_0,
\]
where \(\varepsilon_0>0\) is chosen so that
\(\varepsilon_0\|\vk w\|^2<c_V/2\). Then choose \(N\ge1\) such that
\(2^NU\ge S_0\). Since
\[
   \sup_{0<\|\vk x\|\le\eta}
   H_{\vk x}(2^{N-1}U)
   \to0,
   \qquad \eta\downarrow0,
\]
we may decrease \(\theta_0\) so that
\(H_{\vk x}(2^{N-1}U)<1\) whenever
\(0<\|\vk x\|<\theta_0\). By the dyadic monotonicity established above,
the first crossing then satisfies
\(s\ge2^NU\ge S_0\). Put
\(\vk z=\vk v(s)\odot\vk x\). By
\Cref{l.V.regularity}, compactness of the shell, and
the last display
\[
   s^2\Lambda(\vk x)
   \ge
   \vk w^\top V(\vk z)\vk w
   -\varepsilon_0\|\vk w\|^2
   \ge c_V/2,
\]  
whereas
\[
   s^2\|\Gamma(\vk x)\|
   \le
   \sup_{\vk z\in K}\|V(\vk z)\|+\varepsilon_0.
\]
Thus there exist constants \(0<m<M<\infty\) such that
\begin{equation}
\label{eq:intrinsic-bounds}
   m
   \le s^2\Lambda(\vk x),
   \qquad
   s^2\|\Gamma(\vk x)\|
   \le M,
   \qquad 0<\|\vk x\|<\theta_0 .
\end{equation}
The two estimates in \eqref{eq:intrinsic-bounds} prove
\eqref{eq:Gamma-Lambda-comparable}.

We now prove \textup{(iii)}. Fix \(p\in(0,1)\), and take
\(\varepsilon\in(0,2)\) so small that \(2/(2+\varepsilon)>p\). 
For \(\vk x=\vk t/\vk v(u)\), choose the intrinsic scale
\(s=s(\vk x)\) above and put    $\vk z=\vk v(s)\odot\vk x$,
$\lambda=u/s$.
Let \(S_{\mathrm P}\ge U\) be a common threshold beyond which Potter's bounds
with tolerance \(\varepsilon\) hold for the regularly varying functions
\(v_i^{\alpha_i}\), \(i=1,\ldots,k\). Choose \(N_{\mathrm P}\ge1\) such
that \(2^{N_{\mathrm P}}U\ge S_{\mathrm P}\). Since
\[
   \sup_{0<\|\vk x\|\le\eta}
   H_{\vk x}(2^{N_{\mathrm P}-1}U)
   \to 0,
   \qquad \eta\downarrow0,
\]
we may decrease \(\theta_p\le\theta_0\) so that
\(H_{\vk x}(2^{N_{\mathrm P}-1}U)<1\) whenever
\(0<\|\vk x\|<\theta_p\). Dyadic monotonicity then gives
\[
   s(\vk x)\ge2^{N_{\mathrm P}}U\ge S_{\mathrm P}.
\]
Finally, take \(u_p\ge S_{\mathrm P}\). Thus both arguments \(s\) and \(u\)
exceed the Potter threshold whenever
\(\|\vk x\|<\theta_p\) and \(u\ge u_p\).
Since
\[
   [\vk t]_{\vk\alpha}
   =
   \sum_{i=1}^k |z_i|^{\alpha_i}\frac{v_i(u)^{\alpha_i}}{v_i(s)^{\alpha_i}},
\]
Potter's bounds, see e.g., \cite{resnick2007heavy} and \eqref{eq:intrinsic-shell} yield
\[
   [\vk t]_{\vk\alpha}
   \le
   \begin{cases}
      C\lambda^{2+\varepsilon},&\lambda\ge1,\\
      C\lambda^{2-\varepsilon},&0<\lambda<1.
   \end{cases}
\]
If \(\lambda\ge1\), then on \([\vk t]_{\vk\alpha}\ge1\)
\[
   \lambda^2
   \ge c[\vk t]_{\vk\alpha}^{2/(2+\varepsilon)}
   \ge c[\vk t]_{\vk\alpha}^p.
\]
If \(0<\lambda<1\) and \([\vk t]_{\vk\alpha}\ge1\), then
\[
   1\le [\vk t]_{\vk\alpha}\le C\lambda^{2-\varepsilon}\le C.
\]
The first and second inequalities also imply
$\lambda\ge C^{-1/(2-\varepsilon)}$.
Thus both \([\vk t]_{\vk\alpha}\) and \(\lambda^{-1}\) are uniformly bounded, and
hence again \(\lambda^2\ge c_p[\vk t]_{\vk\alpha}^p\). Consequently, in both cases
\[
   \lambda^2\ge c_p[\vk t]_{\vk\alpha}^p,
   \qquad [\vk t]_{\vk\alpha}\ge1 .
\]
Since $u^2\Lambda(\vk x)=\lambda^2s^2\Lambda(\vk x)$,
\eqref{eq:intrinsic-bounds} proves \eqref{e:scalar-potter-lower}.

It remains to prove the relative increment estimate in \textup{(iv)}. Fix
\(p\in(0,1)\), let \(\theta_p\) be the constant from \textup{(iii)}, and put
$\theta_*:=\min\{\theta_0,\theta_p\}$.
Since \qq{in the limit \eqref{e:relative-matrix-increment} we suppose}
\(\theta\downarrow0\), it suffices throughout to consider
\(0<\theta\le\theta_*\).
For
\(i,j\in\{1,\ldots,d\}\), stationarity implies
\[
\begin{split}
   \RR_{ij}(\vk x+\vk y)-\RR_{ij}(\vk x)
   =
   \RR_{ij}(\vk y)-\Sigma_{ij}
   +\Cov\bigl(
      X_i(\vk x+\vk y)-X_i(\vk x),
      X_j(\vk0)-X_j(\vk x)
   \bigr).
\end{split}
\]
Consequently, Cauchy-Schwarz gives
\begin{equation}
\label{eq:covariance-increment-bound}
\begin{split}
   |\RR_{ij}(\vk x+\vk y)-\RR_{ij}(\vk x)|
   \le
   |\Gamma_{ij}(\vk y)|
   +
   \sqrt{
   \Var(X_i(\vk y)-X_i(\vk0))
   \Var(X_j(\vk x)-X_j(\vk0))
   } .
\end{split}
\end{equation}
 By norm equivalence and
\eqref{eq:Gamma-Lambda-comparable}
\[
   \Var(X_j(\vk x)-X_j(\vk0))
   =
   \Gamma_{jj}(\vk x)+\Gamma_{jj}(-\vk x)
   \le C\|\Gamma(\vk x)\|
   \le C\Lambda(\vk x),
   \qquad 0<\|\vk x\|<\theta_0 .
\]
For \([\vk r]_{\vk\alpha}\le R\), \eqref{e:local-holder-matrix} also gives
\[
   u^2\|\Gamma(\vk r/\vk v(u))\|
   +u^2\Var(X_i(\vk r/\vk v(u))-X_i(\vk0))
   \le C_R .
\]
Here the variance in the last display equals
\[
   \Gamma_{ii}(\vk r/\vk v(u))
   +\Gamma_{ii}(-\vk r/\vk v(u)),
\]
so the estimate follows directly from \eqref{e:local-holder-matrix}.
\qq{Next, we apply} \eqref{eq:covariance-increment-bound} with
\(\vk x=\vk h/\vk v(u)\) and \(\vk y=\vk r/\vk v(u)\).
For \([\vk h]_{\vk\alpha}\ge1\), \(0<\theta\le\theta_*\), and all sufficiently large
\(u\), the constraint \(\|\vk h/\vk v(u)\|<\theta\) permits both
\eqref{eq:Gamma-Lambda-comparable} and
\eqref{e:scalar-potter-lower}; in particular, the denominator below is
positive. By norm equivalence and the preceding estimates
\[
\begin{aligned}
&\frac{
\|\RR((\vk h+\vk r)/\vk v(u))-\RR(\vk h/\vk v(u))\|
}{
\Lambda(\vk h/\vk v(u))
}     \le
C_R\left\{
\frac1{u^2\Lambda(\vk h/\vk v(u))}
+
\frac1{\sqrt{u^2\Lambda(\vk h/\vk v(u))}}
\right\}.
\end{aligned}
\]
Using \eqref{e:scalar-potter-lower}, we therefore obtain, uniformly for
\([\vk h]_{\vk\alpha}\ge a\ge1\) in a sufficiently small local window and for all
sufficiently large \(u\),
\[
\frac{
\|\RR((\vk h+\vk r)/\vk v(u))-\RR(\vk h/\vk v(u))\|
}{
\Lambda(\vk h/\vk v(u))
}
\le
C_{R,p}\left(a^{-p}+a^{-p/2}\right).
\]
For two perturbations \(\vk r,\vk r'\), use the triangle inequality through
\(\RR(\vk h/\vk v(u))\). The last display is uniform for sufficiently small
\(\theta\) and all sufficiently large \(u\); in particular, its denominator
is positive there. Taking first the \(\limsup\) in \(u\), then letting
\(\theta\downarrow0\), and finally \(a\to\infty\), proves
\eqref{e:relative-matrix-increment} in the stated order.
The two additional local restrictions involving \(\vk h+\vk r\) and
\(\vk h+\vk r'\) are not needed for this estimate, but are retained in
\(\mathcal D(u,\theta,a,R)\) for its later applications.
\end{proof}

\BEL \label{lemCont}
Let $Z(\vk t)$, $\vk t\in\R^k$, be a centered GRF   with
continuous covariance kernel
\[
R(\vk t,\vk s)=\Cov(Z(\vk t),Z(\vk s)),
\qquad \vk t,\vk s\in\R^k .
\]
Let $f:\R^k\to\R$ be continuous. Fix $m\in\N$ and
$x_1,\ldots,x_m\in\R$, and define
\[
q(\vk t_1,\ldots,\vk t_m)
\coloneqq
\pk{
Z(\vk t_i)-f(\vk t_i)\le x_i,\ 1\le i\le m
}.
\]
Then $q$ is continuous at a point
$(\vk t_1,\ldots,\vk t_m)\in\R^{mk}$ provided that
\begin{equation}
\label{e:cont-point-condition}
\pk{
Z(\vk t_i)=x_i+f(\vk t_i)
}
=0,
\qquad i=1,\ldots,m .
\end{equation}
In particular, $q$ is continuous on $\R^{mk}$ if
$\Var(Z(\vk t))>0$,
$\vk t\in\R^k$.
 
\EEL
\prooflem{lemCont}
Fix $(\vk t_1,\ldots,\vk t_m)\in\R^{mk}$ satisfying
\eqref{e:cont-point-condition}, and let
$(\vk s_{n1},\ldots,\vk s_{nm})\in\R^{mk}$ be such that
\[
\vk s_{ni}\to \vk t_i,
\qquad n\to\infty,\quad i=1,\ldots,m .
\]
Set
\[
X_n
=
\bigl(Z(\vk s_{n1}),\ldots,Z(\vk s_{nm})\bigr),
\qquad
X
=
\bigl(Z(\vk t_1),\ldots,Z(\vk t_m)\bigr),
\]
and
\[
\vk a_n
=
\bigl(x_1+f(\vk s_{n1}),\ldots,x_m+f(\vk s_{nm})\bigr),
\qquad
\vk a
=
\bigl(x_1+f(\vk t_1),\ldots,x_m+f(\vk t_m)\bigr).
\]
Since $f$ is continuous, $\vk a_n\to \vk a$.

For each $i=1,\ldots,m$, the continuity of the covariance kernel gives
\[
\E{
\bigl(Z(\vk s_{ni})-Z(\vk t_i)\bigr)^2
}
=
R(\vk s_{ni},\vk s_{ni})
+
R(\vk t_i,\vk t_i)
-
2R(\vk s_{ni},\vk t_i)
\to0 
\]
for every \(i\). Therefore,
$X_n\to X$
in probability as $n\to\infty$, and hence
$X_n-\vk a_n
\Rightarrow X-\vk a$.
Setting  $
A=(-\infty,0]^m 
$
we have  
\[
q(\vk s_{n1},\ldots,\vk s_{nm})
=
\pk{X_n-\vk a_n\in A},
\qquad
q(\vk t_1,\ldots,\vk t_m)
=
\pk{X-\vk a\in A}.
\]
The boundary of $A$ satisfies 
$\partial A\subset
\bigcup_{i=1}^m
\{\vk y\in\R^m:\ y_i=0\}$.
Hence
\[
\pk{X-\vk a\in\partial A}
\le
\sum_{i=1}^m
\pk{Z(\vk t_i)=x_i+f(\vk t_i)}
=
0
\]
by \eqref{e:cont-point-condition}. Thus $A$ is a continuity set for the law
of $X-\vk a$. By the Portmanteau theorem,
\[
\pk{X_n-\vk a_n\in A}
\to
\pk{X-\vk a\in A}.
\]
Equivalently, we have 
\[
q(\vk s_{n1},\ldots,\vk s_{nm})
\to
q(\vk t_1,\ldots,\vk t_m).
\]
This proves the claimed continuity at
$(\vk t_1,\ldots,\vk t_m)$.

If $\Var(Z(\vk t))>0$ for every $\vk t\in\R^k$, then each
$Z(\vk t_i)$ has a nondegenerate Gaussian distribution, so
\eqref{e:cont-point-condition} holds for all choices of
$(\vk t_1,\ldots,\vk t_m)$. This gives continuity of $q$ on $\R^{mk}$.
 
\QED

The next result is a modification of \cite[Lem 1.6.1]{Berman92}. For
$\vk t\in[0,T]^k$, $\lambda>0$, and $u>0$, define
\[
L_u(T,\vk t,\lambda)
\coloneqq
\int_{[0,T]^k\cap
\{\vk s:\, v_i(u)|s_i-t_i|\le\lambda,\ i=1,\ldots,k\}}
\ind{\vk X(\vk s)>u\vk b}\,d\vk s .
\]
We write
\[
\{\abs{\vk x}\not\le \vk y\}
\coloneqq
\{\abs{\vk x}\le \vk y\}^{c},
\qquad
\abs{\vk x}=(|x_1|,\ldots,|x_k|),
\]
where vector inequalities are understood componentwise,
\qq{and recall that $\eta(u)=\prod_{i=1}^k v_i(u)$.}

\begin{lemma}\label{lem:separation}
Let $\vk X(\vk t)$, $\vk t\in\R^k$, be a stationary $\R^d$-valued GRF 
    with jointly measurable sample paths. Let
$\vk b\in\R^d\setminus(-\IF,0]^d$, and let $T_u>0$, $u>0$, be such that
\[
T_u v_i(u)\to\infty,
\qquad i=1,\ldots,k,
\]
and
$\pk{\vk X(\vk0)>u\vk b}>0$, $u>0$.
Assume that
\begin{equation}\label{separation}
\lim_{\lambda\to\infty}
\limsup_{u\to\infty}
\frac{
\eta(u)
\int_{[-T_u,T_u]^k\cap\{\abs{\vk s}\not\le \lambda/\vk v(u)\}}
\pk{\vk X(\vk s)>u\vk b,\ \vk X(\vk0)>u\vk b}\,d\vk s
}{ \pk{\vk X(\vk0)>u\vk b}  
}
=0 .
\end{equation}
Assume further that, for every $x\ge0$
\[
\lim_{\lambda\to\infty}
\lim_{u\to\infty}
T_u^{-k}
\int_{[0,T_u]^k}
\pk*{
\eta(u)L_u(T_u,\vk t,\lambda)\le x
\mid
\vk X(\vk t)>u\vk b
}
\,d\vk t
=
F(x).
\]
If   $x\ge0$ is a continuity point of $F$, 
then
\begin{eqnarray}
\lim_{u\to\infty}
T_u^{-k}
\int_{[0,T_u]^k}
\pk*{
\eta(u)L_u(T_u)\le x
\mid
\vk X(\vk t)>u\vk b
}
\,d\vk t
=
F(x).\label{loc.limit}
\end{eqnarray}
\end{lemma}

\prooflem{lem:separation}
Put
\[
 D_u=[0,T_u]^k,\qquad
 E_u(\vk t)=\{\vk X(\vk t)>u\vk b\},\qquad
 p_u=\pk{E_u(\vk0)}
\]
and for  \(y\ge0\) define
\[
\begin{aligned}
 A_u(y)
 &\coloneqq
 T_u^{-k}\int_{D_u}
 \pk*{\eta(u)L_u(T_u)\le y\mid E_u(\vk t)}
 \,d\vk t,\\
 A_{u,\lambda}(y)
 &\coloneqq
 T_u^{-k}\int_{D_u}
 \pk*{\eta(u)L_u(T_u,\vk t,\lambda)\le y\mid E_u(\vk t)}
 \,d\vk t .
\end{aligned}
\]
Since \(L_u(T_u,\vk t,\lambda)\le L_u(T_u)\),
then
\[
 A_u(y)\le A_{u,\lambda}(y).
\]

Let $M_{u,\lambda}(\vk t)
 \coloneqq
 L_u(T_u)-L_u(T_u,\vk t,\lambda)$.
For \(x>\delta>0\),
\[
 A_u(x)
 \ge
 A_{u,\lambda}(x-\delta)-B_{u,\lambda,\delta},
\]
where
\[
 B_{u,\lambda,\delta}
 \coloneqq 
 T_u^{-k}\int_{D_u}
 \pk*{
 \eta(u)M_{u,\lambda}(\vk t)>\delta
 \mid E_u(\vk t)}
 \,d\vk t .
\]
Markov's inequality and Fubini's theorem imply
\[
\begin{aligned}
B_{u,\lambda,\delta}
&\le
\frac{\eta(u)}{\delta T_u^kp_u}
\int_{D_u}\int_{D_u\cap
\{\abs{\vk s-\vk t}\not\le\lambda/\vk v(u)\}}
\pk{E_u(\vk s),E_u(\vk t)}
\,d\vk s\,d\vk t\\
&\le
\frac{\eta(u)}{\delta p_u}
\int_{[-T_u,T_u]^k\cap
\{\abs{\vk h}\not\le\lambda/\vk v(u)\}}
\pk{\vk X(\vk h)>u\vk b,\ \vk X(\vk0)>u\vk b}
\,d\vk h .
\end{aligned}
\]
The last inequality follows from stationarity. Hence, \qq{by \eqref{separation}} we have 
\[
\lim_{\lambda\to\infty}\limsup_{u\to\infty}
B_{u,\lambda,\delta}=0.
\]
The preceding upper and lower bounds, followed first by \(u\to\infty\)
and then by \(\lambda\to\infty\), yield
\[
\limsup_{u\to\infty}A_u(x)\le F(x)
\]
and for \(x>\delta>0\)
\[
\liminf_{u\to\infty}A_u(x)\ge F(x-\delta).
\]
Letting \(\delta\downarrow0\) proves the assertion at every positive
continuity point of \(F\).

For completeness, consider \(x=0\). For each sample path, the Lebesgue
density theorem implies
\[
\int_{D_u}
\ind{E_u(\vk t)}
\ind{L_u(T_u,\vk t,\lambda)=0}
\,d\vk t=0.
\]
The same statement holds with \(L_u(T_u)\) in place of the local
sojourn time. Taking expectations and using stationarity shows that
$A_{u,\lambda}(0)=A_u(0)=0$.

Consequently \(F(0)=0\), and the assertion also holds at \(x=0\).
\QED

\newcommand{\varnothing}{\emptyset}
 The following integral estimate is Lemma~8 of \cite{PI23}, in the form
 needed below; see also \cite{DHW20}.
\begin{lemma} 
\label{lemma:integral_estimate2}
Let \(D\subset\R^n\) be a non-empty compact set, and fix
\(\bm w\in(0,\infty)^d\). Let $\{\bm\chi_{\bm z}(t):\bm z\in\R^d,\ t\in D\}$
be a family of jointly measurable and separable \(\R^d\)-valued GRFs. 
Suppose that for some \(0<\varepsilon\le \qq{\frac14\min_{1\le j\le d}w_j}\),   
there exist constants
\[
  \sigma^2>0,\qquad G\in[0,\infty),\qquad H\in[0,\infty),\qquad
  \gamma\in(0,2]
\]
such that, for all \(\bm z\in\R^d\)
\begin{equation}
\label{eq:int-est-mean}
  \sup_{\varnothing\ne F\subset\{1,\ldots,d\}}\sup_{t\in D}
  \bm w_F^\top\E{\bm\chi_{\bm z,F}(t)}
  \le
  G+\varepsilon\sum_{j=1}^d |z_j|,
\end{equation}
\begin{equation}
\label{eq:int-est-var}
  \sup_{\varnothing\ne F\subset\{1,\ldots,d\}}\sup_{t\in D}
  \Var\bigl(\eta_{\bm z,F}(t)\bigr)
  \le
  \sigma^2, \qquad \eta_{\bm z,F}(t)
  \coloneqq
  \bm w_F^\top
  \bigl(
    \bm\chi_{\bm z,F}(t)-\E{\bm\chi_{\bm z,F}(t)}
  \bigr)
\end{equation}
and
\begin{equation}
\label{eq:int-est-inc}
  \sup_{\varnothing\ne F\subset\{1,\ldots,d\}}\sup_{\bm z\in\R^d}
  \E*{
    \bigl(\eta_{\bm z,F}(t)-\eta_{\bm z,F}(s)\bigr)^2
  }
  \le
  H\|t-s\|^\gamma,
  \qquad s,t\in D .
\end{equation}
Then there exist constants \(C,c>0\), depending only on
\(d,\bm w,D,H,\gamma\), but not on
\(\varepsilon,G,\sigma^2\), such that
\[
  \int_{\R^d}
  e^{\bm w^\top\bm z}
  \pk{\exists\,t\in D:\ \bm\chi_{\bm z}(t)>\bm z}
  \,d\bm z
  \le
  C\exp\{c(G+\sigma^2)\}.
\] 
\end{lemma}

\begin{remark}\label{rem:uniform-integral-estimate}
Let \(K\subset(0,\infty)^d\) and \(D_0\subset\R^n\) be non-empty compact
sets, and put
\[
   w_-\coloneqq\inf_{\bm w\in K}\min_{1\le j\le d}w_j>0.
\]
If \(\bm w\in K\), \(D\subset D_0\) is non-empty and compact, and
\eqref{eq:int-est-mean}--\eqref{eq:int-est-inc} hold with common
\(H,\gamma\) and some \(0<\varepsilon\le w_-/4\), then the constants
\(C,c\) in \Cref{lemma:integral_estimate2} may be chosen depending only on
\(d,K,D_0,H,\gamma\). In particular, \(G\) and \(\sigma^2\) may vary with
\(\bm w\) and \(D\).
 
By \eqref{eq:int-est-inc} 
\[ 
  \left(\E*{(\eta_{\bm z,F}(t)-\eta_{\bm z,F}(s))^2}\right)^{1/2}
  \le \sqrt H\,\|s-t\|^{\gamma/2}.
\]
For any \(t_0\in D\)
\[
  \E{\sup_{t\in D}\eta_{\bm z,F}(t)}
  =
  \E*{\sup_{t\in D}
  (\eta_{\bm z,F}(t)-\eta_{\bm z,F}(t_0))}.
\]
Since every Euclidean cover of \(D_0\) also covers \(D\), Dudley's entropy
bound now gives \(\E{\sup_{t\in D}\eta_{\bm z,F}(t)}\le M\), uniformly in all admissible
\(\bm z,F,\bm w,D\), where \(M\) depends only on \(D_0,H,\gamma\).
Together with \eqref{eq:int-est-var}, the Borell--TIS inequality and
integration of its tail bound therefore yield
\[
  \E{\exp\{2\sup_{t\in D}\eta_{\bm z,F}(t)\}}
  \le C_0\exp\{c_0\sigma^2\}
\]
with uniform constants \(C_0,c_0\).

On the orthant
\[
  \Omega_F
  :=
  \{\bm z:z_i>0,\ i\in F,\ z_j<0,\ j\notin F\},
\]
write
\[
  A\coloneqq\sum_{i\in F}w_i z_i,
  \qquad
  B\coloneqq -\sum_{j\notin F}w_j z_j,
  \qquad
  \delta\coloneqq \frac{\varepsilon}{w_-}\le\frac14.
\]
Then \(\sum_{j=1}^d|z_j|\le(A+B)/w_-\). If
\(\bm\chi_{\bm z}(t)>\bm z\) for some \(t\in D\), then
\[
  \sup_{t\in D}\eta_{\bm z,F}(t)
  >
  A-G-\varepsilon\sum_{j=1}^d|z_j|
  \ge
  (1-\delta)A-\delta B-G.
\]
Consequently, Markov's inequality yields
\[
  \pk{\exists\,t\in D:\bm\chi_{\bm z}(t)>\bm z}
  \le
  C_0\exp\left\{
  2G+c_0\sigma^2-2(1-\delta)A+2\delta B
  \right\}.
\]
Multiplication by \(e^{\bm w^\top\bm z}=e^{A-B}\) leaves the factor
\[
  \exp\{-(1-2\delta)(A+B)\}.
\]
Under the changes of variables \(u_i=w_i z_i\), \(i\in F\), and
\(v_j=-w_jz_j\), \(j\notin F\), the Jacobian is at most \(w_-^{-d}\).
Since \(1-2\delta\ge1/2\), the integral over every such orthant is bounded
by \(C\exp\{c(G+\sigma^2)\}\), uniformly in \(\bm w,D\). The all-negative
orthant is bounded directly by
\[
  \int_{(-\infty,0)^d}e^{\bm w^\top\bm z}\,d\bm z
  =
  \prod_{j=1}^dw_j^{-1}
  \le w_-^{-d}.
\]
Summing over the \(2^d\) orthants proves the asserted uniform estimate.

This uniform form will be used for
\[
  \bm w_{u,0}=\Sigma_{u,0}^{-1}\bm b,
\]
which remains in a compact subset of \((0,\infty)^d\), and for translates
of the unit rectangles \(B_1\times B_2\), all contained in the fixed compact
set \([-1,1]^{2k}\).
\end{remark}

The following lemma is a random-field version of the uniform double-sum
estimate developed in \cite{DHW20}. Related field versions under explicit
additive local covariance assumptions appear in \cite[Lem 3]{PI23}.  
\qq{Let $x_+\coloneqq\max(x,0)$.}
\begin{lemma}\label{lem.double}
Under the Savage condition~\eqref{Savage}, suppose that
\(\vk X(\vk t)\), \(\vk t\in\R^k\), has continuous sample paths and
satisfies Assumptions~\ref{I:B0} and~\ref{I:B2}.
Let
\(\delta\in(0,1)\) and, for \(S>1\) and
\(\boldsymbol\ell\in\mathbb{Z}^k\), put
$  A_{\boldsymbol\ell}(S)\coloneqq \prod_{i=1}^k [S\ell_i,S(\ell_i+1)]$,
and
\[
  \rho_S(\boldsymbol\ell)
  \coloneq
  \inf_{\vk s\in A_{\vk0}(S),\ \vk t\in A_{\boldsymbol\ell}(S)}
  [\vk t-\vk s]_{\vk\alpha}
  =
  \sum_{i=1}^k
  \bigl(S(|\ell_i|-1)_+\bigr)^{\alpha_i}.
\]
Then there exists \(\theta_0=\theta_0(\delta)>0\) such that, for every
\(\theta\in(0,\theta_0)\), there exist constants
\(\mathcal C,\mathcal D>0\), independent of \(u\), \(S\), and
\(\boldsymbol\ell\) (but possibly depending on \(\delta\) and the field \zm{\(\vk X\)}),
with the following property: for every fixed \(S>1\)
there exists \(u_0(S,\theta)<\infty\) such that, for all
\(u\ge u_0(S,\theta)\)
\[
  \frac{\pk*{
  \exists\,\vk s\in A_{\vk0}(S),\
  \exists\,\vk t\in A_{\boldsymbol\ell}(S):
  \vk X(\vk s/\vk v(u))>u\vk b,\,
  \vk X(\vk t/\vk v(u))>u\vk b
  } }
  { \pk{\vk X(\vk0)>u\vk b}}
  \le
  \mathcal C
  \exp\left\{
  -\mathcal D\,\rho_S(\boldsymbol\ell)^{1-\delta}
  \right\}
\]
uniformly over all
$  \boldsymbol\ell\in\mathbb{Z}^k\setminus\{-1,0,1\}^k$
satisfying
$  |\ell_i|
  \le
  \left\lceil\frac{\theta v_i(u)}{S}\right\rceil$,
$i=1,\ldots,k$.

Moreover, there exist
\(\vartheta_\delta,C_\delta,D_\delta>0\) and
\(u_\delta<\infty\) 
\qq{such that} for every \(u\ge u_\delta\),
every \(r\ge1\), and all nonempty compact axis-parallel rectangles
\(B_1,B_2\subset\R^k\) with side lengths at most one, satisfying
\[
  \inf_{\vk s\in B_1,\ \vk t\in B_2}
  [\vk t-\vk s]_{\vk\alpha}
  \ge r,
  \qquad
  \sup_{\vk s\in B_1,\ \vk t\in B_2}
  \|(\vk t-\vk s)/\vk v(u)\|
  <\vartheta_\delta,
\]
\qq{the following inequality holds} %\zm{such that} %we have
\[
  \frac{
  \pk*{
  \exists\,\vk s\in B_1,\ \exists\,\vk t\in B_2:
  \vk X(\vk s/\vk v(u))>u\vk b,\,
  \vk X(\vk t/\vk v(u))>u\vk b
  }
  }{
  \pk{\vk X(\vk0)>u\vk b}
  }
  \le
  C_\delta\exp\{-D_\delta r^{1-\delta}\}.
\]
\end{lemma}

\begin{proof}
All constants below are positive and finite, and may change from line to
line. They are independent of \(u\), \(S\), and \(\boldsymbol\ell\), unless
explicitly stated otherwise. 
\qq{Let $\delta_0\in(0,\delta)$,
and for \(\vk h\in\R^k\), set}
\[
  \Lambda_u(\vk h)
  \coloneqq 
  \vk w^\top
  \bigl(\Sigma-\RR(\vk h/\vk v(u))\bigr)
  \vk w .
\]
Since \(\RR(-\vk h)=\RR(\vk h)^\top\), we have
$\Lambda_u(-\vk h)=\Lambda_u(\vk h)$.
By continuity of \(\RR\) at the origin, we may choose a deterministic
\(\bar\vartheta>0\) so small that, whenever \(\|\vk x\|<\bar\vartheta\),
the matrix
\[
   \frac14\bigl(2\Sigma+\RR(\vk x)+\RR(\vk x)^\top\bigr)
\]
is uniformly positive definite and its inverse times \(\vk b\) belongs to a
fixed compact set \(K\subset(0,\infty)^d\). We also take
\[
   \bar\vartheta
   \le
   \min\{\theta_{\mathrm{cmp}},\theta_{1-\delta_0}\},
\]
where \(\theta_{\mathrm{cmp}}>0\) denotes the constant called \(\theta_0\)
in \eqref{eq:Gamma-Lambda-comparable}, and
\(\theta_{1-\delta_0}\) is the constant in
\eqref{e:scalar-potter-lower}.

We first record two covariance estimates. By
\eqref{e:local-holder-matrix}, \qq{there exists $\vk \beta<\vk \alpha$ such that} for every \(R<\infty\), there exist
constants \(C_R<\infty\) and \(u_R<\infty\) such that, for \(u\ge u_R\)
\begin{equation}
\label{eq:local-matrix-bound-in-proof}
u^2
\left\|
\Sigma-\RR(\vk q/\vk v(u))
\right\|
\le
C_R[\vk q]_{\vk\beta},
\qquad
[\vk q]_{\vk\alpha}\le R .
\end{equation}
Fix \(R<\infty\) and \(\varepsilon_1>0\). For fixed \(u\) the
supremum in \eqref{e:relative-matrix-increment} is nondecreasing in
\(\theta\) and nonincreasing in \(a\); these monotonicities are preserved by
the \(\limsup\) in \(u\). Hence the iterated limit in
\eqref{e:relative-matrix-increment} allows us to choose, in this order,
\(a_0\ge1\), \(\vartheta\in(0,\bar\vartheta]\), and
\(u_0<\infty\) such that, for all \(u\ge u_0\),
\[
u^2
\left\|
\RR((\vk h+\vk r)/\vk v(u))
-
\RR((\vk h+\vk r')/\vk v(u))
\right\|
\le
\varepsilon_1 u^2\Lambda_u(\vk h)
\]
uniformly over all \(\vk h,\vk r,\vk r'\in\R^k\) satisfying
\[
[\vk h]_{\vk\alpha}\ge a_0,\qquad
[\vk r]_{\vk\alpha}\le R,\qquad
[\vk r']_{\vk\alpha}\le R,
\]
and
\begin{eqnarray}
\|\vk h/\vk v(u)\|<\vartheta,\qquad
\|(\vk h+\vk r)/\vk v(u)\|<\vartheta,\qquad
\|(\vk h+\vk r')/\vk v(u)\|<\vartheta .\label{3.local}
\end{eqnarray}
On the complementary bounded region
$1\le[\vk h]_{\vk\alpha}<a_0$,
\qq{we use that for}
\([\vk r]_{\vk\alpha},[\vk r']_{\vk\alpha}\le R\),  with
\(\alpha^*\coloneqq \max_{1\le i\le k}\alpha_i\), the elementary inequality
\(|x+y|^{\alpha_i}\le
2^{\alpha^*}(|x|^{\alpha_i}+|y|^{\alpha_i})\) gives the common bound
\[
  [\vk h+\vk r]_{\vk\alpha}
  \vee
  [\vk h+\vk r']_{\vk\alpha}
  \le
  2^{\alpha^*}(a_0+R)
  =:R' .
\]
After increasing \(u_0\), if necessary, so that \(u_0\ge u_{R'}\),
\eqref{eq:local-matrix-bound-in-proof} implies
\[
\begin{aligned}
&u^2
\left\|
\RR((\vk h+\vk r)/\vk v(u))
-
\RR((\vk h+\vk r')/\vk v(u))
\right\|                                      \\
&\quad\le
u^2
\left\|
\Sigma-\RR((\vk h+\vk r)/\vk v(u))
\right\|
+
u^2
\left\|
\Sigma-\RR((\vk h+\vk r')/\vk v(u))
\right\|
\le C .
\end{aligned}
\]
\qq{Thus,} for every
\(R<\infty\) and every \(\varepsilon_1>0\), there exist
\(\vartheta>0\), \(C<\infty\), and \(u_0<\infty\) such that, for all
\(u\ge u_0\),
\begin{equation}
\label{eq:matrix-increment-in-proof}
u^2
\left\|
\RR((\vk h+\vk r)/\vk v(u))
-
\RR((\vk h+\vk r')/\vk v(u))
\right\|
\le
C+\varepsilon_1 u^2\Lambda_u(\vk h)
\end{equation}
uniformly over all \(\vk h,\vk r,\vk r'\in\R^k\) satisfying
\[
[\vk h]_{\vk\alpha}\ge1,\qquad
[\vk r]_{\vk\alpha}\le R,\qquad
[\vk r']_{\vk\alpha}\le R,
\]
and the  \qq{three local restrictions (\ref{3.local})}. We shall later choose
\(\varepsilon_1>0\) sufficiently small and then keep the corresponding
\(\vartheta\le\bar\vartheta\) fixed. The compact set \(K\), and hence the
coefficient supplied by \Cref{lemma:integral_estimate2}, do not depend on
this choice of \(\varepsilon_1\).

We next prove a unit-block estimate. Let \(B_1,B_2\subset\R^k\) be nonempty
compact axis-parallel rectangles with side lengths at most one. Suppose that
\begin{equation}
\label{eq:unit-block-separation}
\inf_{\vk s\in B_1,\ \vk t\in B_2}
[\vk t-\vk s]_{\vk\alpha}
\ge r\ge1
\end{equation}
and that
\begin{equation}
\label{eq:unit-block-locality}
\|(\vk t-\vk s)/\vk v(u)\|<\vartheta,
\qquad
\vk s\in B_1,\quad \vk t\in B_2 .
\end{equation}
We claim that there exist constants \(C,D>0\), independent of
\(u\), \(B_1\), \(B_2\), and \(r\), such that, uniformly over all such
rectangles and all sufficiently large \(u\)
\begin{equation}
\label{eq:unit-block-double-new}
\pk*{
\exists\,\vk s\in B_1,\
\exists\,\vk t\in B_2:
\vk X(\vk s/\vk v(u))>u\vk b,\ 
\vk X(\vk t/\vk v(u))>u\vk b
}
\le
Cp_u\exp\{-Dr^{1-\delta_0}\}.
\end{equation}
Let $E=B_1\times B_2$.
For \(\vk p=(\vk s,\vk t)\in E\), define
\[
  \vk Z_u(\vk p)
  =
  \frac12
  \left(
  \vk X(\vk s/\vk v(u))
  +
  \vk X(\vk t/\vk v(u))
  \right).
\]
Choose an arbitrary \(\vk p_0=(\vk s_0,\vk t_0)\in E\), and set
$\vk h_0=\vk t_0-\vk s_0$, $\Sigma_{u,0}=\Var(\vk Z_u(\vk p_0))$.
Then
\[
  \Sigma_{u,0}
  =
  \frac14
  \left(
  2\Sigma+\RR(\vk h_0/\vk v(u))+\RR^\top(\vk h_0/\vk v(u))
  \right).
\]
By \eqref{eq:unit-block-locality},
$\|\vk h_0/\vk v(u)\|<\vartheta $.
Since \(\vartheta\le\bar\vartheta\), the choice of
\(\bar\vartheta\) implies that \(\Sigma_{u,0}\) is uniformly positive
definite. Moreover, the vector
\[
  \vk w_{u,0}:=\Sigma_{u,0}^{-1}\vk b
\]
belongs, uniformly in \(u\) and \(\vk p_0\), to the fixed compact set
\(K\subset(0,\infty)^d\) chosen above.

For \(\vk p\in E\), define
\[
  R_u(\vk p,\vk p_0)
  =
  \Cov(\vk Z_u(\vk p),\vk Z_u(\vk p_0)),
  \qquad
  A_u(\vk p)
  =
  \Sigma_{u,0}-R_u(\vk p,\vk p_0),
\]
and let
\[
   \vk W_u(\vk p)
   :=
   \vk Z_u(\vk p)
   -R_u(\vk p,\vk p_0)\Sigma_{u,0}^{-1}\vk Z_u(\vk p_0).
\]
\qq{Clearly,} \(\vk W_u\) is independent of
\(\vk Z_u(\vk p_0)\). For \(\vk z\in\R^d\), set
\begin{equation}
\label{eq:chi-regression-new}
   \vk\chi_{u,\vk z}(\vk p)
   =
   u\vk W_u(\vk p)
   -u^2A_u(\vk p)\Sigma_{u,0}^{-1}\vk b
   +A_u(\vk p)\Sigma_{u,0}^{-1}\vk z.
\end{equation}
By the Gaussian regression formula, this is a jointly measurable version of
the conditional law of
$ u\bigl(\vk Z_u(\cdot)-u\vk b\bigr)+\vk z$
given \(\vk Z_u(\vk p_0)=u\vk b-\vk z/u\). It is separable and continuous
in \(\vk p\), so the conditional path event used below is measurable. Since
\[
\{
\vk X(\vk s/\vk v(u))>u\vk b,\ 
\vk X(\vk t/\vk v(u))>u\vk b
\}
\subset
\{
\vk Z_u(\vk s,\vk t)>u\vk b
\},
\]
it suffices to estimate
\[
  \pk*{\exists\,\vk p\in E:\vk Z_u(\vk p)>u\vk b}.
\]
The events over \(E\) are measurable because \(\vk X\) has continuous sample
paths and \(E\) is compact. Conditioning on \(\vk Z_u(\vk p_0)\) and using
the change of variables
$\vk y=u\vk b-\frac{\vk z}{u}$
gives
\[
\begin{aligned}
&\pk*{\exists\,\vk p\in E:\vk Z_u(\vk p)>u\vk b}   =
u^{-d}
\int_{\R^d}
\pk*{\exists\,\vk p\in E:\vk\chi_{u,\vk z}(\vk p)>\vk z}
\varphi_{\Sigma_{u,0}}
\left(u\vk b-\frac{\vk z}{u}\right)
d\vk z,
\end{aligned}
\]
where \(\vk\chi_{u,\vk z}\) is defined in
\eqref{eq:chi-regression-new}.
The assumption  \eqref{Savage} implies, 
see e.g., \cite{HA2005}
\[
\qq{ p_u\coloneqq \pk{\vk X(\vk0)>u\vk b} } 
  \ge
  c u^{-d}
  \exp\left\{
  -\frac12u^2\vk b^\top\Sigma^{-1}\vk b
  \right\}
\]
for all sufficiently large \(u\). Moreover,
\[
\begin{aligned}
&\varphi_{\Sigma_{u,0}}
\left(u\vk b-\frac{\vk z}{u}\right)   =
C(\Sigma_{u,0})
\exp\left\{
-\frac12u^2\vk b^\top\Sigma_{u,0}^{-1}\vk b
+\vk b^\top\Sigma_{u,0}^{-1}\vk z
-\frac{1}{2u^2}\vk z^\top\Sigma_{u,0}^{-1}\vk z
\right\}.
\end{aligned}
\]
Since \(\Sigma_{u,0}\) stays in a compact neighbourhood of \(\Sigma\), the
determinant factor is uniformly bounded, and the final quadratic term in
\(\vk z\) is non-positive. Hence, uniformly in \(\vk z\in\R^d\) for all large \(u\)
\[
u^{-d}
\varphi_{\Sigma_{u,0}}
\left(u\vk b-\frac{\vk z}{u}\right)
\le
C p_u
\exp\left\{
-\frac12u^2
\vk b^\top(\Sigma_{u,0}^{-1}-\Sigma^{-1})\vk b
\right\}
e^{\vk w_{u,0}^{\top}\vk z}
\]
Consequently, we obtain 
\begin{equation}
\label{eq:density-reduction-new}
\pk*{\exists\,\vk p\in E:\vk Z_u(\vk p)>u\vk b}
\le
Cp_u
\exp\left\{
-\frac12u^2
\vk b^\top(\Sigma_{u,0}^{-1}-\Sigma^{-1})\vk b
\right\}
I_u ,
\end{equation}
where
\[
  I_u
  =
  \int_{\R^d}
  e^{\vk w_{u,0}^{\top}\vk z}
  \pk*{
  \exists\,\vk p\in E:
  \vk\chi_{u,\vk z}(\vk p)>\vk z
  }
  d\vk z .
\]
We estimate \(I_u\) using \Cref{lemma:integral_estimate2} and its uniform
form in \Cref{rem:uniform-integral-estimate}.
\qq{
We note that 
\(\vk w_{u,0}\in K\subset(0,\infty)^d\),
and after the translation
\[
  (\vk s,\vk t)\mapsto(\vk s-\vk s_0,\vk t-\vk t_0),
\]
the index set \(E\) is contained in a fixed compact subset of
\(\R^{2k}\).}
Thus, the constants in
\Cref{lemma:integral_estimate2} are uniform over all unit rectangles.
The representation \eqref{eq:chi-regression-new} yields
\begin{equation}
\label{eq:chi-mean-new}
\E{\vk\chi_{u,\vk z}(\vk p)}
=
-u^2A_u(\vk p)\Sigma_{u,0}^{-1}\vk b
+
A_u(\vk p)\Sigma_{u,0}^{-1}\vk z .
\end{equation}

We first note that
\begin{equation}
\label{eq:Au-small-uniform}\limit{u}
  \sup_{\vk p\in E}\|A_u(\vk p)\|=0 .
\end{equation}
The convergence in \eqref{eq:Au-small-uniform} is uniform over all admissible
unit rectangles \(B_1,B_2\), all \(\vk p_0\in B_1\times B_2\), and all
\(\vk p\in E\).
Indeed, the perturbations
\[
  (\vk s-\vk s_0)/\vk v(u),
  \qquad
  (\vk t-\vk t_0)/\vk v(u)
\]
converge uniformly to \(\vk0\), because the side lengths of \(B_1,B_2\) are
at most one, \qq{hence} both norms are bounded by
\[
   \frac{\sqrt{k}}{\min_{1\le i\le k}v_i(u)}
   \rightarrow0.
\]
The remaining covariance arguments differ from
\(\pm\vk h_0/\vk v(u)\) by these uniformly vanishing perturbations. Moreover,
by \eqref{eq:unit-block-locality}, all these arguments remain in the compact
ball \(\{\|\vk x\|\le \vartheta\}\). Uniform continuity of \(\RR\) on this
ball proves the asserted uniform form of \eqref{eq:Au-small-uniform}. Put
\[
   \varepsilon_K
   :=
   \frac18\inf_{\vk a\in K}\min_{1\le j\le d}a_j>0.
\]
This is smaller than the uniform admissible value in
\Cref{rem:uniform-integral-estimate}. Since \(K\) is compact and
\(\Sigma_{u,0}^{-1}\) is uniformly bounded, we may use
\eqref{eq:Au-small-uniform} to choose \(u\) so large that
\begin{equation}
\label{eq:Auz-small}
\sup_{\vk p\in E}
\left\|
A_u(\vk p)\Sigma_{u,0}^{-1}
\right\|
\le
\varepsilon_K/C_K,
\end{equation}
where \(C_K\) bounds the norms of the vectors in \(K\). Thus, for every
\(F\subset\{1,\ldots,d\}\), the contribution of the second term in
\eqref{eq:chi-mean-new} satisfies
\[
\sup_{\vk p\in E}
\left|
\vk w_{u,0,F}^{\top}
\bigl[A_u(\vk p)\Sigma_{u,0}^{-1}\vk z\bigr]_F
\right|
\le
\varepsilon_K\sum_{j=1}^d|z_j|.
\]

We now estimate the first term in \eqref{eq:chi-mean-new}. A direct
calculation gives, for \(\vk p=(\vk s,\vk t)\),
\[
\begin{aligned}
A_u(\vk p)
=
\frac14\{&
\Sigma-\RR((\vk s-\vk s_0)/\vk v(u))
+
\Sigma-\RR((\vk t-\vk t_0)/\vk v(u))        \\
&+
\RR(-\vk h_0/\vk v(u))
-
\RR((\vk s-\vk t_0)/\vk v(u))               
+
\RR(\vk h_0/\vk v(u))
-
\RR((\vk t-\vk s_0)/\vk v(u))
\}.
\end{aligned}
\]
The first two terms are bounded after multiplication by \(u^2\), uniformly
over \(E\), by \eqref{eq:local-matrix-bound-in-proof}. For the third term,
write
$  \vk s-\vk t_0=-\vk h_0+(\vk s-\vk s_0)$,
and apply \eqref{eq:matrix-increment-in-proof} with
$\vk h=-\vk h_0$, $\vk r=\vk0$, $\vk r'=\vk s-\vk s_0$.
For the fourth term, write
$\vk t-\vk s_0=\vk h_0+(\vk t-\vk t_0)$,
and apply \eqref{eq:matrix-increment-in-proof} with
$\vk h=\vk h_0$, $\vk r=\vk0$, $\vk r'=\vk t-\vk t_0$.
Because the rectangles have side lengths at most one
\qq{
\[
   \max\left([\vk s-\vk s_0]_{\vk \alpha}, [\vk t-\vk t_0]_{\vk \alpha}\right)
   \le k=:R_{\square}.
\]
}
Thus, both applications use the fixed
\(R=R_{\square}\) in \eqref{eq:matrix-increment-in-proof}. The three local
\qq{restrictions in \eqref{3.local}} follow from
\eqref{eq:unit-block-locality}, applied respectively to the pairs
$(\vk s_0,\vk t_0)$, $(\vk s,\vk t_0)$, $(\vk s_0,\vk t)$.
\qq{Recall that \(\Sigma_{u,0}^{-1}\vk b=\vk w_{u,0}\in K\), and} put
\[
   C_{\mathrm m}
   :=
   \max\left\{1,\sup_{\vk a\in K}\|\vk a\|^2\right\}
   <\infty.
\]
Using also \(\Lambda_u(-\vk h_0)=\Lambda_u(\vk h_0)\), and applying
\eqref{eq:matrix-increment-in-proof} with the same \(\varepsilon_1\) in both
cases, we obtain
\begin{equation}
\label{eq:Au-bound-new}
u^2
\sup_{\vk p\in E}
\|A_u(\vk p)\|
\le
C+\varepsilon_1u^2\Lambda_u(\vk h_0).
\end{equation}

By the definition of \(C_{\mathrm m}\),
\eqref{eq:chi-mean-new}, \eqref{eq:Auz-small}, and
\eqref{eq:Au-bound-new} imply that, for every
\(F\subset\{1,\ldots,d\}\)
\[
\sup_{\vk p\in E}
\vk w_{u,0,F}^{\top}
\E{\vk\chi_{u,\vk z,F}(\vk p)}
\le
C+C_{\mathrm m}\varepsilon_1u^2\Lambda_u(\vk h_0)
+
\varepsilon_K\sum_{j=1}^d |z_j| .
\]
Thus the mean condition \qq{(\ref{eq:int-est-mean})} in \Cref{lemma:integral_estimate2} holds with
$G=C+C_{\mathrm m}\varepsilon_1u^2\Lambda_u(\vk h_0)$.

We next verify the variance and increment conditions. We shall use the
following consequence of \eqref{eq:local-matrix-bound-in-proof}: for every
bounded set of \(\vk q\)'s and every bounded set of vectors \(a\in\R^d\)
\begin{equation}
\label{eq:linear-increment-bound}
u^2
\Var\left(
a^\top
\left[
\vk X(\vk q/\vk v(u))-\vk X(\vk0)
\right]
\right)
\le
C[\vk q]_{\vk\beta}.
\end{equation}
Indeed, we have 
\[
\begin{aligned}
\Var\left(
a^\top
\left[
\vk X(\vk q/\vk v(u))-\vk X(\vk0)
\right]
\right)                                      
&=
a^\top
\left(
2\Sigma-\RR(\vk q/\vk v(u))-\RR^\top(\vk q/\vk v(u))
\right)
a                                           \\
& =
a^\top
\left[
\Sigma-\RR(\vk q/\vk v(u))
+
\bigl(\Sigma-\RR(\vk q/\vk v(u))\bigr)^\top
\right]
a ,
\end{aligned}
\]
and \eqref{eq:linear-increment-bound} follows from
\eqref{eq:local-matrix-bound-in-proof}.

For \(F\subset\{1,\ldots,d\}\), let \(\vk c_F\in\R^d\) agree with
\(\vk w_{u,0}\) on \(F\) and vanish off \(F\). By
\eqref{eq:chi-regression-new}, the centered conditional field is
\(u\vk W_u\). Moreover, \(\vk c_F^\top\vk W_u(\vk p)\) is the residual
from the linear regression of \(\vk c_F^\top\vk Z_u(\vk p)\) on the full
vector \(\vk Z_u(\vk p_0)\). This residual has the smallest variance among
all linear prediction errors. Taking
\(\vk c_F^\top\vk Z_u(\vk p_0)\) as a particular predictor therefore implies
\[
\begin{aligned}
\Var\left(
\vk w_{u,0,F}^{\top}
\left[
\vk\chi_{u,\vk z,F}(\vk p)
-
\E{\vk\chi_{u,\vk z,F}(\vk p)}
\right]
\right)                                      \le
u^2
\Var\left(
\vk w_{u,0,F}^{\top}
[\vk Z_u(\vk p)-\vk Z_u(\vk p_0)]_F
\right).
\end{aligned}
\]
Moreover, we have 
\[
\vk Z_u(\vk p)-\vk Z_u(\vk p_0)
=
\frac12
\left[
\vk X(\vk s/\vk v(u))-\vk X(\vk s_0/\vk v(u))
\right]
+
\frac12
\left[
\vk X(\vk t/\vk v(u))-\vk X(\vk t_0/\vk v(u))
\right].
\]
Using \eqref{eq:linear-increment-bound}, the boundedness of
\(\vk w_{u,0}\), and the inequality
\(\Var(Y+Z)\le2\Var(Y)+2\Var(Z)\), we get
\begin{eqnarray}
\sup_{\vk z\in\R^d}
\sup_{F\subset\{1,\ldots,d\}}
\sup_{\vk p\in E}
\Var\left(
\vk w_{u,0,F}^{\top}
\left[
\vk\chi_{u,\vk z,F}(\vk p)
-
\E{\vk\chi_{u,\vk z,F}(\vk p)}
\right]
\right)
\le C .\label{var.c}
\end{eqnarray}
Thus the variance condition \qq{(\ref{eq:int-est-var})} in \Cref{lemma:integral_estimate2} holds.

For the increment condition \qq{(\ref{eq:int-est-inc})}, put
\[
\eta_{\vk z,F}(\vk p)
=
\vk w_{u,0,F}^{\top}
\left(
\vk\chi_{u,\vk z,F}(\vk p)
-
\E{\vk\chi_{u,\vk z,F}(\vk p)}
\right).
\]
Let \(\vk p=(\vk s,\vk t)\) and
\(\vk q=(\vk s',\vk t')\) be points of \(E\). The same projection argument,
now applied to
\(\vk c_F^\top[\vk Z_u(\vk p)-\vk Z_u(\vk q)]\), shows that the residual
variance is no larger than the unconditional variance, by taking the zero
predictor. Moreover,
\[
\vk Z_u(\vk p)-\vk Z_u(\vk q)
=
\frac12
\left[
\vk X(\vk s/\vk v(u))-\vk X(\vk s'/\vk v(u))
\right]
+
\frac12
\left[
\vk X(\vk t/\vk v(u))-\vk X(\vk t'/\vk v(u))
\right].
\]
Therefore, by \eqref{eq:linear-increment-bound},
\[
\E{
\left[
\eta_{\vk z,F}(\vk p)-\eta_{\vk z,F}(\vk q)
\right]^2}
\le
C
\left(
[\vk s-\vk s']_{\vk\beta}
+
[\vk t-\vk t']_{\vk\beta}
\right).
\]
Let
$\gamma=\min_{1\le i\le k}\beta_i>0 $.
Since \(B_1\) and \(B_2\) have side lengths at most one,
\[
[\vk s-\vk s']_{\vk\beta}
+
[\vk t-\vk t']_{\vk\beta}
\le
C\|\vk p-\vk q\|^\gamma
\]
whenever \(\|\vk p-\vk q\|\le1\). For
\(\|\vk p-\vk q\|>1\), the same bound follows after increasing the constant,
using the preceding uniform variance bound. Hence there exists \(H<\infty\)
such that
\begin{eqnarray}
\sup_{F\subset\{1,\ldots,d\}}
\sup_{\vk z\in\R^d}
\E{
\left[
\eta_{\vk z,F}(\vk p)-\eta_{\vk z,F}(\vk q)
\right]^2}
\le
H\|\vk p-\vk q\|^\gamma,
\qquad
\vk p,\vk q\in E .\label{var.d}
\end{eqnarray}
Thus the increment condition \qq{(\ref{eq:int-est-inc})} in \Cref{lemma:integral_estimate2} holds.

Applying \Cref{lemma:integral_estimate2} in the uniform form of
\Cref{rem:uniform-integral-estimate},
we obtain
\begin{equation}
\label{eq:Iu-bound-new}
  I_u
  \le
  C_1
  \exp\left\{
  C_2\bigl(1+C_{\mathrm m}\varepsilon_1
  u^2\Lambda_u(\vk h_0)\bigr)
  \right\}.
\end{equation}
Here \(C_2\) depends only on the fixed compact set, the fixed compact weight
set \(K\), and the \qq{constants in (\ref{var.c}) and (\ref{var.d})}. In particular,
neither \(C_2\) nor \(C_{\mathrm m}\) depends on \(\varepsilon_1\); the
additive constant \(C_1\) may depend on the corresponding bounded-region
cutoff.

It remains to estimate the density exponent in
\eqref{eq:density-reduction-new}. For \(\vk h\in\R^k\) define
\[
\Sigma_{u,\vk h}
=
\frac14
\left(
2\Sigma+\RR(\vk h/\vk v(u))+\RR^\top(\vk h/\vk v(u))
\right),
\qquad
M_{u,\vk h}:=\Sigma-\Sigma_{u,\vk h}.
\]
Then \(M_{u,\vk h}\) is symmetric and, for every \(\vk x\in\R^d\)
\[
\vk x^\top M_{u,\vk h}\vk x
=
\frac14
\Var\left(
\vk x^\top
\bigl[\vk X(\vk h/\vk v(u))-\vk X(\vk0)\bigr]
\right)
\ge0 .
\]
Thus \(M_{u,\vk h}\) is non-negative definite.
Since \(\vartheta\le\bar\vartheta\), \(\Sigma_{u,\vk h}\) is positive
definite whenever \(\|\vk h/\vk v(u)\|<\vartheta\). 
\qq{Using that}
\[
\Sigma_{u,\vk h}^{-1}-\Sigma^{-1}
=
\Sigma^{-1}M_{u,\vk h}\Sigma^{-1}
+
\Sigma^{-1}M_{u,\vk h}\Sigma_{u,\vk h}^{-1}
M_{u,\vk h}\Sigma^{-1},
\]
we observe that both terms on the right are non-negative definite. Hence, with
\(\vk w=\Sigma^{-1}\vk b\),
\[
\begin{aligned}
\vk b^\top
\left(
\Sigma_{u,\vk h}^{-1}-\Sigma^{-1}
\right)
\vk b
&\ge
\vk w^\top M_{u,\vk h}\vk w       
=
\frac12
\vk w^\top
\left(
\Sigma-\RR(\vk h/\vk v(u))
\right)
\vk w
=
\frac12\Lambda_u(\vk h).
\end{aligned}
\]
Applying this with \(\vk h=\vk h_0\), we obtain
\begin{equation}
\label{eq:density-lambda-lower}
u^2
\vk b^\top
\left(
\Sigma_{u,0}^{-1}-\Sigma^{-1}
\right)
\vk b
\ge
\frac12u^2\Lambda_u(\vk h_0).
\end{equation}

Combining \eqref{eq:density-reduction-new}, \eqref{eq:Iu-bound-new}, and
\eqref{eq:density-lambda-lower}, we get
\[
\pk*{\exists\,\vk p\in E:\vk Z_u(\vk p)>u\vk b}
\le
Cp_u
\exp\left\{
-\frac14u^2\Lambda_u(\vk h_0)
+
C_2C_{\mathrm m}\varepsilon_1u^2\Lambda_u(\vk h_0)
+
C
\right\}.
\]
Choose \(\varepsilon_1>0\) so small that
\(C_2C_{\mathrm m}\varepsilon_1<1/8\), and
then fix the corresponding \(\vartheta\). It follows that
\begin{equation}
\label{eq:unit-before-potter}
\pk*{\exists\,\vk p\in E:\vk Z_u(\vk p)>u\vk b}
\le
Cp_u
\exp\left\{
-Du^2\Lambda_u(\vk h_0)
\right\}.
\end{equation}

By \eqref{eq:unit-block-separation},
\[
  [\vk h_0]_{\vk\alpha}
  =
  [\vk t_0-\vk s_0]_{\vk\alpha}
  \ge r\ge1 .
\]
Moreover, by \eqref{eq:unit-block-locality} we have 
\[
  \|\vk h_0/\vk v(u)\|<\vartheta\le\theta_{1-\delta_0}.
\]
Therefore \eqref{e:scalar-potter-lower}, with \(p=1-\delta_0\) implies 
for all large \(u\)
\[
u^2\Lambda_u(\vk h_0)
\ge
c_{1-\delta_0}
[\vk h_0]_{\vk\alpha}^{1-\delta_0}
\ge
c_{1-\delta_0}r^{1-\delta_0}.
\]
Substituting this into \eqref{eq:unit-before-potter} proves
\eqref{eq:unit-block-double-new}.
All lower bounds on \(u\) used in this unit-block argument depend only on the
fixed field parameters, \(\delta_0,\varepsilon_1,\vartheta,K\), and
\(R_{\square}\); in particular, they are independent of \(B_1,B_2\), and
\(r\). Since \(r\ge1\) and \(\delta_0<\delta\), the stronger bound
\eqref{eq:unit-block-double-new} also proves the unit-rectangle assertion in
the statement, with \(\vartheta_\delta=\vartheta\) after relabelling the
constants.

We now pass from unit blocks to \(S\)-blocks. Let
$\alpha_*=\min_{1\le i\le k}\alpha_i$.
Choose \(\theta_0>0\) so small that
$2\sqrt{k}\,\theta_0<\vartheta$.
Fix \(\theta\in(0,\theta_0)\) and \(S>1\). Let next 
\[
\boldsymbol\ell\in\mathbb{Z}^k\setminus\{-1,0,1\}^k,
\qquad
|\ell_i|
\le
\left\lceil\frac{\theta v_i(u)}{S}\right\rceil,
\qquad i=1,\ldots,k .
\]
Partition each of \(A_{\vk0}(S)\) and \(A_{\boldsymbol\ell}(S)\) into
subrectangles contained in the corresponding block and having side lengths at
most one. The number of subrectangles needed for each block
is bounded by \(C(1+S)^k\). Hence, the number of pairs of rectangles is
bounded by \(C(1+S)^{2k}\).

For every
$\vk s\in A_{\vk0}(S)$, $\vk t\in A_{\boldsymbol\ell}(S)$,
we have coordinatewise
\[
  |t_i-s_i|
  \ge
  S(|\ell_i|-1)_+,
  \qquad i=1,\ldots,k 
\]
implying 
\[
[\vk t-\vk s]_{\vk\alpha}
\ge
\sum_{i=1}^k
\bigl(S(|\ell_i|-1)_+\bigr)^{\alpha_i}
=
\rho_S(\boldsymbol\ell).
\]
Since
\(\boldsymbol\ell\notin\{-1,0,1\}^k\), there exists \(j\) such that
\(|\ell_j|\ge2\) and hence   
\begin{equation}
\label{eq:rho-lower-S}
  \rho_S(\boldsymbol\ell)
  \ge S^{\alpha_j}
  \ge S^{\alpha_*}
  >1 .
\end{equation}

Moreover, we have 
\[
|t_i-s_i|
\le
S(|\ell_i|+1)
\le
\theta v_i(u)+2S .
\]
Since \(S\) is fixed and \(v_i(u)\to\infty\), there exists
\(u_0(S,\theta)<\infty\) such that for all \(u\ge u_0(S,\theta)\)
\[
\frac{|t_i-s_i|}{v_i(u)}
\le
2\theta,
\qquad i=1,\ldots,k 
\] and therefore 
\[
\|(\vk t-\vk s)/\vk v(u)\|
\le
2\sqrt{k}\,\theta
<
\vartheta .
\]
Thus the unit-block estimate \eqref{eq:unit-block-double-new} applies to
every pair of unit rectangles, with
$r=\rho_S(\boldsymbol\ell)$.
Summing over all pairs of rectangles yields
\[
\pk*{
  \exists\,\vk s\in A_{\vk0}(S),\
  \exists\,\vk t\in A_{\boldsymbol\ell}(S):
  \vk X(\vk s/\vk v(u))>u\vk b,\,
  \vk X(\vk t/\vk v(u))>u\vk b
  } 
\le
C(1+S)^{2k}p_u
\exp\{-D\rho_S(\boldsymbol\ell)^{1-\delta_0}\}.
\]

It remains to absorb the polynomial factor. By \eqref{eq:rho-lower-S}
\[
  (1+S)^{2k}
  \le
  C\rho_S(\boldsymbol\ell)^{2k/\alpha_*}.
\]
Since \(\delta_0<\delta\), there exist constants
\(\widetilde C_\delta,\widetilde D_\delta>0\) such that for all \(x\ge1\)
\[
  x^{2k/\alpha_*}e^{-Dx^{1-\delta_0}}
  \le
  \widetilde C_\delta
  e^{-\widetilde D_\delta x^{1-\delta}} .
\]
Taking \(x=\rho_S(\boldsymbol\ell)\), we obtain
\[ \pk*{
  \exists\,\vk s\in A_{\vk0}(S),\
  \exists\,\vk t\in A_{\boldsymbol\ell}(S):
  \vk X(\vk s/\vk v(u))>u\vk b,\,
  \vk X(\vk t/\vk v(u))>u\vk b
  } 
\le
\mathcal C p_u
\exp\left\{
-\mathcal D\,\rho_S(\boldsymbol\ell)^{1-\delta}
\right\}.
\]
The constants \(\mathcal C,\mathcal D\) are independent of
\(u\), \(S\), and \(\boldsymbol\ell\). The lower threshold on \(u\) may
depend on the fixed value of \(S\) and on \(\theta\). This proves the lemma.
\end{proof}

\subsection{Proof of Theorem~\ref{contBx}}\label{proofconBx}
We first present a lemma needed in the proof of Theorem~\ref{contBx}.

\begin{lemma}\label{lem:levelset-null-shift}
Let \(K=[a,b]^k\), \(a<b\), or \(K=\R^k\), and let
\(f:K\to\R^d\) be continuous. Then, for every
\(x\in(0,\lambda_k(K))\),
\[
  \lambda_d\left(
  \left\{\vk s\in(0,\infty)^d:
  \int_K \ind{f(\vk t)+\vk s>\vk0}\,d\vk t=x
  \right\}
  \right)=0,
\]
where, if \(K=\R^k\), the condition \(x\in(0,\lambda_k(K))\) means
\(x\in(0,\infty)\). Consequently, if \(h:\R^d\to\R\) is Lebesgue measurable
and finite almost everywhere, then
\[
  \int_{(0,\infty)^d}
  \ind{\int_K \ind{f(\vk t)+\vk s>\vk0}\,d\vk t=x}\,
  h(\vk s)\,d\vk s=0 .
\]
\end{lemma}

\begin{proof}
All vector inequalities are understood componentwise. It is enough to prove
the stronger assertion
\[
  \lambda_d\left(
  \left\{\vk s\in\R^d:
  \int_K \ind{f(\vk t)+\vk s>\vk0}\,d\vk t=x
  \right\}
  \right)=0 .
\]
The level set is measurable, since
$(\vk t,\vk s)\mapsto \ind{f(\vk t)+\vk s>\vk0}$
is Borel measurable on \(K\times\R^d\), and hence
\[
  \vk s\mapsto
  \int_K \ind{f(\vk t)+\vk s>\vk0}\,d\vk t
\]
is measurable by Fubini's theorem.

We first consider the case \(d=1\). Put next 
\[
  J(s)\coloneqq \lambda_k\{\vk t\in K:f(\vk t)>-s\},
  \qquad s\in\R .
\]
Then \(J\) is nondecreasing. We show that, for every
\(x\in(0,\lambda_k(K))\), the equation \(J(s)=x\) has at most one solution.
Assume that \(J(s)=J(s')=x\) for some \(s<s'\). Since \(-s'<-s\) we have 
\[
  0
  =
  J(s')-J(s)
  =
  \lambda_k\{\vk t\in K:-s'<f(\vk t)\le -s\}
\]
and hence 
\[
  U\coloneqq\{\vk t\in K:-s'<f(\vk t)<-s\}
\]
has Lebesgue measure zero. The set \(U\) is relatively open in \(K\). Since
every non-empty relatively open subset of either \([a,b]^k\) or \(\R^k\) has
positive \(k\)-dimensional Lebesgue measure, we must have \(U=\emptyset\).
Consequently,
\[
  f(K)\cap(-s',-s)=\emptyset .
\]

Since \(K\) is connected and \(f\) is continuous, \(f(K)\) is a connected
subset of \(\R\). Hence either
\[
  f(\vk t)\le -s'
  \qquad\text{for all }\vk t\in K,
\]
or
\[
  f(\vk t)\ge -s
  \qquad\text{for all }\vk t\in K.
\]
In the first case \(J(s')=0\), contradicting \(x>0\). In the second case
\(J(s')=\lambda_k(K)\), contradicting \(x<\lambda_k(K)\). Therefore
\(J(s)=x\) has at most one solution, and so
$\lambda_1\{s\in\R:J(s)=x\}=0$.

Now let \(d\ge2\). Write
\[
  \vk s
  =
  (r,r+q_2,\ldots,r+q_d),
  \qquad
  r\in\R,\quad \vk q=(q_2,\ldots,q_d)\in\R^{d-1}.
\]
This linear change of variables has determinant \(1\). For fixed
\(\vk q\in\R^{d-1}\) define
\[
  g_{\vk q}(\vk t)
  :=
  \min\bigl(
    f_1(\vk t),
    f_2(\vk t)+q_2,
    \ldots,
    f_d(\vk t)+q_d
  \bigr),
  \qquad \vk t\in K .
\]
Then \(g_{\vk q}\) is continuous and
\[
  f(\vk t)+(r,r+q_2,\ldots,r+q_d)>\vk0
  \quad\Longleftrightarrow\quad
  g_{\vk q}(\vk t)+r>0 .
\]
By the scalar case for every fixed \(\vk q\in\R^{d-1}\)
\[
  \lambda_1\left(
  \left\{
  r\in\R:
  \lambda_k\{\vk t\in K:g_{\vk q}(\vk t)>-r\}=x
  \right\}
  \right)=0 .
\]
Therefore, by the above change of variables and Tonelli's theorem
\[
\begin{aligned}
  &\lambda_d\left(
  \left\{\vk s\in\R^d:
  \int_K \ind{f(\vk t)+\vk s>\vk0}\,d\vk t=x
  \right\}
  \right)                                                     
   =
  \int_{\R^{d-1}}
  \int_{\R}
  \ind{
  \lambda_k\{\vk t\in K:g_{\vk q}(\vk t)>-r\}=x
  }
  \,dr\,d\vk q
  =0 .
\end{aligned}
\]
Restricting \(\vk s\) to \((0,\infty)^d\) preserves the nullity.

Finally, let
\[
  A_x:=
  \left\{\vk s\in(0,\infty)^d:
  \int_K \ind{f(\vk t)+\vk s>\vk0}\,d\vk t=x
  \right\}.
\]
We have shown that \(\lambda_d(A_x)=0\). Since \(h\) is measurable and finite almost everywhere, \(\ind{A_x}h=0\)
almost everywhere. Hence
\[
  \int_{(0,\infty)^d}\ind{\vk s \in A_x}h(\vk s)\,d\vk s=0 
\]
completing the proof.
\end{proof}
 
\begin{proof}[Proof of \Cref{contBx}]
Let
\[
  S_\infty
  =
  \int_{\R^k}
  \ind{
  \diag(\vk w)(\vk Y(\vk t)-V(\vk t)\vk w)+\vk E>\vk0
  }\,d\vk t ,
\]
where \(E_1,\ldots,E_d\) are independent unit exponential rvs,
independent of \(\vk Y\).

Fix \(x>0\) and let \(\Omega_0\) be the event on which the sample paths of
\(\vk Y\) are continuous and thus  \(\pk{\Omega_0}=1\). For
\(\omega\in\Omega_0\), define
\[
  f_\omega(\vk t)
  =
  \diag(\vk w)
  ( \vk Y(\vk t,\omega)-V(\vk t)\vk w),
  \qquad \vk t\in\R^k .
\]
The function \(f_\omega\) is continuous. Hence, by
\Cref{lem:levelset-null-shift}, applied with \(K=\R^k\) and
$h(\vk s)=\prod_{j=1}^d e^{-s_j}$,
we have
\[
  \int_{(0,\infty)^d}
  \ind{
  \int_{\R^k}
  \ind{f_\omega(\vk t)+\vk s>\vk0}\,d\vk t
  =x
  }
  \prod_{j=1}^d e^{-s_j}\,d\vk s
  =
  0 .
\]
Using the independence of \(\vk E\) and \(\vk Y\), it follows that
\[
\begin{aligned}
  \pk{S_\infty=x}
  &=
  \E*{
  \int_{(0,\infty)^d}
  \ind{
  \int_{\R^k}
  \ind{
  \diag(\vk w)(\vk Y(\vk t)-V(\vk t)\vk w)+\vk s>\vk0
  }\,d\vk t
  =x
  }
  \prod_{j=1}^d e^{-s_j}\,d\vk s
  } =0 
\end{aligned}
\]
which implies
$\pk{S_\infty=x}=0$,
$x>0$.

It remains to treat the point \(x=0\). Since
\(\vk Y(\vk0)=\vk0\), \(V(\vk0)=0\), \(\vk w>\vk0\), and
\(\vk E>\vk0\) almost surely,
\[
  \vk Y(\vk0)-V(\vk0)\vk w+\vk E/\vk w
  =
  \vk E/\vk w
  >
  \vk0
  \qquad\text{a.s.}
\]
By continuity, the integrand defining \(S_\infty\) is equal to one on a
neighbourhood of \(\vk0\), almost surely. Hence
$ S_\infty>0$ {a.s.}, \qq{and therefore}
$\pk{S_\infty=0}=0$.

Since \(\mathfrak F_{\vk Y,\vk w}\) is the distribution function of
\(S_\infty\), the preceding two displays show that
\(\mathfrak F_{\vk Y,\vk w}\) has no atoms on \([0,\infty)\). Hence
\(\mathfrak F_{\vk Y,\vk w}\) is continuous on \([0,\infty)\).

It remains to prove the continuity of
\(\mathcal B_{\vk Y,\vk w}\) under the assumption
\(\mathcal B_{\vk Y,\vk w}(0)<\infty\). Let \(\mu\) be the Stieltjes measure
induced by \(\mathfrak F_{\vk Y,\vk w}\) on the finite part of
\([0,\infty)\), and define
\[
  \nu(A)
  =
  \int_{A\cap(0,\infty)}\frac1y\,\mu(dy),
  \qquad A\in\mathcal B([0,\infty)).
\]
The assumption \(\mathcal B_{\vk Y,\vk w}(0)<\infty\) means precisely that
$\nu((0,\infty))<\infty$.
Moreover, for \(x>0\),
$\mathcal B_{\vk Y,\vk w}(x)
  =
  \nu([x,\infty))$,
with no ambiguity at the endpoint, because
\[
  \nu(\{x\})
  =
  \frac1x\mu(\{x\})
  =
  \frac1x\pk{S_\infty=x}
  =
  0 .
\]

Let \(x>0\) and let \(x_n\to x\), with \(x_n\ge0\).   
As \(\nu\) is finite, dominated convergence gives
\[
  \mathcal B_{\vk Y,\vk w}(x_n)
  =
  \nu([x_n,\infty))
  \to
  \nu([x,\infty))
  =
  \mathcal B_{\vk Y,\vk w}(x).
\]
Thus \(\mathcal B_{\vk Y,\vk w}\) is continuous on \((0,\infty)\).

At \(x=0\), if \(x_n\downarrow0\), then
\[
  \mathcal B_{\vk Y,\vk w}(x_n)
  =
  \nu([x_n,\infty))
  \uparrow
  \nu((0,\infty))
  =
  \mathcal B_{\vk Y,\vk w}(0)
\]
by monotone convergence. Since
\(\mathcal B_{\vk Y,\vk w}\) is nonincreasing, the same limit holds for any
sequence \(x_n\to0\) with \(x_n\ge0\). Therefore
\(\mathcal B_{\vk Y,\vk w}\) is continuous at \(0\). Together with
continuity on \((0,\infty)\), this proves continuity on \([0,\infty)\).
\end{proof}

\subsection{Proof of Theorem~\ref{T:stat}}\label{s.th1}

\begin{proof}
Throughout the proof, the Savage condition~\eqref{Savage},
Assumptions~\ref{I:B0}--\ref{I:B2}, and the conditions \(T_u>0\),
\eqref{ta}, and~\eqref{tb} remain in force. When
\(\limsup_{u\to\infty}T_u=\infty\), Assumption~\ref{I:B3} also remains
in force.
\qq{Recall notation}
$p_u=\pk{\vk X(\vk0)>u\vk b}$, $\eta(u)=\prod_{i=1}^k v_i(u)$,
and let \(F_u\) denote the distribution function of
\[
  \eta(u)L_u(T_u),
  \qquad
  L_u(T_u)=
  \int_{[0,T_u]^k}\ind{\vk X(\vk s)>u\vk b}\,d\vk s .
\]

By Fubini's theorem for every \(u>0\) and \(x\ge0\)
\begin{align}
\int_0^x y\,dF_u(y)
&=
\E*{
\eta(u)L_u(T_u)
\ind{\eta(u)L_u(T_u)\le x}
}
\nonumber\\
&=
\eta(u)\int_{[0,T_u]^k}
\pk{
\eta(u)L_u(T_u)\le x,\
\vk X(\vk t)>u\vk b
}
\,d\vk t .
\end{align}
By stationarity we have 
\[
\E{L_u(T_u)}=T_u^kp_u,
\qquad
\pk{\vk X(\vk t)>u\vk b}=p_u .
\]
Dividing the preceding identity by
\(\eta(u)\E{L_u(T_u)}\) therefore gives Berman's identity
\begin{equation}\label{eq:Berman-identity}
\frac{\int_0^x y\,dF_u(y)}
{\eta(u)\E{L_u(T_u)}}
=
T_u^{-k}
\int_{[0,T_u]^k}
\pk{
\eta(u)L_u(T_u)\le x
\mid
\vk X(\vk t)>u\vk b
}
\,d\vk t .
\end{equation}
For fixed \(\vk t\in\R^k\)
define
\begin{eqnarray}
  \vk Z_{u,\vk t}(\vk s)
  =
  u\left(
  \vk X(\vk t+\vk s/\vk v(u))
  -
  \RR(\vk s/\vk v(u))\Sigma^{-1}\vk X(\vk t)
  \right),
  \qquad \vk s\in\R^k .\label{reg.1}
\end{eqnarray}

By Gaussian regression, \(\vk Z_{u,\vk t}\) is independent of
\(\vk X(\vk t)\) with  covariance function \[
  K_u(\vk x,\vk y)
  =
  u^2\left[
  \RR\!\left(\frac{\vk x-\vk y}{\vk v(u)}\right)
  -
  \RR\!\left(\frac{\vk x}{\vk v(u)}\right)
  \Sigma^{-1}
  \RR\!\left(\frac{\vk y}{\vk v(u)}\right)^\top
  \right].
\]
For \(\vk q\in\R^k\), write
\[
A_u(\vk q)
=
\Sigma-\RR(\vk q/\vk v(u)).
\]
Since
\(\RR(\vk y/\vk v(u))^\top=\RR(-\vk y/\vk v(u))\),
\qq{we have}
\[
\begin{aligned}
&\RR(\vk x/\vk v(u))
\Sigma^{-1}
\RR(\vk y/\vk v(u))^\top
=
\Sigma-A_u(\vk x)-A_u(-\vk y)
+
A_u(\vk x)\Sigma^{-1}A_u(-\vk y)
\end{aligned}
\] implying 
\[
\begin{aligned}
K_u(\vk x,\vk y)
=
u^2\bigl[
A_u(\vk x)+A_u(-\vk y)-A_u(\vk x-\vk y)
\bigr]
-
u^2A_u(\vk x)\Sigma^{-1}A_u(-\vk y).
\end{aligned}
\]
By \eqref{e:stat-all-r}, the first term converges locally uniformly to
\[
V(\vk x)+V(-\vk y)-V(\vk x-\vk y),
\]
whereas \(A_u(\vk q)=O(u^{-2})\) locally uniformly, so the product term is
\(O(u^{-2})\). Hence
\[
K_u(\vk x,\vk y)
\rightarrow
V(\vk x)+V(-\vk y)-V(\vk x-\vk y)
=
R_V(\vk x,\vk y)
\]
locally uniformly in \((\vk x,\vk y)\) \zm{as $u\rightarrow\infty$}. This implies convergence of the
finite-dimensional distributions of \(\vk Z_{u,\vk t}\) to those of the
centered GRF \(\vk Y\) with covariance function \(R_V\). Moreover, uniformly for \(\vk s\)
in compact subsets of \(\R^k\), as $u\to\infty$
\[
  u^2\left(\mathcal I_d-\RR(\vk s/\vk v(u))\Sigma^{-1}\right)\vk b
  =
  u^2\left(\Sigma-\RR(\vk s/\vk v(u))\right)\vk w
  \to
  V(\vk s)\vk w .
\]

Fix \(m\ge1\) and \(\vk s_1,\ldots,\vk s_m\in\R^k\), and put
\[
\vk H_u\coloneqq u\bigl(\vk X(\vk t)-u\vk b\bigr).
\]
\qq{Representation (\ref{reg.1})} gives, for \(i=1,\ldots,m\),
\begin{align}
u\bigl(
\vk X(\vk t+\vk s_i/\vk v(u))-u\vk b
\bigr)
=
\vk Z_{u,\vk t}(\vk s_i)
+
\RR(\vk s_i/\vk v(u))\Sigma^{-1}\vk H_u
-
u^2\left[
\mathcal I_d-
\RR(\vk s_i/\vk v(u))\Sigma^{-1}
\right]\vk b .
\label{eq:conditional-regression-identity}
\end{align}
For every \(u\), the vector
\[
\bigl(
\vk Z_{u,\vk t}(\vk s_1),\ldots,
\vk Z_{u,\vk t}(\vk s_m)
\bigr)
\]
is independent of \(\vk H_u\). Since
\(\{\vk X(\vk t)>u\vk b\}=\{\vk H_u>\vk0\}\), this independence remains
valid under the conditional law, and the residual vector retains its
unconditional distribution. Moreover,
\cite{hashorva2007asymptotic} gives
\[
\vk H_u
\mid
\{\vk X(\vk t)>u\vk b\}
\xrightarrow{d}
\vk E/\vk w, \qquad u\to \infty
\]
where \(\vk E\) has independent unit exponential components. Combining this
with the finite-dimensional convergence of \(\vk Z_{u,\vk t}\), we obtain
the joint convergence
\[
\left(
(\vk Z_{u,\vk t}(\vk s_i))_{i=1}^m,\vk H_u
\right)
\Bigm|
\{\vk X(\vk t)>u\vk b\}
\Rightarrow
\left(
(\vk Y(\vk s_i))_{i=1}^m,\vk E/\vk w
\right),
\]
as $u\rightarrow\infty$ where \(\vk E\) is independent of \(\vk Y\).

Furthermore,  as $u\to \infty$
\[
\RR(\vk s_i/\vk v(u))\Sigma^{-1}\rightarrow\mathcal I_d
\]
and
\[
u^2\left[
\mathcal I_d-\RR(\vk s_i/\vk v(u))\Sigma^{-1}
\right]\vk b
\rightarrow
V(\vk s_i)\vk w .
\]
Thus \eqref{eq:conditional-regression-identity} and the continuous mapping
theorem yield
\[
\left(
u\bigl[
\vk X(\vk t+\vk s_i/\vk v(u))-u\vk b
\bigr]
\right)_{i=1}^m
\Bigm|
\{\vk X(\vk t)>u\vk b\}
\Rightarrow
\left(
\vk Y(\vk s_i)-V(\vk s_i)\vk w+\vk E/\vk w
\right)_{i=1}^m .
\]
Each coordinate of the limiting vector has no atom at zero: conditional on
\(\vk Y\), it is a constant plus \(E_j/w_j\), which has a density. Hence the
boundary of the positive orthant has probability zero, and therefore
\begin{align}
  &\lim_{u\to\infty}
  \pk{
  \vk X(\vk t+\vk s_i/\vk v(u))>u\vk b,\ i=1,\ldots,m
  \mid
  \vk X(\vk t)>u\vk b
  }
  \nonumber\\
  &\qquad =
  \pk{
  \vk Y(\vk s_i)-V(\vk s_i)\vk w+\vk E/\vk w>\vk0,\
  i=1,\ldots,m
  } .
  \label{eq:conditional-fidi}
\end{align}
For \(\lambda>0\), define the local sojourn time
\[
  L_u(T_u,\vk t,\lambda)
  =
  \int_{[0,T_u]^k
  \cap
  \{\vk s:\ |s_i-t_i|\le \lambda/v_i(u),\ i=1,\ldots,k\}}
  \ind{\vk X(\vk s)>u\vk b}\,d\vk s 
\]
and put
\[
  S_\lambda
  =
  \int_{[-\lambda,\lambda]^k}
  \ind{
  \vk Y(\vk s)-V(\vk s)\vk w+\vk E/\vk w>\vk0
  }\,d\vk s .
\]

We pass from finite-dimensional convergence to the local sojourn functional.
Suppose first that the local rectangle is contained in \([0,T_u]^k\) and hence   
\[
\eta(u)L_u(T_u,\vk t,\lambda)
=
\int_{[-\lambda,\lambda]^k}
\ind{\vk X(\vk t+\vk s/\vk v(u))>u\vk b}\,d\vk s .
\]
For every integer \(n\ge1\), conditioning on
\(\{\vk X(\vk t)>u\vk b\}\) and using Fubini's theorem gives
\[
\begin{aligned}
&\E*{
\left(
\eta(u)L_u(T_u,\vk t,\lambda)
\right)^n
\mid
\vk X(\vk t)>u\vk b
}  \\
&\quad =
\int_{([-\lambda,\lambda]^k)^n}
\pk{
\vk X(\vk t+\vk s_j/\vk v(u))>u\vk b,\ 1\le j\le n
\mid
\vk X(\vk t)>u\vk b
}
\,d\vk s_1\cdots d\vk s_n .
\end{aligned}
\]
By \eqref{eq:conditional-fidi}, the integrand converges pointwise to
\[
\pk{
\vk Y(\vk s_j)-V(\vk s_j)\vk w+\vk E/\vk w>\vk0,\ 1\le j\le n
}.
\]
Since the integrand is bounded by \(1\), dominated convergence yields
convergence of all moments to the moments of \(S_\lambda\). Both random
variables are bounded by \((2\lambda)^k\). Hence, moment convergence implies
\[
\eta(u)L_u(T_u,\vk t,\lambda)
\mid
\{\vk X(\vk t)>u\vk b\}
\Rightarrow
S_\lambda .
\]
By stationarity, the conditional law on the left-hand side is independent of
$\vk t$ whenever the local rectangle is untruncated. Hence, this convergence
is uniform over $\vk t\notin B_{u,\lambda}$, where $B_{u,\lambda}$ is defined
below.

For fixed \(\lambda\) the boundary contribution is negligible. Indeed let
\[
  B_{u,\lambda}
  =
  \left\{
  \vk t\in[0,T_u]^k:
  \exists i,\ 
  t_i<\lambda/v_i(u)
  \text{ or }
  T_u-t_i<\lambda/v_i(u)
  \right\}.
\]
Then by the union bound
\[
  \frac{\lambda_k(B_{u,\lambda})}{T_u^k}
  \le
  2\lambda\sum_{i=1}^k\frac{1}{T_uv_i(u)}
  \to0, \qquad u\to\infty
\]
because \(T_uv_i(u)\to\infty\) for \(i=1,\ldots,k\). Therefore, for every
continuity point \(x\) of the distribution of \(S_\lambda\),
\begin{equation}\label{eq:local-avg-limit}
  \lim_{u\to\infty}
  T_u^{-k}
  \int_{[0,T_u]^k}
  \pk{
  \eta(u)L_u(T_u,\vk t,\lambda)\le x
  \mid
  \vk X(\vk t)>u\vk b
  }\,d\vk t
  =
  \pk{S_\lambda\le x}.
\end{equation}
Conditioning on $\vk Y$, \Cref{lem:levelset-null-shift}, applied with
$K=[-\lambda,\lambda]^k$ and the density
\[
   \prod_{j=1}^d w_j e^{-w_jz_j}\ind{\vk z>\vk0}
\]
of $\vk E/\vk w$, shows that $S_\lambda$ has no atoms in
$(0,(2\lambda)^k)$. Moreover, at $\vk s=\vk0$ the vector inside the
indicator equals $\vk E/\vk w>\vk0$ almost surely; continuity therefore
implies $S_\lambda>0$ almost surely. Thus, for every fixed \(x\ge0\), there
exists \(\lambda_x<\infty\) such that \(x\) is a continuity point of the
distribution of \(S_\lambda\) for every \(\lambda\ge\lambda_x\).
Consequently, for every such \(\lambda\), the inner \(u\)-limit in
\eqref{eq:local-avg-limit} exists and equals \(\pk{S_\lambda\le x}\).

Let
\[
  S_\infty
  =
  \int_{\R^k}
  \ind{
  \vk Y(\vk s)-V(\vk s)\vk w+\vk E/\vk w>\vk0
  }\,d\vk s .
\]
Since every \(w_i>0\)
\[
\vk Y(\vk s)-V(\vk s)\vk w+\vk E/\vk w>\vk0
\quad\Longleftrightarrow\quad
\diag(\vk w)\bigl(\vk Y(\vk s)-V(\vk s)\vk w\bigr)+\vk E>\vk0.
\]
Hence \(S_\infty=\mathcal J(\vk Y,\vk w)\) almost surely, and its
distribution function is \(\mathfrak F_{\vk Y,\vk w}\).

Moreover, \(S_\lambda\uparrow S_\infty\), and therefore
\[
\{S_\lambda\le x\}\downarrow\{S_\infty\le x\}.
\]
Continuity of probability from above gives, for every \(x\ge0\),
\[
\pk{S_\lambda\le x}
\rightarrow
\pk{S_\infty\le x}
=
\mathfrak F_{\vk Y,\vk w}(x).
\]
Hence, for every fixed \(x\ge0\),
\[
  \lim_{\lambda\to\infty}
  \lim_{u\to\infty}
  T_u^{-k}
  \int_{[0,T_u]^k}
  \pk{
  \eta(u)L_u(T_u,\vk t,\lambda)\le x
  \mid
  \vk X(\vk t)>u\vk b
  }\,d\vk t
  =
  \mathfrak F_{\vk Y,\vk w}(x)
\]
which verifies the local-limit hypothesis of
\Cref{lem:separation}.

It remains to verify the separation condition in \Cref{lem:separation},
namely
\begin{equation}\label{eq:sep-to-prove}
  \lim_{\lambda\to\infty}
  \limsup_{u\to\infty}
  \frac{
  \eta(u)
  \int_{[-T_u,T_u]^k
  \cap
  \{|\vk h|\not\le \lambda/\vk v(u)\}}
  \pk{\vk X(\vk h)>u\vk b,\ \vk X(\vk0)>u\vk b}\,d\vk h
  }
  {p_u}
  =0 .
\end{equation}
Here
\[
  |\vk h|\le \lambda/\vk v(u)
  \quad\text{means}\quad
  |h_i|\le \lambda/v_i(u),\qquad i=1,\ldots,k .
\]

For \(\vk h\in\R^k\), define
\[
  \zeta(\vk h)
  =
  \vk X(\vk h)-\RR(\vk h)\Sigma^{-1}\vk X(\vk0),
  \qquad
  \chi(\vk h)
  =
  \vk w^\top\zeta(\vk h),
\]
and
$L(\vk h)
  =
  \vk w^\top(\Sigma-\RR(\vk h))\vk w$.
The random vector \(\zeta(\vk h)\) is independent of \(\vk X(\vk0)\), and
\(\chi(\vk h)\) is centered Gaussian with variance
\[
  b^2(\vk h)
  \coloneqq
  \Var(\chi(\vk h))
  =
  \vk w^\top
  \left[
  \Sigma-\RR(\vk h)\Sigma^{-1}\RR(\vk h)^\top
  \right]
  \vk w .
\] 
Conditioning on \(\vk X(\vk0)\), using the change of variables
\(\vk y=u\vk b+\vk z/u\), \(\vk z>\vk0\), 
there exists \(C>0\) such that for all large \(u\)
\begin{align}
  \pk{\vk X(\vk h)>u\vk b,\ \vk X(\vk0)>u\vk b}
  %\nonumber\\
  %&\quad
  \le
  C p_u
  \int_{\vk z>\vk0}
  \pk{
  \chi(\vk h)
  >
  uL(\vk h)
  -
  u^{-1}\vk w^\top\RR(\vk h)\Sigma^{-1}\vk z
  }
  \prod_{j=1}^d w_j e^{-w_jz_j}\,d\vk z .
\end{align}
In the preceding integral, set
$
e_j=w_jz_j$,
$j=1,\ldots,d$, so

\[
\prod_{j=1}^d w_je^{-w_jz_j}\,d\vk z
=
\prod_{j=1}^d e^{-e_j}\,d\vk e,
\]
which is the joint density of independent unit exponential rvs
\(E_1,\ldots,E_d\). Since the entries of \(\RR(\vk h)\) are uniformly bounded,
there is \(c_*>0\), independent of \(\vk h\), such that
\[
\left|
\vk w^\top\RR(\vk h)\Sigma^{-1}\vk z
\right|
\le
c_*\sum_{j=1}^d z_j
\le
c_*\left(\max_{1\le j\le d}w_j^{-1}\right)
\sum_{j=1}^d e_j .
\]
Set
$c_0=c_*\max_{1\le j\le d}w_j^{-1}$,
$\mathcal E=\sum_{j=1}^dE_j$,
and take \(\vk E\) independent of \(\chi(\vk h)\). It follows that
\begin{equation}\label{eq:joint-tail-bound}
  \pk{\vk X(\vk h)>u\vk b,\ \vk X(\vk0)>u\vk b}
  \le
  C p_u\,
  \pk{
  \chi(\vk h)
  >
  uL(\vk h)-\frac{c_0}{u}\mathcal E
  },
\end{equation}
where \(E_1,\ldots,E_d\) are independent unit exponential rvs.

The crucial bound is
\begin{equation}\label{eq:b2-bound}
  b^2(\vk h)=2L(\vk h)
  -
  \vk w^\top
  (\Sigma-\RR(\vk h))
  \Sigma^{-1}
  (\Sigma-\RR(\vk h))^\top
  \vk w  \le 2L(\vk h).
\end{equation} 

Fix \(p_0\in(0,1)\), and let \(\theta_{p_0}\) be the corresponding constant
in \eqref{e:scalar-potter-lower}. Choose \(\rho>0\)
so small that   
$\sqrt{k}\rho<\theta_{p_0}$.
Split the domain in \eqref{eq:sep-to-prove} into
\[
  D_{u,\lambda}^{(0)}
  =
  [-T_u,T_u]^k
  \cap
  \{|\vk h|\not\le \lambda/\vk v(u)\}
  \cap
  \{\|\vk h\|_\infty<\rho\}
\]
and
\[
  D_{u,\lambda}^{(1)}
  =
  [-T_u,T_u]^k
  \cap
  \{|\vk h|\not\le \lambda/\vk v(u)\}
  \cap
  \{\|\vk h\|_\infty\ge\rho\}.
\]

We first treat \(D_{u,\lambda}^{(1)}\).
If \(T_*=\limsup_{u\to\infty} T_u<\infty\), then, for all large \(u\),
\(D_{u,\lambda}^{(1)}\subset\{\rho\le\|\vk h\|_\infty\le T_*+1\}\).
Hence, by Assumption~\ref{I:B1} and the continuity of \(\RR\)
\[
   \inf_{\rho\le\|\vk h\|_\infty\le T_*+1} L(\vk h)>0 .
\]
If \(T_*=\infty\), the same conclusion follows from
Assumption~\ref{I:B3}. Thus, in either case, there exists \(m_\rho>0\)
such that
\[
  L(\vk h)\ge m_\rho,
  \qquad
  \vk h\in D_{u,\lambda}^{(1)}
\]
for all large \(u\). Moreover, the entries of the covariance function \(\RR\)
are uniformly bounded by Cauchy's inequality, and hence
\[
\bar b^2:=\sup_{\vk h\in\R^k}b^2(\vk h)<\infty .
\]
It follows that, for some constants \(C_1,c_1>0\)
\[
  \pk{
  \chi(\vk h)
  >
  uL(\vk h)-\frac{c_0}{u}\mathcal E
  }
  \le
  C_1e^{-c_1u^2},
  \qquad
  \vk h\in D_{u,\lambda}^{(1)} .
\]
To see this, observe that on the event
$\{  \mathcal E \le \frac{u^2m_\rho}{2c_0}\}$,
the threshold is at least \(um_\rho/2\), giving a Gaussian tail of order
\(e^{-cu^2}\) whenever \(b^2(\vk h)>0\), uniformly in \(\vk h\). If
\(b^2(\vk h)=0\), then \(\chi(\vk h)=0\) almost surely, so the corresponding
positive-threshold probability is zero. The complementary event is also
bounded by \(Ce^{-c'u^2}\), because \(\mathcal E\) has a Gamma distribution.
Consequently,
\[
\limsup_{u\to\infty}
  \eta(u)
  \int_{D_{u,\lambda}^{(1)}}
  \pk{
  \chi(\vk h)
  >
  uL(\vk h)-\frac{c_0}{u}\mathcal E
  }\,d\vk h
  \le\limsup_{u\to\infty}
  C_1\eta(u)(2T_u)^k e^{-c_1u^2}
  =0 .
\]
Here we used
\[
   \log\bigl(\eta(u)T_u^k\bigr)=o(u^2),
\]
which follows from the regular variation of the \(v_i\)'s together with the assumption 
\(\log T_u=o(u^2)\).

It remains to treat \(D_{u,\lambda}^{(0)}\). Make the componentwise change of
variables
\[
  r_i=h_i v_i(u),
  \qquad i=1,\ldots,k .
\]
  The local contribution is bounded by
\[
  \int_{A_{u,\lambda}}
  \pk{
  u\chi(\vk r/\vk v(u))
  >
  u^2L(\vk r/\vk v(u))-c_0\mathcal E
  }\,d\vk r ,
\]
where
\[
  A_{u,\lambda}
  =
  \left\{
  \vk r:
  |\vk r|\not\le \lambda\vk1,\ 
  \|\vk r/\vk v(u)\|_\infty<\rho,\ 
  |r_i|\le T_uv_i(u),\ i=1,\ldots,k
  \right\}.
\]
For \(\vk r\in A_{u,\lambda}\),
\[
  \|\vk r/\vk v(u)\|
  \le
  \sqrt{k}\|\vk r/\vk v(u)\|_\infty
  <
  \sqrt{k}\rho
  <
  \theta_{p_0}.
\]
After increasing \(\lambda\), if necessary, the condition
$|\vk r|\not\le\lambda\vk1$
implies
$
  [\vk r]_{\vk\alpha}\ge1$, 
  and hence  uniformly on \(A_{u,\lambda}\),
\begin{equation}\label{eq:L-lower-new}
u^2L(\vk r/\vk v(u))
\ge
c_{p_0}[\vk r]_{\vk\alpha}^{p_0},
\end{equation}
by \eqref{e:scalar-potter-lower}. Put next 
\[
G_{u,\vk r}\coloneqq u\chi(\vk r/\vk v(u)),
\qquad
a_{u,\vk r}\coloneqq u^2L(\vk r/\vk v(u))>0.
\]
In view of  \eqref{eq:b2-bound} we have 
\[
   \Var(G_{u,\vk r})
   =
   u^2b^2(\vk r/\vk v(u))
   \le
   2a_{u,\vk r}
\]
and hence 
\[
\begin{aligned}
&\pk{
G_{u,\vk r}
>
a_{u,\vk r}-c_0\mathcal E
}                                                \le
\pk{
G_{u,\vk r}
>
\frac12a_{u,\vk r}
}
+
\pk{
\mathcal E>
\frac{a_{u,\vk r}}{2c_0}
}.
\end{aligned}
\]
If \(\Var(G_{u,\vk r})=0\), then
$\pk{G_{u,\vk r}>a_{u,\vk r}/2}=0$.
Otherwise, the Gaussian tail estimate and the variance bound give
\[
\pk{G_{u,\vk r}>a_{u,\vk r}/2}
\le
\exp\left\{-\frac{a_{u,\vk r}^2}
{8\Var(G_{u,\vk r})}\right\}
\le
e^{-a_{u,\vk r}/16}.
\]
Thus the same estimate holds in both the degenerate and nondegenerate cases.
The second term is bounded in the same way, up to a multiplicative constant,
since \(\mathcal E\) has a Gamma distribution. Hence, by
\eqref{eq:L-lower-new},
there exist constants \(C_2,c_2>0\) such that
\begin{equation}\label{eq:local-exp-bound-new}
\pk{
G_{u,\vk r}
>
u^2L(\vk r/\vk v(u))-c_0\mathcal E
}
\le
C_2\exp\{-c_2[\vk r]_{\vk\alpha}^{p_0}\},
\qquad \vk r\in A_{u,\lambda}
\end{equation}
and hence 
\[
\begin{aligned}
  &\limsup_{u\to\infty}
  \eta(u)
  \int_{D_{u,\lambda}^{(0)}}
  \pk{
  \chi(\vk h)
  >
  uL(\vk h)-\frac{c_0}{u}\mathcal E
  }\,d\vk h                                     \le
  C_2
  \int_{\{|\vk r|\not\le \lambda\vk1\}}
  \exp\{-c_2[\vk r]_{\vk\alpha}^{p_0}\}\,d\vk r .
\end{aligned}
\]
Indeed, since \(p_0\in(0,1)\), \qq{we have}
\[
[\vk r]_{\vk\alpha}^{p_0}
=
\left(\sum_{i=1}^k|r_i|^{\alpha_i}\right)^{p_0}
\ge
k^{p_0-1}\sum_{i=1}^k|r_i|^{\alpha_ip_0}
\]
implying 
\[
e^{-c_2[\vk r]_{\vk\alpha}^{p_0}}
\le
\prod_{i=1}^k
e^{-c_2k^{p_0-1}|r_i|^{\alpha_ip_0}},
\]
which is integrable on \(\R^k\). Dominated convergence therefore shows that
the integral over \(\{|\vk r|\not\le\lambda\vk1\}\) tends to zero as
\(\lambda\to\infty\). Combining this with the estimate on
\(D_{u,\lambda}^{(1)}\), and then using
\eqref{eq:joint-tail-bound}, proves \eqref{eq:sep-to-prove}.

The separation estimate was carried out on the full cube
\([-T_u,T_u]^k\). Thus, all \(2^k\) orthants are treated simultaneously, and no
separate argument for positive and negative lags is needed.

By \Cref{lem:separation}, we obtain
\[
  \lim_{u\to\infty}
  T_u^{-k}
  \int_{[0,T_u]^k}
  \pk{
  \eta(u)L_u(T_u)\le x
  \mid
  \vk X(\vk t)>u\vk b
  }\,d\vk t
  =
  \mathfrak F_{\vk Y,\vk w}(x)
\]
at every continuity point \(x\ge0\) of
\(\mathfrak F_{\vk Y,\vk w}\). Together with
\eqref{eq:Berman-identity}, this gives
\[
  \lim_{u\to\infty}
  \frac{\int_0^x y\,dF_u(y)}
  {\eta(u)\E{L_u(T_u)}}
  =
  \mathfrak F_{\vk Y,\vk w}(x)
\]
at every continuity point of \(\mathfrak F_{\vk Y,\vk w}\). By
\Cref{contBx}, the distribution function \(\mathfrak F_{\vk Y,\vk w}\) is
continuous on \([0,\infty)\). Hence, the convergence holds for every
\(x\ge0\).
\end{proof}

\subsection{Proof of Theorem~\ref{T:doublesum}}
\label{s.th2}
Throughout the proof, all assumptions inherited from \Cref{T:stat}, including
the Savage condition remain in force.

For a Borel set \(A\subset\R^k\) define
\[
L_u^*(A)
\coloneq
\eta(u)\int_A\ind{\vk X(\vk t)>u\vk b}\,d\vk t 
\] 
and \qq{let $D_u\coloneq[0,T_u]^k.$}
%hence 
%\[
%L_u^*(D_u)
%=
%\eta(u)\int_{[0,T_u]^k}
%$\ind{\vk X(\vk t)>u\vk b}\,d\vk t , \qquad D_u\coloneq[0,T_u]^k.
%$\]
By stationarity, we have further 
\[
\E{L_u(T_u)}
=
T_u^k p_u,
\qquad
\E{L_u^*(D_u)}
=
T_u^k\eta(u)p_u .
\]

For \(S>0\), define
\begin{equation}
\label{eq:BS-def}
\mathcal B_S(x)
\coloneq
\int_{\R^d}
\pk*{
\int_{[0,S]^k}
\ind{
\operatorname{diag}(\vk w)
\bigl[\vk Y(\vk t)-V(\vk t)\vk w\bigr]
>
\vk z
}
\,d\vk t
>x
}
e^{\vk z^\top\vk1}\,d\vk z ,
\qquad x\ge0.
\end{equation}
 
Let \(F_u\) denote the distribution function of \(L_u^*(D_u)\). In view of 
\Cref{T:stat}, for every \(x\ge0\) we have 
\[
\frac{\int_0^x y\,dF_u(y)}
     {\int_0^\infty y\,dF_u(y)}
=
\frac{\int_0^x y\,dF_u(y)}
     {T_u^k\eta(u)p_u}
\rightarrow
\mathfrak F_{\vk Y,\vk w}(x), \qquad  u\to\infty.
\]
By \Cref{contBx}, the distribution function
\(\mathfrak F_{\vk Y,\vk w}\) is continuous on \([0,\infty)\). Put
\[
   m_u=T_u^k\eta(u)p_u,
   \qquad
   \mu_u(dy)=\frac{y\,dF_u(y)}{m_u}.
\]
The preceding convergence states that the probability measures $\mu_u$
converge weakly on the compactified half-line $[0,\infty]$ to the law of
$\mathcal J(\vk Y,\vk w)$, whose finite-part distribution function is
$\mathfrak F_{\vk Y,\vk w}$.

Moreover, Savage's asymptotics, the regular variation of \(\eta\), and
\eqref{tb} give, with \(a=\vk b^\top\Sigma^{-1}\vk b>0\),
\[
   \log m_u=-\frac12au^2+o(u^2)\to-\infty.
\]
Thus \(m_u\to0\), which is the normalization in
\cite[Lem~1.2.1]{Berman92}. The passage to tail probabilities can also be
seen directly: for \(x>0\), the identity
\[
   \frac{\pk{L_u^*(D_u)>x}}{m_u}
   =
   \int_{(x,\infty]}\frac1y\,\mu_u(dy),
   \qquad \frac1\infty:=0,
\]
has a bounded integrand whose only discontinuity on $[0,\infty]$ is at $x$.
The limiting law has no atom there, so weak convergence gives
\begin{equation}
\label{eq:Berman-positive-x}
\lim_{u\to\infty}
\frac{\pk{L_u^*(D_u)>x}}
     {T_u^k\eta(u)p_u}
=
\int_x^\infty \frac1y\,d\mathfrak F_{\vk Y,\vk w}(y)
\eqqcolon
\mathcal B_{\vk Y,\vk w}(x).
\end{equation}
Any possible mass of the limiting law at \(\infty\) contributes zero to the
integral in \eqref{eq:Berman-positive-x}; thus this step does not presuppose
that \(\mathcal J(\vk Y,\vk w)\) is almost surely finite.
For \(x>0\), this integral is finite since
\[
\mathcal B_{\vk Y,\vk w}(x)
\le
\frac1x .
\]
Moreover, since \(\mathfrak F_{\vk Y,\vk w}\) is continuous, the function
\(x\mapsto\mathcal B_{\vk Y,\vk w}(x)\) is continuous on \((0,\infty)\).

\medskip
\noindent Next fix 
 \(S>1\). For \(\vk l\in\mathbb{Z}^k\), set
\[
A_{\vk l}(S)
=
[S\vk l,S(\vk l+\vk1)]
=
\bigtimes_{i=1}^k [Sl_i,S(l_i+1)],
\qquad 
\Delta_{\vk l}
=
A_{\vk l}(S)/\vk v(u)
=
\bigtimes_{i=1}^k
\left[
\frac{Sl_i}{v_i(u)},
\frac{S(l_i+1)}{v_i(u)}
\right]
\]
and let
\[
N_{u,i}(T_u,S)
\coloneq
\left\lfloor\frac{T_uv_i(u)}{S}\right\rfloor,
\qquad
\vk N_u(T_u,S)
=
(N_{u,1}(T_u,S),\ldots,N_{u,k}(T_u,S)).
\]
Since \( \limit{u}T_uv_i(u)=\infty\), we have as \(u\to\infty\)
\begin{equation}
\label{eq:block-number}
\frac{\prod_{i=1}^k N_{u,i}(T_u,S)}
     {T_u^k\eta(u)}
\to
\frac1{S^k},
\qquad
\frac{\prod_{i=1}^k (N_{u,i}(T_u,S)+1)}
     {T_u^k\eta(u)}
\to
\frac1{S^k}.
\end{equation}
\qq{Moreover, we have
\[
D_u^-
\coloneqq
\bigcup_{\vk0\le\vk l\le\vk N_u(T_u,S)-\vk1}
\Delta_{\vk l}\subseteq D_u 
\subseteq
\bigcup_{\vk0\le\vk l\le\vk N_u(T_u,S)}
\Delta_{\vk l}=:D_u^+.
\]
}
For two blocks define
\[
q_{\vk i,\vk j}(u)
\coloneqq
\pk*{
\exists\,\vk s\in\Delta_{\vk i},\
\exists\,\vk t\in\Delta_{\vk j}:
\vk X(\vk s)>u\vk b,\ \vk X(\vk t)>u\vk b
}.
\]
Set next
\[
  E_{\vk l}(x)
  =
  \{L_u^*(\Delta_{\vk l})>x\}.
\]
For the inner blocks,
\(\bigcup_{\vk l}E_{\vk l}(x)\subset\{L_u^*(D_u)>x\}\), while
\(E_{\vk i}(x)\cap E_{\vk j}(x)\) is contained in the event defining
\(q_{\vk i,\vk j}(u)\). Hence Bonferroni's inequality applied to the
events \(E_{\vk l}(x)\) yields 
\begin{eqnarray}
I_1(u,S)
\le
\pk{L_u^*(D_u)>x}
\le
I_2(u,S),\label{Bonf}    
\end{eqnarray}

where
\[
I_1(u,S)
=
\sum_{\vk0\le\vk l\le\vk N_u(T_u,S)-\vk1}
\pk{L_u^*(\Delta_{\vk l})>x}
-
\sum_{\substack{
\vk0\le\vk i,\vk j\le\vk N_u(T_u,S)-\vk1\\
\vk i\ne\vk j
}}
q_{\vk i,\vk j}(u),
\]
and
\[
I_2(u,S)
=
\sum_{\vk0\le\vk l\le\vk N_u(T_u,S)}
\pk{L_u^*(\Delta_{\vk l})>x}
+
\sum_{\substack{
\vk0\le\vk i,\vk j\le\vk N_u(T_u,S)\\
\vk i\ne\vk j
}}
q_{\vk i,\vk j}(u).
\]
For the upper bound, if the total occupation over the outer block union
exceeds \(x\), but no single block contributes more than \(x\), then at
least two blocks have positive occupation, which is covered by the
double-sum term. Distinct blocks meet only on Lebesgue-null boundaries, so
occupation is additive over the unions. The sums over
\(\vk i\ne\vk j\) are ordered; this only enlarges the pair-error term in the
upper bound and weakens the lower bound.

\medskip
\noindent
\textit{Single sum.}
By stationarity, for all $u>0$
\[
\sum_{\vk0\le\vk l\le\vk N_u(T_u,S)}
\pk{L_u^*(\Delta_{\vk l})>x}
=
\prod_{i=1}^k (N_{u,i}(T_u,S)+1)
\pk{L_u^*(\Delta_{\vk0})>x},
\] 
where   
$L_u^*(\Delta_{\vk0})
=
\int_{[0,S]^k}
\ind{\vk X(\vk t/\vk v(u))>u\vk b}\,d\vk t$.

We claim that, for fixed \(S>0\)
\[
\lim_{u\to\infty}
\frac{\pk{L_u^*(\Delta_{\vk0})>x}}{p_u}
=
\mathcal B_S(x)
\]
at every continuity point \(x\in(0,S^k)\) of \(x\mapsto\mathcal B_S(x)\),
and also at \(x=0\).
We apply Lemma~\ref{lemma:Pickands} with the singleton index set
\(Q_u=\{\ast\}\)
\[
  E=[0,S]^k,\qquad
  \bm t_0=\vk0,\qquad
  \bm Z_{u,\ast}(\vk t)=\vk X(\vk t/\vk v(u))
\]
and with limiting field \(\bm Z=\vk Y\). The Savage condition
\(\vk w=\Sigma^{-1}\vk b>\vk0\) implies that the solution of
\(\Pi_\Sigma(\vk b)\) is \(\widetilde{\bm b}=\vk b\) and that the
determining set is \(I=\{1,\ldots,d\}\).
Condition~\textup{\ref{C1}} holds because \(\Sigma\) is positive definite and
\[
  \Sigma_{u,\ast}=\Var(\bm Z_{u,\ast}(\vk0))=\Sigma .
\]
For \(\vk q\in \qq{[-S,S]^k=:\tilde{E}}\) 
put
$  G_u(\vk q)
  :=
  u^2\bigl(\Sigma-\RR(\vk q/\vk v(u))\bigr)$.
By \eqref{e:stat-all-r}, \(G_u\to V\) uniformly on the compact set \(\qq{\tilde{E}}\) as $u\to \infty$.
Consequently, as $u\to \infty$
\[
  \sup_{\vk t\in E}
  u\left\|\Sigma-\RR(\vk t/\vk v(u))\right\|
  =
  \frac1u\sup_{\vk t\in E}\|G_u(\vk t)\|
  \rightarrow0,
\]
which is the first assertion in \textup{\ref{C2}}. Moreover, uniformly for
\(\vk t\in E\) as $u\to \infty$
\[
\begin{aligned}
&u^2\left[
\mathcal I_d-\RR(\vk t/\vk v(u))\Sigma^{-1}
\right]\vk b 
=
G_u(\vk t)\vk w
\rightarrow
V(\vk t)\vk w .
\end{aligned}
\]
By stationarity 
$R_{u,\ast}(\vk s,\vk t)
  =
  \RR((\vk s-\vk t)/\vk v(u))$.
Hence, uniformly for \(\vk s,\vk t\in E\)
\[
\begin{aligned}
&u^2\Bigl[
\RR((\vk s-\vk t)/\vk v(u))
-\RR(\vk s/\vk v(u))\Sigma^{-1}
 \RR(-\vk t/\vk v(u))
\Bigr]  =
G_u(\vk s)+G_u(-\vk t)-G_u(\vk s-\vk t)
-u^{-2}G_u(\vk s)\Sigma^{-1}G_u(-\vk t) \\
&\quad\rightarrow
V(\vk s)+V(-\vk t)-V(\vk s-\vk t)
=R_V(\vk s,\vk t), \qquad u \to \infty.
\end{aligned}
\]
Here the product term tends to zero uniformly because \(G_u\) is uniformly
bounded on \(\qq{\tilde{E}}\).
Thus \textup{\ref{C2}} holds with
\(\bm f(\vk t)=V(\vk t)\vk w\) and \(K=R_V\).

Finally, since \(\qq{\tilde{E}}\) is compact, choose \(R<\infty\) such that
\([\vk q]_{\vk\alpha}\le R\) for \(\vk q\in \qq{\tilde{E}}\).
Condition~\textup{\ref{C3}} follows from the local matrix bound
\eqref{e:local-holder-matrix}, with
$\gamma_i=\beta_i\in(0,2]$, $i=1,\ldots,k$.
Indeed, stationarity and \(\RR(-\vk q)=\RR(\vk q)^\top\) give
\[
\begin{aligned}
&u^2\E*{\left\|
\vk X(\vk t/\vk v(u))-\vk X(\vk s/\vk v(u))
\right\|^2} \\
&\quad=
u^2\operatorname{tr}\left[
2\Sigma-\RR((\vk t-\vk s)/\vk v(u))
-\RR((\vk t-\vk s)/\vk v(u))^\top
\right]
\le C[\vk t-\vk s]_{\vk\beta}.
\end{aligned}
\]
Because \(\Sigma_{u,\ast}=\Sigma\), the exponential correction in
Lemma~\ref{lemma:Pickands} is identically one. The limiting field is
\(\vk Y\), its covariance kernel is \(R_V\), and the limiting drift is
\(V(\vk t)\vk w\); hence the limiting constant is precisely
\(\mathcal B_S(x)\). Lemma~\ref{lemma:Pickands} therefore gives the
claimed convergence at every continuity point and at \(x=0\).

The assertion at \(x=0\) follows directly from
Lemma~\ref{lemma:Pickands}; by sample-path continuity,
\[
\{L_u^*(\Delta_{\vk0})>0\}
=
\{\exists\,\vk t\in[0,S]^k:
\vk X(\vk t/\vk v(u))>u\vk b\}.
\]
Indeed, a strict componentwise exceedance persists on a relative
neighbourhood of \(\vk t\) in the cube and therefore gives positive
\(k\)-dimensional occupation; the converse is immediate.

Therefore, using \eqref{eq:block-number} we obtain
\begin{equation}
\label{eq:single-sum}
\lim_{u\to\infty}
\frac{1}{T_u^k\eta(u)p_u}
\sum_{\vk0\le\vk l\le\vk N_u(T_u,S)}
\pk{L_u^*(\Delta_{\vk l})>x}
=
\frac{\mathcal B_S(x)}{S^k}.
\end{equation}
The same limit holds for the \qq{indices} 
\(\vk0\le\vk l\le\vk N_u(T_u,S)-\vk1\).

\qq{{\it Double sum.}}
\qq{We show that the double sums in $I_1(u,S)$ and $I_2(u,S)$,
normalised as in \eqref{eq:single-sum}, are negligible after first letting
$u\to\infty$ and then $S\to\infty$.}
By stationarity, the number of pairs of blocks with fixed displacement
\(\vk l\) is bounded by
$\prod_{i=1}^k (N_{u,i}(T_u,S)+1).$
\qq{In view of (\ref{Bonf}), with
$q_{\vk l}(u)\coloneq q_{\vk l,\vk0}(u)$,
it suffices to show that}
\begin{equation}
\label{eq:double-target}
\lim_{S\to\infty}
\limsup_{u\to\infty}
\frac1{S^kp_u}
\sum_{\substack{
-\vk N_u(T_u,S)\le\vk l\le\vk N_u(T_u,S)\\
\vk l\ne\vk0
}}
q_{\vk l}(u)
=0.
\end{equation}

Fix \(\delta\in(0,1)\), choose
\(\theta\in(0,\theta_0(\delta))\) as in \Cref{lem.double}, and put
$
N_{u,i}(\theta,S)
\coloneq
\left\lceil\frac{\theta v_i(u)}S\right\rceil$.
Split the sum in \eqref{eq:double-target} into
$
\Sigma_{\mathrm{far}}+\Sigma_{\mathrm{nn}}+\Sigma_{\mathrm{nb}},
$
where
\[
\Sigma_{\mathrm{far}}
=
\sum_{\substack{
-\vk N_u(T_u,S)\le\vk l\le\vk N_u(T_u,S)\\
\exists i:\ |l_i|>N_{u,i}(\theta,S)
}}
q_{\vk l}(u),
\qquad 
\Sigma_{\mathrm{nn}}
=
\sum_{\substack{
-\vk N_u(T_u,S)\le\vk l\le\vk N_u(T_u,S)\\
|l_i|\le N_{u,i}(\theta,S),\ i=1,\ldots,k\\
\vk l\notin\{-1,0,1\}^k
}}
q_{\vk l}(u),
\]
and
\[
\Sigma_{\mathrm{nb}}
=
\sum_{\vk l\in\{-1,0,1\}^k\setminus\{\vk0\}}
q_{\vk l}(u).
\]
Let
$
a
\coloneq
\vk b^\top\Sigma^{-1}\vk b
=
\vk w^\top\vk b
>0$.
Since $\vk w>\vk0$, for all $u>0$ we have 
\[
q_{\vk l}(u)
\le
\pk*{
\sup_{(\vk s,\vk t)\in\Delta_{\vk l}\times\Delta_{\vk0}}
\vk w^\top
\frac{\vk X(\vk s)+\vk X(\vk t)}2
>
ua
}.
\]
For \(\vk h=\vk s-\vk t\) we have 
\[
\Var\left(
\vk w^\top
\frac{\vk X(\vk s)+\vk X(\vk t)}2
\right)
=
a-\frac12\vk w^\top(\Sigma-\RR(\vk h))\vk w .
\]
\qq{{\it The sum $\Sigma_{\mathrm{far}}$.}}
For indices in \(\Sigma_{\mathrm{far}}\), all sufficiently large \(u\), and
all \((\vk s,\vk t)\in\Delta_{\vk l}\times\Delta_{\vk0}\), one has
$ \|\vk s-\vk t\|_\infty\ge \theta/2 $.
If \(T_* = \limsup_{u\to\infty} T_u<\infty\), then for fixed \(S\) and all
large \(u\), every index in \eqref{eq:double-target} satisfies
\[
  |s_i-t_i|
  \le
  \frac{S(|l_i|+1)}{v_i(u)}
  \le
  T_u+\frac{S}{v_i(u)}
  \le
  T_*+1,
  \qquad i=1,\ldots,k.
\]
Thus the same lags also satisfy
\(\|\vk s-\vk t\|_\infty\le T_*+1\).
Therefore, by Assumption~\ref{I:B1} and
by the continuity of \(\RR\),
\[
   \inf_{\theta/2\le\|\vk h\|_\infty\le T_*+1}
   \vk w^\top(\Sigma-\RR(\vk h))\vk w>0 .
\]
If \(T_*=\infty\), the same lower bound follows from
Assumption~\ref{I:B3}. Thus, in either case, there exists \(m_\theta>0\)
such that
\[
   \vk w^\top(\Sigma-\RR(\vk h))\vk w\ge m_\theta
\]
for all far-block lags under consideration. Decreasing \(m_\theta\), if
necessary, we also assume that \(m_\theta\le a\).
Consequently,
\[
\sup_{(\vk s,\vk t)\in\Delta_{\vk l}\times\Delta_{\vk0}}
\Var\left(
\vk w^\top
\frac{\vk X(\vk s)+\vk X(\vk t)}2
\right)
\le
a-\frac12m_\theta
<
\bar\sigma_\theta^2
\coloneq
a-\frac14m_\theta
<a .
\]
The expected supremum of the above centered GRF    is uniformly
bounded in \(\vk l\) and \(u\), for fixed \(S\). Indeed,
\[
\sup_{(\vk s,\vk t)\in\Delta_{\vk l}\times\Delta_{\vk0}}
\vk w^\top\frac{\vk X(\vk s)+\vk X(\vk t)}2
\le
\frac12
\sup_{\vk s\in\Delta_{\vk l}}\vk w^\top\vk X(\vk s)
+
\frac12
\sup_{\vk t\in\Delta_{\vk0}}\vk w^\top\vk X(\vk t).
\]
By stationarity, the expectations of the two suprema on the right-hand side
are equal. For fixed \(S\) and all large \(u\) each block is contained in a
fixed compact set after translation. By stationarity and Fernique's
theorem applied to $\vk w^\top\vk X$ on this compact set, the expected
suprema are finite and bounded uniformly in $\vk l$ and $u$ for fixed $S$.

By Borell--TIS, there exists \(C_S>0\) such that, uniformly over all far
indices we have 
\[
q_{\vk l}(u)
\le
\exp\left\{
-\frac{(ua-C_S)^2}{2\bar\sigma_\theta^2}
\right\}.
\]
It is well known, see, e.g., \cite{HA2005}, that
\[
p_u
=
u^{-d+o(1)}
\exp\left\{-\frac12au^2\right\}.
\]
The exponent gained relative to the Gaussian factor in \(p_u\) is
\[
\begin{aligned}
&\frac{(ua-C_S)^2}{2\bar\sigma_\theta^2}
-\frac12au^2 
=
\frac12\left(\frac{a^2}{\bar\sigma_\theta^2}-a\right)u^2
-\frac{aC_S}{\bar\sigma_\theta^2}u
+\frac{C_S^2}{2\bar\sigma_\theta^2}.
\end{aligned}
\]
Its quadratic coefficient is positive because
\(\bar\sigma_\theta^2<a\). We may therefore choose
\(\gamma_\theta>0\), independently of \(S\), such that, for every fixed
\(S\) and all sufficiently large \(u\)
\[
q_{\vk l}(u)
\le
p_u\,u^{d+o(1)}e^{-\gamma_\theta u^2}
\]
uniformly over all far indices. The number of far indices is bounded by
\(C_ST_u^k\eta(u)\). Since \(\eta(u)\) is regularly varying at infinity and
\(\ln T_u=o(u^2)\),  we have 
\[
\limit{u}T_u^k\eta(u)u^{d+o(1)}e^{-\gamma_\theta u^2}=0
\]
implying
\begin{equation}
\label{eq:far-negligible}
\lim_{S\to\infty}
\limsup_{u\to\infty}
\frac{\Sigma_{\mathrm{far}}}{S^kp_u}
=0 .
\end{equation}
\qq{{\it The sum $\Sigma_{\mathrm{nn}}$.}}
 For the near but non-neighbouring part, where
$\vk l\notin\{-1,0,1\}^k$,
recall that
\[
\rho_S(\vk l)
=
\sum_{i=1}^k
\bigl(S(|l_i|-1)_+\bigr)^{\alpha_i}.
\]
By \Cref{lem.double}, \qq{with  \(\delta\) as in \Cref{lem.double}}, there exist
\(\mathcal C,\mathcal D>0\) such that
\[
\frac{q_{\vk l}(u)}{p_u}
\le
\mathcal C
\exp\left\{
-\mathcal D\rho_S(\vk l)^{1-\delta}
\right\}
\]
uniformly over all indices in \(\Sigma_{\mathrm{nn}}\) implying 
\[
\limsup_{u\to\infty}
\frac{\Sigma_{\mathrm{nn}}}{S^kp_u}
\le
\frac{\mathcal C}{S^k}
\sum_{\vk l\in\Z^k\setminus\{-1,0,1\}^k}
\exp\left\{
-\mathcal D\rho_S(\vk l)^{1-\delta}
\right\}.
\]
Let \(p=1-\delta\in(0,1)\). Since
\[
\left(\sum_{i=1}^k a_i\right)^p
\ge
c\sum_{i=1}^k a_i^p,
\qquad a_i\ge0,
\]
we have
\[
\rho_S(\vk l)^p
\ge
c\sum_{i=1}^k
\bigl(S(|l_i|-1)_+\bigr)^{\alpha_i p}.
\]
To account explicitly for the exclusion
\(\vk l\notin\{-1,0,1\}^k\), put
\[
   a_{i,S}(n)
   :=
   \exp\left\{-c\mathcal D
   \bigl(S(|n|-1)_+\bigr)^{\alpha_i p}\right\},
\]
hence we can write 
\[
\begin{aligned}
&\sum_{\vk l\in\Z^k\setminus\{-1,0,1\}^k}
\exp\{-\mathcal D\rho_S(\vk l)^p\}                              \le
\sum_{j=1}^k
\left(\sum_{|n|\ge2}a_{j,S}(n)\right)
\prod_{i\ne j}\left(\sum_{n\in\Z}a_{i,S}(n)\right).
\end{aligned}
\]
The unrestricted one-dimensional factors are bounded uniformly for
\(S\ge1\), while
\[
   \sum_{|n|\ge2}a_{j,S}(n)
   \le C\exp\{-c'S^{\alpha_jp}\}.
\]
Thus the sum tends to zero faster than any power of \(S\) implying  
\begin{equation}
\label{eq:nn-negligible}
\lim_{S\to\infty}
\limsup_{u\to\infty}
\frac{\Sigma_{\mathrm{nn}}}{S^kp_u}
=0 .
\end{equation} 

\qq{{\it The sum $\Sigma_{\mathrm{nb}}$.}}
It remains to treat the finitely many shifts
$\vk l\in\{-1,0,1\}^k\setminus\{\vk0\}$.
Fix such a shift and choose \(j\) with \(|l_j|=1\). We treat \(l_j=1\); the
other sign reduces to this one because stationarity and interchange of the two
blocks give \(q_{\vk l}(u)=q_{-\vk l}(u)\).
Split
\[
A_{\vk l}(S)
=
A_{\vk l}^{(1)}(S)\cup A_{\vk l}^{(2)}(S),
\]
where
\[
A_{\vk l}^{(1)}(S)
=
\{\vk t\in A_{\vk l}(S):S\le t_j\le S+\sqrt S\},
\qquad 
A_{\vk l}^{(2)}(S)
=
A_{\vk l}(S)\setminus A_{\vk l}^{(1)}(S)
\] 
and hence 
\[
q_{\vk l}(u)
\le
q_{\vk l}^{(1)}(u)+q_{\vk l}^{(2)}(u),
\]
where \(q_{\vk l}^{(1)}(u)\) is bounded by the probability of an exceedance in
\(A_{\vk l}^{(1)}(S)/\vk v(u)\), and \(q_{\vk l}^{(2)}(u)\) is the joint
exceedance probability between \(A_0(S)/\vk v(u)\) and
\(A_{\vk l}^{(2)}(S)/\vk v(u)\).

The strip \(A_{\vk l}^{(1)}(S)\) has side lengths of order
$\sqrt S,\ S,\ldots,S$.
It can therefore be covered by at most \(CS^{k-1/2}\) translates of the unit
cube. By stationarity, the unit-cube \(x=0\) estimate from the single-sum step
is uniform over these translates; hence the union bound implies
\[
\limsup_{u\to\infty}
\frac{q_{\vk l}^{(1)}(u)}{p_u}
\le
CS^{k-1/2},
\]
and hence we obtain 
\[
\lim_{S\to\infty}
\limsup_{u\to\infty}
\frac{q_{\vk l}^{(1)}(u)}{S^kp_u}
=0.
\]

For the separated part, the two local sets are separated by at least
\(\sqrt S\) in the \(j\)-th coordinate. Cover $A_{\vk0}(S)$ and the closure
of $A_{\vk l}^{(2)}(S)$ by closed rectangles of side lengths at most one.
Every resulting pair \(B_1,B_2\) therefore satisfies
\[
  \inf_{\vk s\in B_1,\ \vk t\in B_2}
  [\vk t-\vk s]_{\vk\alpha}
  \ge S^{\alpha_j/2}.
\]
Moreover, for all points in the two blocks
\[
  \|(\vk t-\vk s)/\vk v(u)\|
  \le
  2S\left(\sum_{i=1}^k v_i(u)^{-2}\right)^{1/2}
  \rightarrow0 .
\]
Consequently, for fixed \(S\) and all sufficiently large \(u\), the unit-rectangle
assertion of \Cref{lem.double} applies uniformly to every covering pair.
There are at most \(C(1+S)^{2k}\) such pairs, so  with constants \(C,D>0\)
independent of \(S\)
\[
\limsup_{u\to\infty}
\frac{q_{\vk l}^{(2)}(u)}{p_u}
\le
C(1+S)^{2k}
\exp\left\{
-D S^{\alpha_j(1-\delta)/2}
\right\}
\]
yielding
\[
\lim_{S\to\infty}
\limsup_{u\to\infty}
\frac{q_{\vk l}^{(2)}(u)}{S^kp_u}
=0.
\]
Since there are only finitely many neighbouring shifts
\begin{equation}
\label{eq:nb-negligible}
\lim_{S\to\infty}
\limsup_{u\to\infty}
\frac{\Sigma_{\mathrm{nb}}}{S^kp_u}
=0 .
\end{equation}

Combining \eqref{eq:far-negligible}, \eqref{eq:nn-negligible}, and
\eqref{eq:nb-negligible} proves \eqref{eq:double-target}. Thus the double
sums in the Bonferroni bounds \qq{$I_1(u,S), I_2(u,S)$} are negligible. 
To make the order of limits
explicit  set
\[
e_S
:=
\limsup_{u\to\infty}
\frac1{S^kp_u}
\sum_{\substack{
-\vk N_u(T_u,S)\le\vk l\le\vk N_u(T_u,S)\\
\vk l\ne\vk0}}
q_{\vk l}(u).
\]
Then \eqref{eq:double-target} says that \(e_S\to0\). For every fixed \(S\)
with \(S^k>x\) and every continuity point \(x>0\) of \(\mathcal B_S\), the
Bonferroni bounds, \eqref{eq:block-number}, \eqref{eq:single-sum}, and the
already established global limit \eqref{eq:Berman-positive-x} yield
\begin{equation}
\label{eq:BS-identification-prelim}
\frac{\mathcal B_S(x)}{S^k}-e_S
\le
\mathcal B_{\vk Y,\vk w}(x)
\le
\frac{\mathcal B_S(x)}{S^k}+e_S.
\end{equation}

We now remove the exceptional-set restriction. First take \(S=n\in\N\).
The intersection over \(n\) of the continuity sets of \(\mathcal B_n\) has
full Lebesgue measure. Hence \eqref{eq:BS-identification-prelim} and
\(e_n\to0\) imply
\[
\frac{\mathcal B_n(x)}{n^k}
\rightarrow
\mathcal B_{\vk Y,\vk w}(x)
\]
on this common full-measure set. Both sides are nonincreasing in \(x\), and
\(\mathcal B_{\vk Y,\vk w}\) is continuous on \((0,\infty)\); squeezing
between nearby points of the common set extends the convergence to every
\(x>0\). Finally, if \(n\le S<n+1\), monotonicity of the integration domain
implies
\[
\mathcal B_n(x)\le \mathcal B_S(x)\le \mathcal B_{n+1}(x).
\]
After division by \(S^k\), the factors \(n^k/S^k\) and
\((n+1)^k/S^k\) both tend to one and therefore we have 
\[
\lim_{S\to\infty}
\frac{\mathcal B_S(x)}{S^k}
=
\mathcal B_{\vk Y,\vk w}(x),
\qquad x>0.
\]

Furthermore, cover \([0,S]^k\) by at most \(C(1+S)^k\) unit cubes. The local
\(x=0\) convergence and a union bound, followed by the limit in \(u\) yield  
\[
\mathcal B_S(0)
\le
C(1+S)^k .
\]
Here \(C\) is the fixed unit-cube local constant and is independent of
\(S\).
Consequently,
\[
\sup_{x>0}\mathcal B_{\vk Y,\vk w}(x)
\le
\limsup_{S\to\infty}
\frac{\mathcal B_S(0)}{S^k}
<\infty 
\]
implying thus 
\[
\mathcal B_{\vk Y,\vk w}(0)
=
\lim_{x\downarrow0}\mathcal B_{\vk Y,\vk w}(x)
=
\int_{0+}^\infty \frac1y\,d\mathfrak F_{\vk Y,\vk w}(y)
<\infty .
\]
By \Cref{contBx}, since \(\mathcal B_{\vk Y,\vk w}(0)<\infty\),
the function \(\mathcal B_{\vk Y,\vk w}\) is continuous on \([0,\infty)\).
 
It remains to justify the convergence at \(x=0\), \qq{for which we show 
the estimate}
\begin{equation}
\label{eq:small-sojourn}
\lim_{\varepsilon\downarrow0}
\limsup_{u\to\infty}
\frac{
\pk{0<L_u^*(D_u)\le\varepsilon}
}{
T_u^k\eta(u)p_u
}
=0 .
\end{equation}

Fix \(S>1\), use the inner union \(D_u^-\) and put
$
R_u(S)\coloneq D_u\setminus D_u^- $.
Then for all $u$ large
\[
\{0<L_u^*(D_u)\le\varepsilon\}
\subset
\{L_u^*(R_u(S))>0\}
\cup
\bigcup_{\vk0\le\vk l\le\vk N_u(T_u,S)-\vk1}
\{0<L_u^*(\Delta_{\vk l})\le\varepsilon\}.
\]
Indeed, if the remainder contributes nothing and the total occupation is
positive and at most \(\varepsilon\), then at least one inner block has
positive occupation, and that block's occupation is at most \(\varepsilon\).

The remainder \(R_u(S)\) is covered by at most
\[
\prod_{i=1}^k (N_{u,i}(T_u,S)+1)
-
\prod_{i=1}^k N_{u,i}(T_u,S)
\]
blocks of the form \(\Delta_{\vk l}\). Since
\(T_uv_i(u)\to\infty\) for every \(i\), we have
\[ \limit{u}
\frac{
\prod_{i=1}^k (N_{u,i}(T_u,S)+1)
-
\prod_{i=1}^k N_{u,i}(T_u,S)
}{
T_u^k\eta(u)
}
=0 .
\]
For fixed \(S\) the local \(x=0\) upper bound and stationarity imply for all
sufficiently large \(u\)
\[
\pk{L_u^*(\Delta_{\vk l})>0}\le C_Sp_u
\]
with one constant uniform in \(\vk l\) and  therefore 
\[
\limsup_{u\to\infty}
\frac{\pk{L_u^*(R_u(S))>0}}
     {T_u^k\eta(u)p_u}
=0 .
\]

By stationarity we obtain 
\[
\sum_{\vk0\le\vk l\le\vk N_u(T_u,S)-\vk1}
\pk{0<L_u^*(\Delta_{\vk l})\le\varepsilon}
=
\prod_{i=1}^k N_{u,i}(T_u,S)
\pk{0<L_u^*(\Delta_{\vk0})\le\varepsilon}.
\]
At continuity points \(\varepsilon>0\) of \(x\mapsto\mathcal B_S(x)\), the
finite-window local estimate yields
\[
\lim_{u\to\infty}
\frac{
\pk{0<L_u^*(\Delta_{\vk0})\le\varepsilon}
}{p_u}
=
\mathcal B_S(0)-\mathcal B_S(\varepsilon).
\]
Using \eqref{eq:block-number}, we get
\[
\limsup_{u\to\infty}
\frac{
\pk{0<L_u^*(D_u)\le\varepsilon}
}{
T_u^k\eta(u)p_u
}
\le
\frac{\mathcal B_S(0)-\mathcal B_S(\varepsilon)}{S^k}.
\]
For fixed \(S\) monotone convergence gives
\[
\mathcal B_S(\varepsilon)\uparrow\mathcal B_S(0),
\qquad
\varepsilon\downarrow0 .
\]
Letting \(\varepsilon\downarrow0\) through continuity points, and then using
monotonicity to remove this restriction, proves \eqref{eq:small-sojourn}.

Now, for every \(\varepsilon>0\)
\[
\pk{L_u^*(D_u)>0}
\le
\pk{L_u^*(D_u)>\varepsilon}
+
\pk{0<L_u^*(D_u)\le\varepsilon}.
\]
Using the convergence already proved for positive \(\varepsilon\), the
estimate \eqref{eq:small-sojourn}, and the continuity of
\(\mathcal B_{\vk Y,\vk w}\) at zero, we obtain
\[
\limsup_{u\to\infty}
\frac{\pk{L_u^*(D_u)>0}}
     {T_u^k\eta(u)p_u}
\le
\mathcal B_{\vk Y,\vk w}(0).
\]
The reverse inequality follows from
\[
\pk{L_u^*(D_u)>0}
\ge
\pk{L_u^*(D_u)>\varepsilon}
\]
and then letting \(\varepsilon\downarrow0\) implying 
\[
\lim_{u\to\infty}
\frac{\pk{L_u^*(D_u)>0}}
     {T_u^k\eta(u)p_u}
=
\mathcal B_{\vk Y,\vk w}(0).
\]

Combining \qq{the above} result with the convergence for \(x>0\), we obtain
for all \(x\ge0\)
\[
\lim_{u\to\infty}
\frac{
\pk{
\eta(u)
\int_{[0,T_u]^k}
\ind{\vk X(\vk t)>u\vk b}\,d\vk t
>x
}
}{
T_u^k\eta(u)\pk{\vk X(\vk0)>u\vk b}
}
=
\mathcal B_{\vk Y,\vk w}(x).
\]
It remains only to record the representation at \(x=0\). Applying the same
Bonferroni argument with \(x=0\), together with the single-sum estimate and the
double-sum negligibility already proved above, gives
\[
\lim_{S\to\infty}
\frac{\mathcal B_S(0)}{S^k}
=
\lim_{u\to\infty}
\frac{\pk{L_u^*(D_u)>0}}
     {T_u^k\eta(u)p_u}
=
\mathcal B_{\vk Y,\vk w}(0).
\]
Together with the representation already obtained for \(x>0\), this proves
\[
\mathcal B_{\vk Y,\vk w}(x)
=
\lim_{S\to\infty}
\frac1{S^k}
\int_{\R^d}
\pk*{
\int_{[0,S]^k}
\ind{
\operatorname{diag}(\vk w)
\bigl[\vk Y(\vk t)-V(\vk t)\vk w\bigr]
>
\vk z
}
\,d\vk t
>x
}
e^{\vk z^\top\vk1}\,d\vk z,
\qquad x\ge0.
\]
The finiteness follows from
$\mathcal B_S(0)\le C(1+S)^k$. We finish with a  positivity argument setting next  
\[
   Z(\vk t):=\vk w^\top\vk Y(\vk t),
   \qquad
   q_{\vk w}(\vk t):=\vk w^\top V(\vk t)\vk w.
\]
By \eqref{additional},
$q_{\vk w}(\vk t)\ge c_V[\vk t]_{\vk\alpha}$.
The homogeneity of \(V\) and the formula for \(R_V\) imply the
self-similarity
\[
   \{Z(c^{1/\vk\alpha}\vk t):\vk t\in\R^k\}
   \stackrel{\mathrm{fdd}}{=}
   \{\sqrt c\,Z(\vk t):\vk t\in\R^k\},
   \qquad c>0.
\]
Let
\[
   M=\sup_{[\vk t]_{\vk\alpha}\le1}|Z(\vk t)|,
   \qquad
   M_n=\sup_{[\vk t]_{\vk\alpha}\le2^n}|Z(\vk t)|.
\]
Continuity of $Z$, compactness of $\{[\vk t]_{\vk\alpha}\le1\}$, and Fernique's theorem imply
that $M$ is finite and integrable. Borell--TIS gives a Gaussian tail for $M$, while
self-similarity gives \(M_n\stackrel d=2^{n/2}M\). Consequently, for every
\(\varepsilon>0\)
\[
   \sum_{n=1}^\infty
   \pk{M_n>\varepsilon2^n}
   <\infty.
\]
The Borel--Cantelli lemma, applied to a sequence of positive rational
$\varepsilon$'s, gives $M_n\le\varepsilon2^n$ eventually for every such
$\varepsilon$. If $2^{n-1}<[\vk t]_{\vk\alpha}\le2^n$, then
\[
   \frac{|Z(\vk t)|}{[\vk t]_{\vk\alpha}}
   \le
   2\frac{M_n}{2^n}.
\]
It follows that
\begin{equation}
\label{eq:scalar-tangent-growth}
   Z(\vk t)=o([\vk t]_{\vk\alpha})
   \quad\text{almost surely as }[\vk t]_{\vk\alpha}\to\infty.
\end{equation}

If the indicator defining \(\mathcal J(\vk Y,\vk w)\) equals one at
\(\vk t\), summing its componentwise inequalities gives
\[
   Z(\vk t)-q_{\vk w}(\vk t)+\vk1^\top\vk E>0.
\]
Since \(\vk1^\top\vk E<\infty\) almost surely,
\eqref{eq:scalar-tangent-growth} and
$q_{\vk w}(\vk t)\ge c_V[\vk t]_{\vk\alpha}$ show that
the occupied set is almost surely bounded. Hence
\[
   \mathcal J(\vk Y,\vk w)<\infty
   \qquad\text{almost surely}.
\]
It is also strictly positive almost surely: at \(\vk t=\vk0\), the vector
inside the indicator equals \(\vk E>\vk0\), and continuity makes the same
inequality hold on a random neighbourhood of the origin.

Finally, \(\mathcal J(\vk Y,\vk w)\) has unbounded support. Given \(M_0>0\),
choose a compact rectangle \(K_0\subset\R^k\) with
\(\lambda_k(K_0)>M_0\), and, conditional on \(\vk Y\), set
\[
   A_j
   :=
   1+
   \max\left(0,
   -\min_{\vk t\in K_0}
   w_j\bigl[Y_j(\vk t)-(V(\vk t)\vk w)_j\bigr]
   \right),
   \qquad j=1,\ldots,d.
\]
All $A_j$ are finite by continuity, and independence of the exponential
variables gives
\[
   \pk{E_j>A_j,\ 1\le j\le d\mid\vk Y}
   =
   \exp\left\{-\sum_{j=1}^dA_j\right\}
   >0.
\]
On this event the occupied set contains \(K_0\), so
$\mathcal J(\vk Y,\vk w)>M_0$. Taking expectations gives
\[
   \pk{\mathcal J(\vk Y,\vk w)>M_0}
   \ge
   \E*{\exp\left\{-\sum_{j=1}^dA_j\right\}}
   >0.
\]
Since $\mathcal J(\vk Y,\vk w)<\infty$ almost surely, it follows that
\[
   \pk{M_0<\mathcal J(\vk Y,\vk w)<\infty}>0,
   \qquad M_0>0.
\]
The Stieltjes representation \eqref{B.def.F} now shows that
\(\mathcal B_{\vk Y,\vk w}(x)>0\) for every \(x\ge0\). This completes the
proof.
\QED

\subsection{Compound Poisson extension}
\label{s:compound-Poisson}
\begin{proof}[Proof of the assertion in \Cref{rem:cp-extension}(iii)]
Write
\[
 p_u=\pk{\vk X(\vk0)>u\vk b},\qquad
 q=\vk b^\top\Sigma^{-1}\vk b>0,\qquad
 \mathcal B(x)=\mathcal B_{\vk Y,\vk w}(x).
\]
The Savage condition and the Gaussian density give
\begin{equation}\label{eq:cp-point-tail}
 p_u\sim c_0u^{-d}e^{-qu^2/2},\qquad
 c_0=(2\pi)^{-d/2}|\Sigma|^{-1/2}\prod_{i=1}^d w_i^{-1}.
\end{equation} 

We first establish the comparison estimate needed for blocking. Put
\[
 r(a)=\sup_{\|\vk h\|_\infty\ge a}
 \|\Sigma^{-1/2}\RR(\vk h)\Sigma^{-1/2}\|.
\]
Let $D$ be a finite union of cubes $D_1,\ldots,D_N$, with distinct cubes
separated by at least $a$ in the sup norm, and set
\[
 S_j=\eta(u)\int_{D_j}\ind{\vk X(\vk t)>u\vk b}\,d\vk t.
\]
For each fixed $s>0$, all sufficiently large $u$, and all sufficiently small
$r(a)$, we claim that
\begin{equation}\label{eq:cp-laplace-comparison}
 \left|\E{e^{-s\sum_{j=1}^N S_j}}
       -\prod_{j=1}^N\E{e^{-sS_j}}\right|
 \le C_s\bigl(\eta(u)\lambda_k(D)p_u\bigr)^2
       u^2r(a)e^{Cu^2r(a)}.
\end{equation}
The constants do not depend on the number or sizes of the cubes.

To prove \eqref{eq:cp-laplace-comparison}, first approximate each occupation
integral by a finite Riemann sum with cell volumes $a_l>0$. Its Laplace
functional is a product of factors
\[
 1-c_l\ind{\vk X(\vk t_l)>u\vk b},\qquad
 c_l=1-e^{-s\eta(u)a_l}\le s\eta(u)a_l.
\]
Interpolate the covariance matrix of all grid values between independent
copies on the different cubes and the original field. Only covariances
between different cubes change. The Gaussian interpolation identity bounds
the contribution from sites $l,m$ and coordinates $i,j$ by
$c_lc_m|\RR_{ij}(\vk t_l-\vk t_m)|$ times a two-site Gaussian boundary
integral. The remaining factors lie in $[0,1]$ and may be discarded in an
upper bound.

Here is a uniform bound for that boundary integral. Its covariance matrix is
\[
 C_\tau(\vk h)=
 \begin{pmatrix}\Sigma&\tau\RR(\vk h)\\
                 \tau\RR(\vk h)^\top&\Sigma\end{pmatrix},
 \qquad 0\le\tau\le1,
 \qquad \vk B=(\vk b^\top,\vk b^\top)^\top.
\]
Uniformly in $\tau$ and $\|\vk h\|_\infty\ge a$, for small $r(a)$,
\[
 C_\tau(\vk h)^{-1}\vk B\ge c\vk1,\qquad
 \vk B^\top C_\tau(\vk h)^{-1}\vk B\ge2q-Cr(a),
\]
where $c>0$, and the determinants are bounded away from zero. These facts
follow by perturbing $\operatorname{diag}(\Sigma,\Sigma)$ and using
$\vk w>\vk0$. In the Gaussian density fix coordinates $i,d+j$ at
$ub_i,ub_j$ and integrate the other coordinates above their thresholds.
Writing the integration variable as $u\vk B+\vk z$, where
$z_i=z_{d+j}=0$ and the other $z_l\ge0$, yields the bound
\[
 C e^{-qu^2+Cr(a)u^2}
 \int_{[0,\infty)^{2d-2}}e^{-cu\sum z_l}\,d\vk z
 \le C u^{-(2d-2)}e^{-qu^2+Cr(a)u^2}
 \le C u^2p_u^2e^{Cr(a)u^2}.
\]
For $d=1$ the integral over zero free coordinates equals one. We used
\eqref{eq:cp-point-tail} in the last inequality.

Since $|\RR_{ij}(\vk h)|\le Cr(a)$ between distinct cubes, summing the
interpolation bound and using $\sum_l a_l=\lambda_k(D)$ proves
\eqref{eq:cp-laplace-comparison} for the Riemann sums. The interpolation
identity can first be applied to smooth approximations of the factors and
to covariance matrices made nonsingular by adding arbitrarily small
independent Gaussian noises. For fixed $u$, take the noise standard
deviations and the smoothing widths at most $1/u$, and then let both tend
to zero. The same density estimate holds with a bounded extra factor before
these limits. The mixed derivatives between different sites contain the two
boundary densities and the other $2d-2$ threshold constraints, so the bounds
above are uniform in the full grid size.
No inverse of the covariance matrix of the entire grid is used in the
bound. Finally, sample continuity and
\[
 \E{\lambda_k\{\vk t\in D:X_i(\vk t)=ub_i\}}=0,
 \qquad i=1,\ldots,d,
\]
show that the excursion indicators are Riemann integrable almost surely.
For each fixed $u$ let the mesh tend to zero. Bounded convergence and the
mesh-independent bound prove \eqref{eq:cp-laplace-comparison} for the
occupation integrals.

We now choose the cubes. For $A>1$ put
\[
 \epsilon(A)=\sup_{\|\vk h\|_\infty\ge A}
 \|\Sigma^{-1/2}\RR(\vk h)\Sigma^{-1/2}\|
 \log\|\vk h\|_\infty.
\]
The assumed Berman condition gives $\epsilon(A)\to0$. Set, for large $u$,
\[
 \delta_u=\max\{u^{-1},\sqrt{\epsilon(e^u)}\},\qquad
 a_u=e^{\delta_u u^2},\qquad
 \ell_u=e^{2\delta_u u^2}.
\]
Then $\delta_u\to0$, $a_u/\ell_u\to0$, $\log\ell_u=o(u^2)$, and
\[
 u^2r(a_u)\le\epsilon(e^u)/\delta_u\to0.
\]
The intensity assumption and \eqref{eq:cp-point-tail} imply
$\log T_u/u^2\to q/(2k)>0$, so $\ell_u/T_u\to0$.
Retain the Cartesian array of cubes of side $\ell_u$ with spacing
$\ell_u+a_u$ contained in $[0,T_u]^k$. Their number $N_u$ satisfies
\[
 N_u=\left\lfloor\frac{T_u}{\ell_u+a_u}\right\rfloor^k,
 \qquad
 \frac{N_u\ell_u^k}{T_u^k}\to1.
\]
Writing $D_u^{\rm ret}$ for their union, stationarity gives
\[
 \E*{\eta(u)\int_{[0,T_u]^k\setminus D_u^{\rm ret}}
       \ind{\vk X(\vk t)>u\vk b}\,d\vk t}
 =\eta(u)p_u\bigl(T_u^k-N_u\ell_u^k\bigr)\to0.
\]
Thus discarding the corridors and boundary remainder changes the total
scaled occupation by a quantity tending to zero in $L^1$.

Let $S_u$ be the scaled occupation in a single retained cube and put
$b_u=\ell_u^k\eta(u)p_u$. Theorem~\ref{T:doublesum} applied with
$T=\ell_u$ implies
\[
 b_u\to0,\qquad N_ub_u\to\vartheta,\qquad
 \frac{\pk{S_u>x}}{b_u}\to\mathcal B(x),\qquad x\ge0.
\]
For every fixed $s>0$, using the convergence at zero as a dominating bound,
\[
 \frac{1-\E{e^{-sS_u}}}{b_u}
 =s\int_0^\infty e^{-sx}\frac{\pk{S_u>x}}{b_u}\,dx
 \longrightarrow s\int_0^\infty e^{-sx}\mathcal B(x)\,dx.
\]
Consequently the product of the one-block Laplace transforms converges to
\[
 \exp\left\{-\vartheta s\int_0^\infty e^{-sx}\mathcal B(x)\,dx\right\}.
\]
Since $\eta(u)\lambda_k(D_u^{\rm ret})p_u\to\vartheta$ and
$u^2r(a_u)\to0$, \eqref{eq:cp-laplace-comparison} transfers this limit to the
dependent retained cubes. The discarded occupation tends to zero in $L^1$,
so the same limit holds on $[0,T_u]^k$. Finally, $\mathcal B$ is continuous,
$\mathcal B(0)\in(0,\infty)$, and $\mathcal B(x)\le1/x$ for $x>0$.
Thus $\mathcal B(x)/\mathcal B(0)$ is the survival function of a positive
jump, and the displayed Laplace transform is that of a compound Poisson law
with intensity $\vartheta\mathcal B(0)$ and the claimed jump distribution.
\end{proof}

\def\dist{\mathrm{dist}}

\section{Appendix}

Throughout the appendix, inequalities between vectors are understood
componentwise.

Denote by \(\Pi_{\Sigma}(\bm b)\) the quadratic programming problem
\[
  \Pi_{\Sigma}(\bm b):
  \qquad
  \text{minimise } \bm x^\top\Sigma^{-1}\bm x
  \quad\text{over } \bm x\in\R^d \text{ such that } \bm x\ge \bm b ,
\]
where \(\Sigma\) is a \(d\times d\) symmetric positive definite matrix and
\(\bm b\in\R^d\setminus(-\infty,0]^d\). For a matrix \(A=(a_{ij})\), and
non-empty index sets \(I,J\subset\{1,\dots,d\}\), we write \(A_{IJ}\) for
the submatrix of \(A\) obtained by keeping the rows indexed by \(I\) and the
columns indexed by \(J\).

\begin{lemma}\label{AL}
Let \(d\ge1\), let \(\Sigma\) be a \(d\times d\) symmetric positive definite
matrix, and let \(\bm b\in\R^d\setminus(-\infty,0]^d\). Then
\(\Pi_{\Sigma}(\bm b)\) has a unique solution \(\widetilde{\bm b}\).
Moreover, there exists a unique non-empty determining index set
\(I\subset\{1,\dots,d\}\), with \(m:=|I|\), such that, writing
\(J:=\{1,\dots,d\}\setminus I\),
\begin{align}
  \widetilde{\bm b}^{\top}\Sigma^{-1}\bm z
  &=
  \bm b_I^\top\Sigma_{II}^{-1}\bm z_I,
  \qquad \bm z\in\R^d,                                             \label{Cc}\\
  \widetilde{\bm b}_I
  &=
  \bm b_I\neq\bm0_I,                                                \label{eq:tildeb-active}\\
  \widetilde{\bm b}_J
  &=
  \Sigma_{JI}\Sigma_{II}^{-1}\bm b_I\ge \bm b_J,
  \qquad\text{if }J\neq\emptyset,                                   \label{empt}\\
  \bm w
  &:=
  \Sigma^{-1}\widetilde{\bm b},
  \qquad
  \bm w_I=\Sigma_{II}^{-1}\bm b_I>\bm0_I,
  \qquad
  \bm w_J=\bm0_J,                                                    \label{emptB}\\
  \tau_{\bm b}
  &:=
  \min_{\bm x\ge\bm b}\bm x^\top\Sigma^{-1}\bm x
  =
  \widetilde{\bm b}^{\top}\Sigma^{-1}\widetilde{\bm b}
  =
  \bm b_I^\top\Sigma_{II}^{-1}\bm b_I
  =
  \bm w^\top\widetilde{\bm b}
  =
  \bm w_I^\top\bm b_I>0,                                             \label{eq:alfa}.
\end{align}
Finally, \(m=d\) if and only if \(\Sigma^{-1}\bm b>\bm0\).
\end{lemma}

\begin{proof}
For \(d=1\), the assertion is immediate. For \(d\ge2\), all claims are
stated in \cite[Lem 4.1]{DHW20}; see also \cite[Prop 2.1]{HA2005}.
\end{proof}

We now state the vector-valued Pickands lemma for sojourn functionals. Let
\(E\subset\R^k\) be compact with \(\lambda_k(E)>0\), and fix
\(\bm t_0\in E\). Let \(Q_u\), \(u>0\), be non-empty index sets. Let
\[
  \bm Z_{u,\tau}(\bm t),
  \qquad \bm t\in E,\quad u>0,\quad \tau\in Q_u
\]
be centered \(\R^d\)-valued GRFs with almost surely continuous sample paths. We write
\[
  R_{u,\tau}(\bm s,\bm t)
  :=
  \Cov\bigl(\bm Z_{u,\tau}(\bm s),\bm Z_{u,\tau}(\bm t)\bigr),
  \qquad \bm s,\bm t\in E.
\]
Let \(\bm Z(\bm t)\), \(\bm t\in E\), be a centered \(\R^d\)-valued GRF with
almost surely continuous sample paths and the covariance matrix function
\(K(\bm s,\bm t)\).
Let \(\bm b\in\R^d\setminus(-\infty,0]^d\). The following assumptions are
imposed.

\begin{enumerate}[({C}1)]
\item\label{C1}
For all large \(u\) and all \(\tau\in Q_u\), the matrix
\[
  \Sigma_{u,\tau}:=R_{u,\tau}(\bm t_0,\bm t_0)
\]
is positive definite and
\[
  \lim_{u\to\infty}
  \sup_{\tau\in Q_u}
  u\left\|\Sigma-\Sigma_{u,\tau}\right\|
  =0
\]
for some positive definite matrix \(\Sigma\).
\end{enumerate}

Let \(\widetilde{\bm b}\) be the unique solution of
\(\Pi_\Sigma(\bm b)\), and put
\[
  \bm w:=\Sigma^{-1}\widetilde{\bm b}.
\]
Let \(I\subset\{1,\dots,d\}\) be the determining index set from
Lemma~\ref{AL}, let \(J:=\{1,\dots,d\}\setminus I\), and set
$m:=|I|$.
\begin{enumerate}[({C}2)]
\item\label{C2}
We have
\[
  \lim_{u\to\infty}
  \sup_{\tau\in Q_u,\ \bm t\in E}
  u\left\|
  \Sigma_{u,\tau}-R_{u,\tau}(\bm t,\bm t_0)
  \right\|
  =0 .
\]
Moreover, there exists a continuous function \(\bm f:E\to\R^d\) such that
\[
  \lim_{u\to\infty}
  \sup_{\tau\in Q_u,\ \bm t\in E}
  \left\|
  u^2
  \left[
  \mathcal I_d-R_{u,\tau}(\bm t,\bm t_0)\Sigma_{u,\tau}^{-1}
  \right]
  \widetilde{\bm b}
  -
  \bm f(\bm t)
  \right\|
  =0
\]
and
\[
  \lim_{u\to\infty}
  \sup_{\tau\in Q_u}
  \sup_{\bm s,\bm t\in E}
  \left\|
  u^2
  \left[
  R_{u,\tau}(\bm s,\bm t)
  -
  R_{u,\tau}(\bm s,\bm t_0)
  \Sigma_{u,\tau}^{-1}
  R_{u,\tau}(\bm t_0,\bm t)
  \right]
  -
  K(\bm s,\bm t)
  \right\|
  =0 .
\]
\end{enumerate}

\begin{enumerate}[({C}3)]
\item\label{C3}
There exist constants \(C>0\) and
\(\bm\gamma=(\gamma_1,\dots,\gamma_k)\in(0,2]^k\) such that, for all large
\(u\)
\[
  \sup_{\tau\in Q_u}
  u^2
  \E*{
  \left\|
  \bm Z_{u,\tau}(\bm t)-\bm Z_{u,\tau}(\bm s)
  \right\|^2
  }
  \le
  C\sum_{\ell=1}^k |t_\ell-s_\ell|^{\gamma_\ell},
  \qquad \bm s,\bm t\in E .
\]
\end{enumerate}

For \(x\ge0\), define
\[
  \mathcal B_{\bm Z,\bm w}(x;E)
  :=
  \int_{\R^d}
  e^{\bm1^\top\bm y}
  \pk*{
  \int_E
  \ind{
  \operatorname{diag}(\bm w)
  \bigl[\bm Z(\bm t)-\bm f(\bm t)\bigr]
  >
  \bm y
  }
  d\bm t
  >x
  }
  d\bm y .
\]
In the notation of the main text, if \(\bm f(\bm t)=V(\bm t)\bm w\), this
coincides with the definition involving \(V(\bm t)\bm w\).

\begin{lemma}\label{lemma:Pickands}
In the setup above, assume that \textup{\ref{C1}--\ref{C3}} hold. Let
\(\bm X\sim N(\bm0,\Sigma)\). Then the following convergence holds for
\(x=0\) and for every continuity point \(x\in(0,\lambda_k(E))\) of
\(x\mapsto\mathcal B_{\bm Z,\bm w}(x;E)\)
\[
\lim_{u\to\infty}
\sup_{\tau\in Q_u}
\left|
\frac{
\pk*{
\displaystyle
\int_E
\ind{\bm Z_{u,\tau}(\bm t)>u\bm b}
\,d\bm t>x
}
}{
\pk{\bm X>u\bm b}
\exp\!\left(
-\frac12u^2
\widetilde{\bm b}^{\top}
(\Sigma_{u,\tau}^{-1}-\Sigma^{-1})
\widetilde{\bm b}
\right)
}
-
\mathcal B_{\bm Z,\bm w}(x;E)
\right|
=0 .
\]
\end{lemma}

\prooflem{lemma:Pickands}
The proof follows by combining the uniform conditioning argument in
\cite[Lem.~4.7]{DHW20} with the sojourn argument in
\cite[Lem.~4.1]{DHLM23}, using further \Cref{lemma:integral_estimate2}.
The extension to every continuity point follows by monotonicity.
\QED

\section*{Acknowledgments}
K. D{\k{e}}bicki 
was partially supported by the National Science Centre, Poland, Grant No
2024/55/B/\allowbreak ST1/01062
(2025-2028). 

% Bibliography generated from the unchanged vecPPPP.bib with ieeetr.
% Embedded here so this revision is a standalone LaTeX document.
\def\polhk#1{\setbox0=\hbox{#1}{\ooalign{\hidewidth
  \lower1.5ex\hbox{`}\hidewidth\crcr\unhbox0}}}

\end{document}